\documentclass[12pt,a4paper, twoside]{article}
\usepackage{amssymb,amsfonts,amsmath,amsthm,euscript}
\usepackage[latin1]{inputenc}
\usepackage[english]{babel}
\usepackage[left=2.7cm,right=2.7cm,bottom=2.66cm,top=2.66cm]{geometry}
\usepackage{dsfont}
\usepackage{mathrsfs}
\usepackage{natbib}
\usepackage{enumerate}
\usepackage{color}
\usepackage{sectsty}
\usepackage{setspace}
\usepackage{graphicx}
\usepackage{subfigure}
\usepackage{placeins}
\usepackage{tikz}
\usepackage[final]{pdfpages}

\graphicspath{{figures/}}

\theoremstyle{plain}
\newtheorem{thm}{Theorem}[section]
\newtheorem{prop}{Proposition}[section]

\newtheorem{coro}{Corollary}[section]

\theoremstyle{definition}
\newtheorem{rmk}{Remark}[section]

\newcommand{\E}{\mathds E}
\newcommand{\var}{\mathrm{Var}}
\newcommand{\R}{\mathds R}
\newcommand{\N}{\mathds N}

\newcommand{\F}{\mathscr{F}}

\newcommand{\bs}[1]{\boldsymbol{#1}}
\renewcommand{\qed}{\hfill \mbox{\raggedright \rule{.07in}{.1in}}}
\renewcommand{\proof}{\noindent\textbf{Proof: }}

\long\def\sfootnote[#1]#2{\begingroup%
\def\thefootnote{\fnsymbol{footnote}}\footnote[#1]{#2}\endgroup}

\def\bfootnote{\xdef\@thefnmark{}\@footnotetext}

\begin{document}
\pagestyle{myheadings} 
\markboth{ Chaotic driven beta regression models }{G. Pumi, T.S. Prass and R.R. Souza}

\thispagestyle{empty}
{\centering
\Large{\bf A Dynamic Model for Double Bounded Time Series With Chaotic Driven Conditional Averages}\vspace{.5cm}\\
\large{ {\bf Guilherme Pumi$\!\phantom{i}^{\mathrm{a,}}$\sfootnote[1]{Corresponding author. This Version: \today}\let\thefootnote\relax\footnote{\hskip-.3cm$\phantom{s}^\mathrm{a}$Mathematics and Statistics Institute - Federal University of Rio Grande do Sul
}, Taiane Schaedler Prass$\!\!\phantom{s}^\mathrm{a}$ and  Rafael Rig\~ao Souza$\!\!\phantom{s}^\mathrm{a}$
}
 \\
\let\thefootnote\relax\footnote{E-mail: guilherme.pumi@ufrgs.br (G. Pumi), taiane.prass@ufrgs.br (T.S. Prass), rafars@mat.ufrgs.br (R.R. Souza)
}\\
\vskip.3cm
}}
\begin{abstract}
In this work we introduce a class of dynamic models for time series taking values on the unit interval. The proposed model follows a generalized linear model approach where the random component, conditioned on the past information,  follows a beta distribution, while the conditional mean specification may include covariates and also an extra additive term given by the iteration of a map that can present chaotic behavior.
The resulting  model is very flexible and its systematic component can accommodate short and long range dependence, periodic behavior, laminar phases, etc. We derive easily verifiable conditions for the stationarity of the proposed model, as well as conditions for the law of large numbers and a Birkhoff-type theorem to hold. A Monte Carlo simulation study is performed to assess the finite sample behavior of the partial maximum likelihood approach for parameter estimation in the proposed model. Finally, an application  to the proportion of stored hydroelectrical energy in Southern Brazil is presented. \vspace{.2cm}\\
\noindent \textbf{Keywords:} time series; chaotic processes; generalized linear models.\vspace{.2cm}\\
\noindent \emph{2010 Mathematical Subject Classification}: Primary: 37M10, 62M10, 62J12.
\end{abstract}
\onehalfspacing

\section{Introduction}

Many time series encountered in statistical applications present two important characteristics: bounds, in the sense that its distribution has a bounded support, and serial dependence. Common cases are  rates and proportions observed over time. In these cases, Gaussian based approaches are not adequate. Time series modeling of double bounded time series has been subject of intense research, especially in the last decade, and several approaches to the problem have been proposed. One such approach has received a lot of attention in the last few years. The idea is to include a time dependent structure into a Generalized Linear Model (GLM) framework and has been popularized in the works of \cite{Zeger1988}, \cite{Benjamin2003} and  \cite{Ferrari2004}. Processes following this type of structure are often referred to as Generalized Autoregressive Moving Average (GARMA) models. More specifically, the model's systematic component follows the usual approach of GLM with an additional dynamic term of the form
\begin{align}\label{e:add_spec}
	g(\mu_t)={\eta}_t=\bs{x}_{t}'\bs{\beta}+\tau_t,
\end{align}
where $g$ is a suitable link function, $\mu_t$ is some quantity of interest (usually the (un)con\-di\-tio\-nal mean or median), $\bs{x}_t$ denotes a vector of {(possibly random)} covariates observed at time $t$ and $\tau_t$ is a term responsible to accommodate  any serial correlation in the sequence $\mu_t$. The term $\tau_t$ can take a variety of forms, depending on the model's scope and intended application. In the Beta Autoregressive Moving Average ($\beta$ARMA) model \citep{Rocha2009}, for instance, the model's random component follows a (conditional) beta distribution while in the specification for the conditional mean $\mu_t$, $\tau_t$ follows a classical ARMA process. In \cite{Bayers}, the authors define the Kumaraswamy ARMA model (KARMA), where the model's random component follows a (conditional) Kumaraswamy distribution while in the specification for the conditional median $\mu_t$, $\tau_t$ also follows an ARMA process. Since ARMA models can only accommodate short range dependence, these models can only account for a short range dependence structure on their systematic component. In the case of conditionally beta distributed random component, \cite{Pumi2019} introduce the $\beta$ARFIMA (Beta Autoregressive Fractionally Integrated Moving Average) model generalizing \cite{Rocha2009} by allowing $\tau_t$ to follow a long range dependent ARFIMA  process \cite[see, for instance,][]{Hosking1981,Brockwell1991,palma2007,Box2008}.  Inference for this type of models is done via partial maximum likelihood.

In this work we propose a model where the random component follows a conditional beta distribution, while the systematic component depends on the iterations of a (usually chaotic) map defined on the unit interval.
Let $T:[0,1]\rightarrow[0,1]$ and $U_0$ be a random variable taking values in $(0,1)$, in particular, provided the existence of an absolute continuous invariant measure for $T$, $U_0$ will be distributed according to it {(see Section \ref{ssd})}.  We consider the so-called class of chaotic process defined by setting $Z_t:=h\big(T^t(U_0)\big)$, $t\in\N$, for a suitably smooth link function $h:(0,1)\rightarrow\R$. {A key concept here is the one of invariant measure for a map $T$:
	provided an absolute continuous invariant measure exists, $U_0$ will be distributed according to it {(see Section \ref{ssd})}.} Observe that $T$ does not need to be a chaotic transformation in the usual sense (see next section) for $Z_t$ to be called a chaotic process.
Such processes have been applied in a variety of problems from rock drilling \citep[see][and references therein]{ly} to intermittency in human cardiac rate \citep[see][]{Zebrowsky2001}, econometrics \citep{Gandolfo}, biology and medicine \citep{Jackson}, etc.  However, the goal on these applications usually lie on understanding the dynamics of the process (i.e., intermittence, presence of fixed/attracting/repelling points, {invariant measures}, etc)  rather than statistical inference or forecasting.
\begin{figure}[h!]
	\centering
	\mbox{
		\subfigure[]{\includegraphics[width=0.5\textwidth]{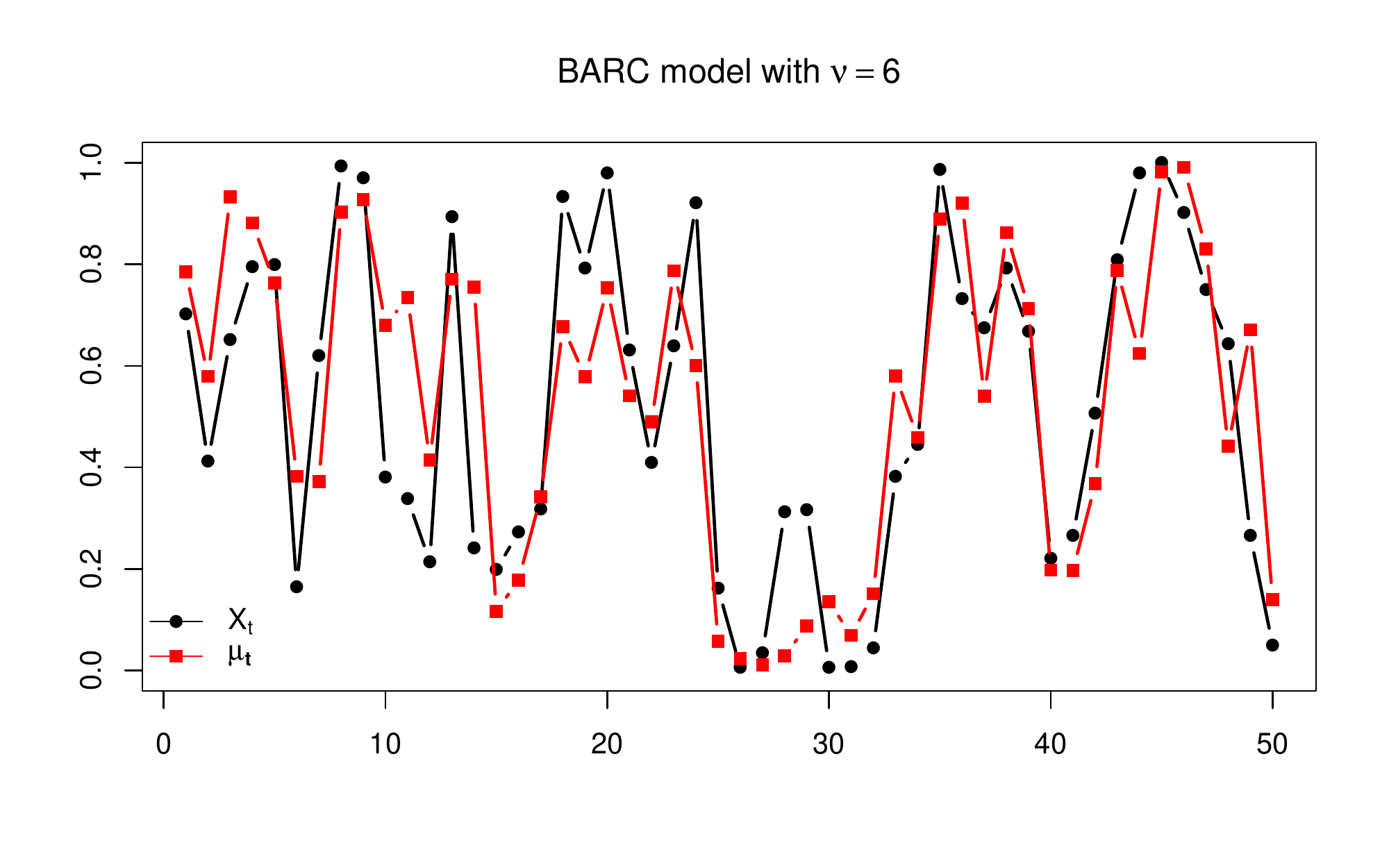}}
		\subfigure[]{\includegraphics[width=0.5\textwidth]{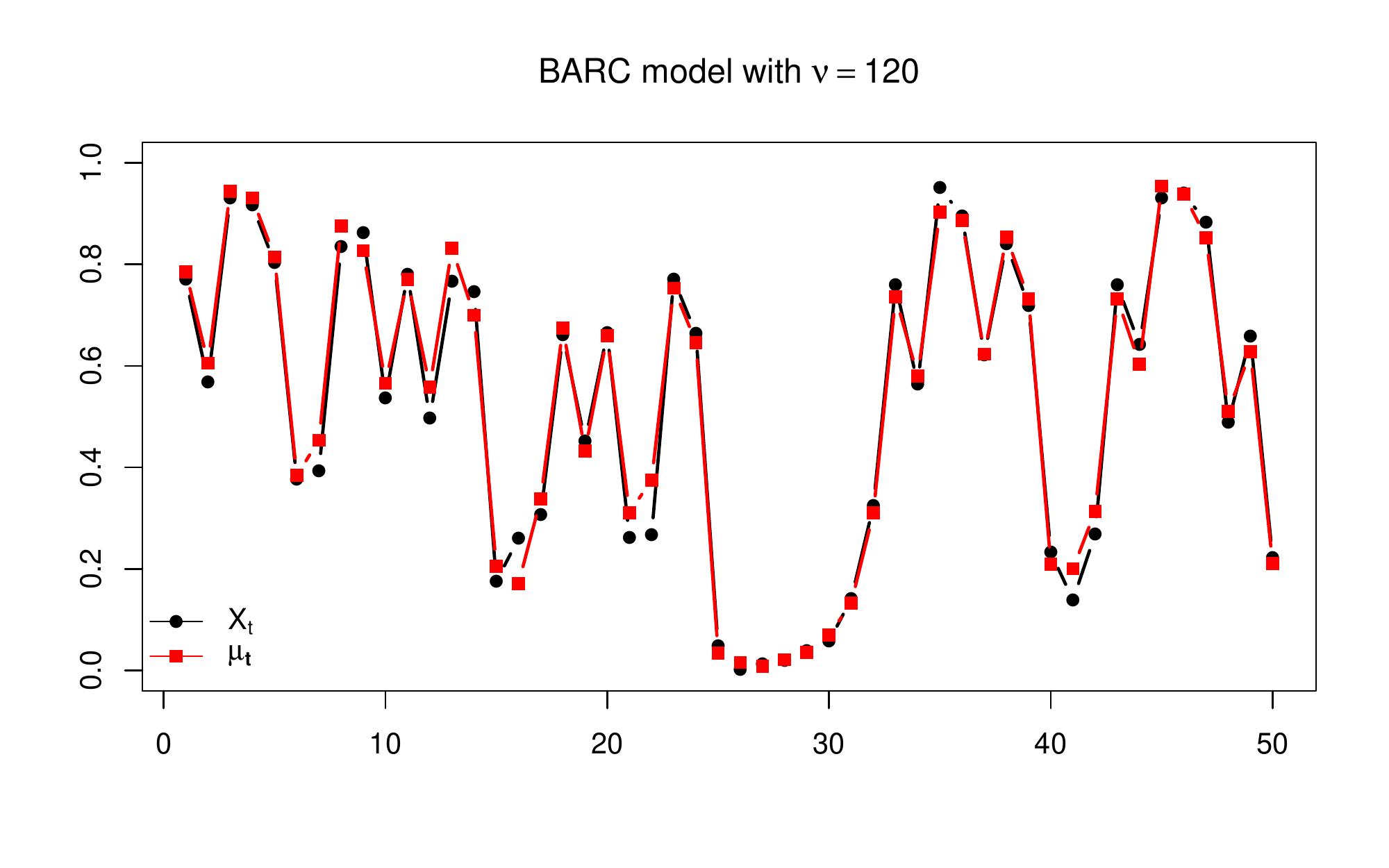}}
	}
	\caption{
		In (a) and (b) we present the sample paths (black) and conditional mean $\mu_t$ (red) of two $\beta$ARC(1) models obtained from the Mannevile-Pomeau transformation starting at $u_0=\pi/4$ with parameters parameter $\phi=0.3$ and $s=0.3$. In (a) we have $\nu=6$ while in (b) $\nu=120$.}\label{f:chaos}
\end{figure}

The novelty of the proposed model lies in 3 different fronts: first, its capability of modeling non-linear behaviors that other GARMA-like models can't; second, its flexibility; and finally, general theoretical results that are not available for other GARMA-like models in the literature can be obtained for $\beta$ARC models under easily verifiable conditions.

Upon changing the transformation $T$, one can drastically change the sample paths properties and dependence behavior of the resulting $\beta$ARC model.
Possibilities include intrinsic periodical behavior (generated by  repelling or absorbing periodic points in the dynamics), laminar phases, histogram control, and many other non-linear behavior, which cannot be mimicked by classical ARMA and ARFIMA structures present in the standard GARMA-like models, such as the $\beta$ARMA/KARMA/$\beta$ARFIMA.  These structures can be obtained simply by changing the transformation $T$, which translates into a very general and flexible class of models capable of modeling a wide variety of non-linear behavior in the systematic component. It also means that we can effective forecast more general dependence structures, especially non-linear ones. Finally, despite its flexibility, the proposed model also allow for the derivation of several general mathematical results absent in the GARMA-like model literature. For instance, easily verifiable conditions for its stationarity are available and its unconditional covariance structure is also obtainable. Furthermore, under very mild conditions, a strong law of large numbers and a Birkhoff-type theorem hold. To the best of our knowledge, similar results are not yet available in this generality for other GARMA-like processes in the literature.

\newpage

The paper is organized as follows. In the next section we define the proposed model and present some basic results from dynamical systems necessary to the work. We also present a miscellany of theoretical results regarding stationarity, law of large numbers and the covariance structure of the proposed model.	In Section \ref{inf} we consider inference on the proposed model via the partial maximum likelihood (PMLE) approach. In Section \ref{sim} we briefly present a Monte Carlo simulation study to assess the finite sample performance of the PMLE approach.  The usefulness of the proposed model is illustrated through an application to real data regarding the proportion of stored hydroelectrical energy in southern Brazil (Section \ref{se:application}). Conclusions are reserved to Section \ref{conc}. This paper is also accompanied by a supplementary material in which we present more details regarding dynamical systems and also a broad Monte Carlo simulation study to assess the finite sample performance of the PMLE approach.

\section{Model Definition and Properties}\label{stat}

In this section we shall define the proposed model and prove a miscellany of theoretical results related to it. We also present some basic definitions from  dynamical systems necessary to the work.
\subsection{Model Definition}

Let $\{Y_t\}_{t\geq1}$ be a time series of interest and let $\{\bs x_t\}_{t\geq 1}$ denote a set of $l$-dimensional exogenous time dependent (possibly random) covariates. Let $\F_{t}$ denote the $\sigma$-field representing the observed history of the model up to time $t$, that is, the sigma-field generated by $(U_0,\bs x_t^\prime,\cdots,\bs x_1^\prime,Y_t,\cdots,Y_{1})$, where $U_0$ {is a random variable taking values in $(0,1)$}. In this work we are concerned with an observation-driven model in which the random component follows a conditional beta distribution, parameterized as \cite{Ferrari2004}:
\begin{align}\label{e:density}
	f(y;\mu_t,\nu|\F_{t-1})=\frac{\Gamma(\nu)}{\Gamma(\nu\mu_t)\Gamma\big(\nu(1-\mu_t)\big)}\,y^{\nu\mu_t-1}(1-y)^{\nu(1-\mu_t)-1},
\end{align}
for $0<y<1$, $0<\mu_t<1$ and $\nu>0$, where $\mu_t:=\E(Y_t|\F_{t-1})$. Observe that  $\var(Y_t|\F_{t-1})=\frac{\mu_t(1-\mu_t)}{1+\nu}$, so that $\nu$ acts as a precision parameter and that the model is conditionally heteroscedastic as the conditional variance depends on $\mu_t$. {However, since $\var(Y_t|\F_{t-1})\leq\frac1{4\nu}$, very high values of $\nu$ can account for conditional homoscedastic behavior in practice, as depicted in Figure \ref{f:chaos}} (this result is proven in Theorem \ref{resemble}).

To define the systematic component of the proposed model, let $T_{\bs\theta}:[0,1]\rightarrow[0,1]$  be a {dynamical system, i.e., a }function, potentially depending on an $r$-dimensional vector of parameters  $\bs\theta=(\theta_1,\cdots,\theta_r)^\prime\in\R^r$. {Let also} $g,h:(0,1)\rightarrow\R$  be two twice continuously differentiable link functions. In the additive specification \eqref{e:add_spec} we consider $\tau_t$  as a process in the form
\begin{align}\label{e:model}
	\eta_t:=g(\mu_t)=\alpha+\bs x_t^\prime\bs\beta+\sum_{j=1}^p\phi_j\big(g(y_{t-j})-\bs{x}_{t-j}'\bs{\beta}\big) + h\big(T_{\bs\theta}^{t-1}(U_0)\big),
\end{align}
where $\alpha\in \R$ is an intercept, $\bs\beta:=(\beta_1,\cdots,\beta_l)^\prime$ is an $l$-dimensional vector of parameter  associated to the covariates, $\bs\phi:=(\phi_1,\cdots,\phi_p)^\prime$ is a $p$-dimensional parameter related to the autoregressive structure in the model and {$U_0\in(0,1)$ is a random variable which will usually follows the absolute continuous invariant measure for the map $T_{\bs \theta}$, that will soon be introduced.
	Here $T_{\bs\theta}^{t-1}$ denotes the {$(t-1)$-th} iterate of the map (see next subsection).}
Specification \eqref{e:density} and \eqref{e:model}  define the proposed model, which we shall call beta autoregressive chaotic of order $p$ and denote by $\beta$ARC$(p)$ model. { As we shall see in the next sections, although the map $T_{\bs\theta}$  is defined in the closed interval $[0,1]$, usually $T_{\bs\theta}^t(U_0)$ takes values on the open interval $(0,1)$, for all $t$, with probability 1.}

{If we consider specification \eqref{e:model} without any covariate and without the autoregressive part,} the behavior of $\mu_t$ {(given by the orbit or sample path of the map $T_{\bs \theta}$)} often defines the overall behavior of the associated sample path. This means that the richness of possible sample paths {$\{T_{\bs\theta}^{t-1}(U_0)\}_{t \geq 1}$} in the class of all possible {chaotic } process {(which means all possible choices of maps $T_{\bs\theta}$)} can also be translated directly into the context of $\beta$ARC models.
Hence, the most interesting case of the proposed model occurs in the absence of covariates and the autoregressive parts. In that case, the {links $g$ and $h$ can be taken as the identity function} and the conditional average $\mu_t$ is driven solely by the behavior of the transformation $T_{\bs\theta}$ with the model's systematic component simplifying to
\begin{equation}\label{e:mutp0}
	\mu_t=T_{\bs\theta}^{t-1}(U_0).
\end{equation}
In what follows we shall refer to the $\beta$ARC model following \eqref{e:mutp0} as {the} pure chaotic $\beta$ARC models,
{while the  $\beta$ARC model following \eqref{e:model} where both $\bs\beta \neq 0$ and $\bs\phi\neq 0$ is called the full $\beta$ARC model.}

\subsection{Some definitions and results on Dynamical Systems}\label{ssd}

The proposed $\beta$ARC models strongly rely on the dynamic $T$.
{For this reason, in this section we introduce some standard definitions and results from one dimensional dynamic systems. More details and some references regarding dynamical systems can be found in the supplementary material.}

Let $T :[0,1] \to [0,1]$ and $x_0\in(0,1)$. We define the $k$-th iterate of $T$ as the $k$-fold composition $x_k:=T^k(x_0)=T\big(T^{k-1}(x_0)\big)$ and the sequence $\{x_0,x_1,x_2,\cdots\}$ is called the orbit (or sample path) of $T$.
A point $x$ is called a fixed point if  $T(x) = x$ and  is called a periodic point with period $s$ (where $s$ is a positive integer) if $T^{s}(x) = x$ and  $T^k(x) \neq x$, for all $k < s$. Fixed and periodic points can be very different in its nature: if $T$ is differentiable at a fixed point  $x$ we say that  $x$ is  attracting if  $|T'(x)| < 1$, repelling if $|T'(x)| > 1$, and indifferent (or neutral) if $|T'(x)| = 1$. Similar definitions hold for periodic points changing $T$ for $T^s$.

For a Borel measurable transformation $T:[0,1] \to [0,1]$, $\lambda_T$ is called a $T$-invariant probability measure (or invariant measure for short) if $\lambda_T$ is a probability measure defined on the Borel sets of $[0,1]$ and satisfies $\lambda_T\big(T^{-1}(A)\big) = \lambda_T(A)$ for all measurable set $A \subset [0,1]$.
If such invariant measure is absolutely continuous with respect to the Lebesgue measure, then we call it an ACIM.
\begin{rmk}
	{Whenever an ACIM $\lambda_T$ exists for a given map $T$, the natural choice for the distribution of $U_0$ is $\lambda_T$. If the map $T$ has an ACIM $\lambda_T$ and $U_0$  is chosen according to $\lambda_T$, then $T^t(U_0)\in(0,1)$, for all $t$, with probability 1, which is especially important for pure $\beta$ARC processes since $\mu_t$ must lie in $(0,1)$ for the model to be well defined. Another reason to favor maps which have ACIM is that, as we shall see in the sequel, this choice of distribution for $U_0$ allows for the derivation of several interesting results. }
\end{rmk}
An invariant measure $\lambda_T$ is called ergodic if the only measurable sets that are invariant for $T$ are sets of full or zero measure, i.e., if $T^{-1}(A)=A$ implies $\lambda_T(A)=0$ or $\lambda_T(A)=1$ (ergodicity implies that it is not possible to split the dynamics into two invariant sets with both having nonzero measure). Birkhoff Ergodic Theorem states that,  if $\lambda_T$ is ergodic for $T$, then  for any $\lambda_T$-integrable function $f:[0,1]\to \R$, and for  $\lambda_T$-almost all $x\in[0,1]$, we have
\[
\lim_{n\to \infty} \frac{1}{n} \sum_{k=0}^{n-1} f\big(T^k(x)\big)  = \int_0^1fd\lambda_T.
\]
In particular, by taking $f(x) = I_{A}(x)$, the indicator of $A$,  Birkhoff's theorem implies the convergence of the histogram (that is, the sample density) to the associated density \citep{DZ2009}.

{As we shall see in Section \ref{2.3}, the most interesting results for $\beta$ARC models are obtained when the map $T$ has an ACIM.
	However, the existence of such measure is not always guaranteed: for instance, if the map has attracting periodic orbits then it usually has no ACIM. Fortunately, a simple and easily verifiable condition ensures the existence of such measure: hyperbolicity.} We say that $T$ is uniformly expanding if $T$ is continuously differentiable and there exists $\rho>1$ such that  $|T'(x)| \geq \rho$, for all $x\in(0,1)$. This kind of maps are also called {\it hyperbolic maps}. Any fixed point of an uniformly expanding map is repelling. Also, if we require the derivative of such maps to be H\"older-continuous, then they present an unique ACIM, which gives positive mass to any open subset, and is also ergodic \citep[see][and references therein]{GB,HK,MVS}. 


\begin{figure}[h!]
	\centering
	\mbox{
		\hskip.1cm
		\subfigure[]{\includegraphics[width=0.3\textwidth]{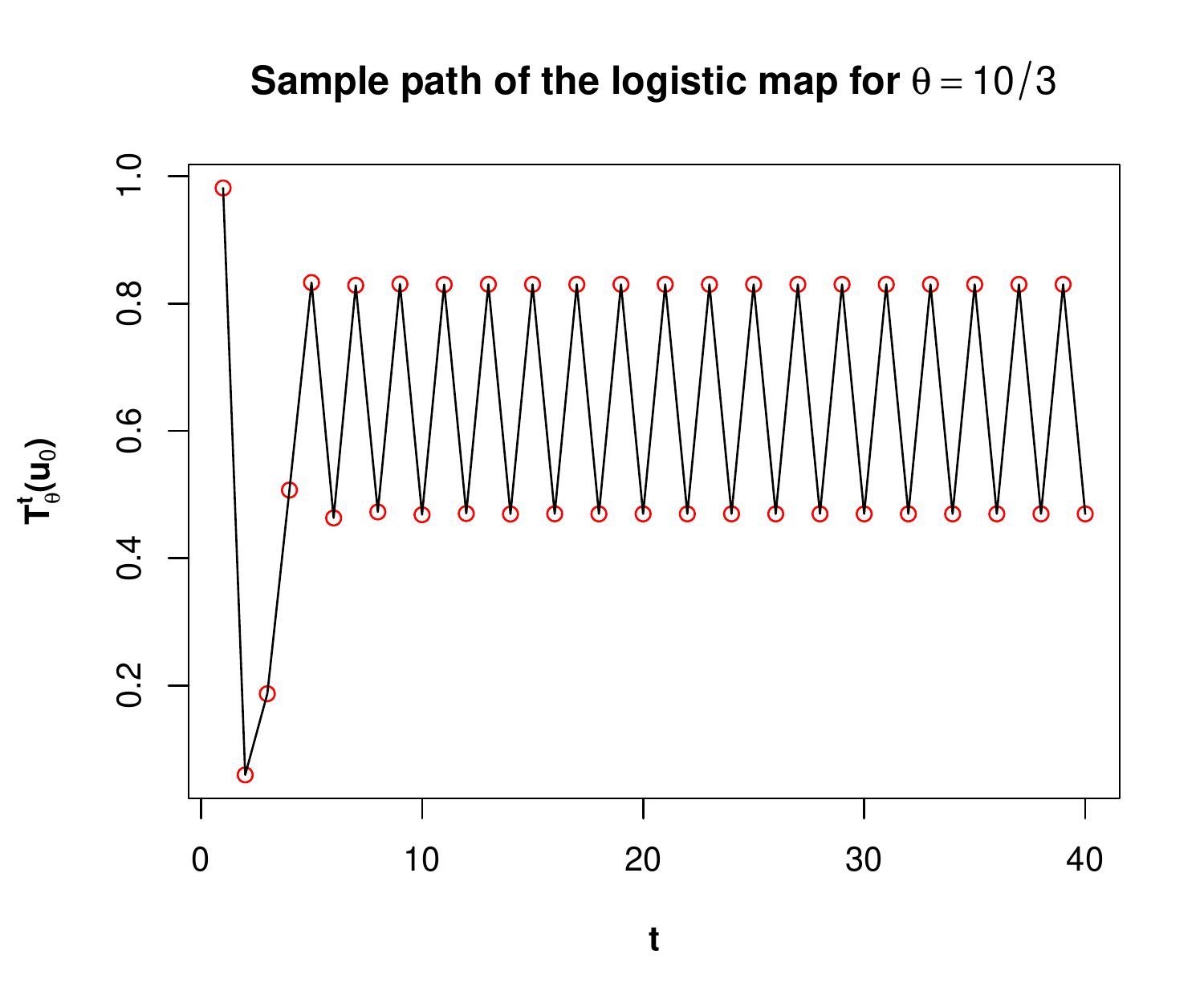}}\hskip.1cm
		\subfigure[]{\includegraphics[width=0.3\textwidth]{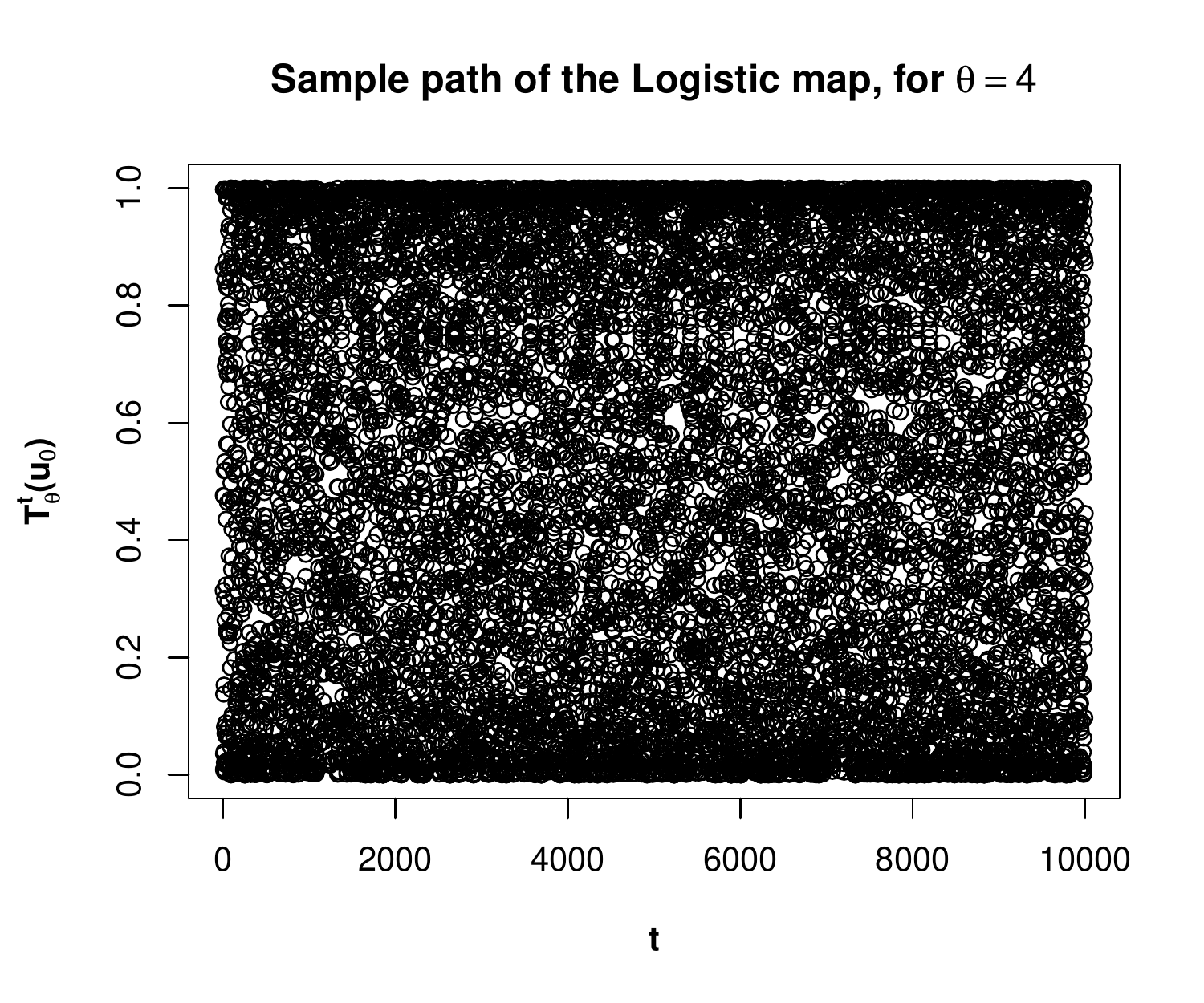}}\hskip.1cm
		\subfigure[]{\includegraphics[width=0.3\textwidth]{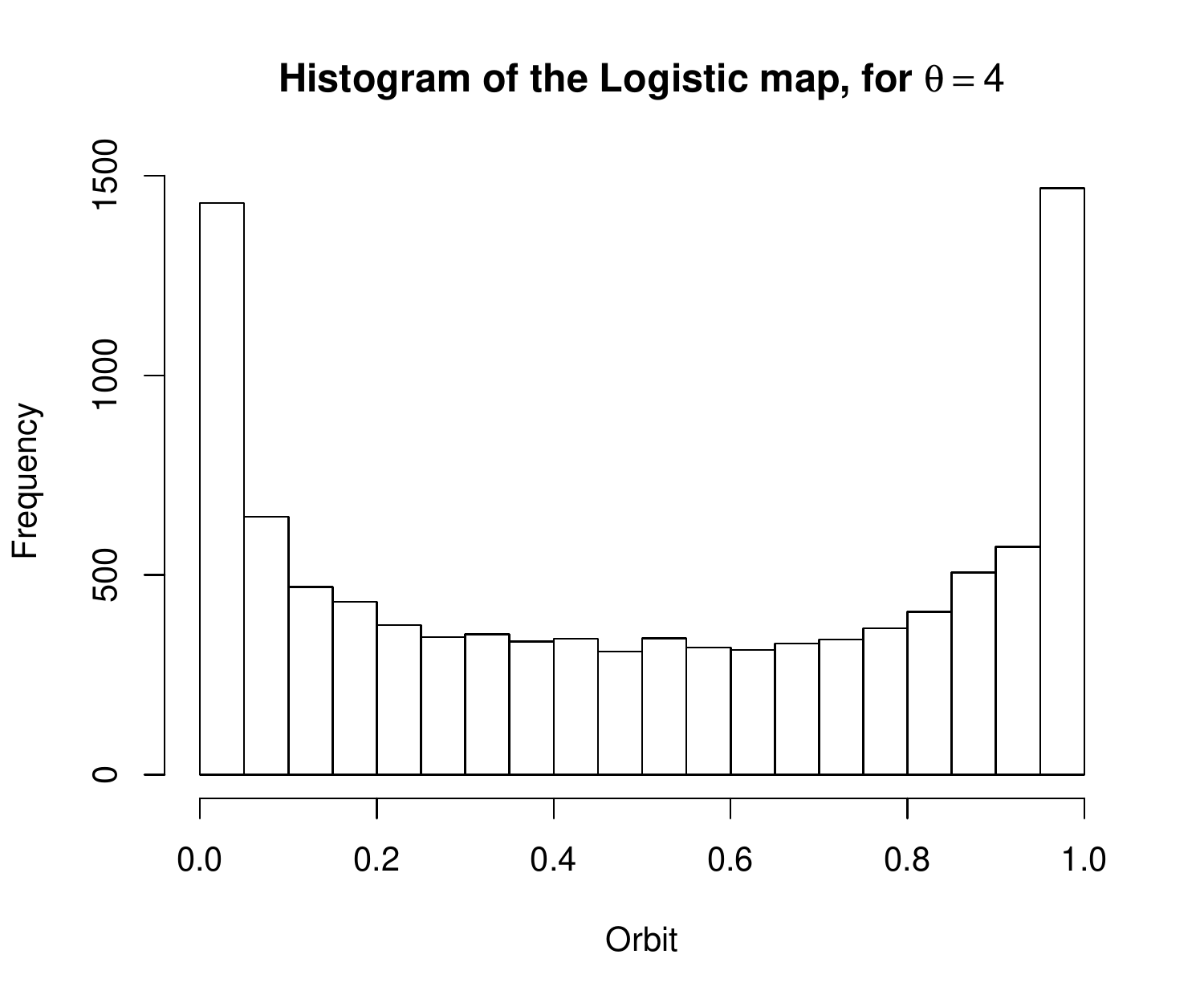}}
		\hskip.1cm
	}
	\caption{  (a) Plot of a sample path   associated  to  the logistic map $T_{\theta}(x)=\theta x(1-x)$ for  $\theta=10/3$ starting at $u_0=\pi/3.2$ showing an attracting periodic orbit of period 2; (b) sample path of the logistic map for $\theta=4$;  (c) Histogram of the sample path in (b).}\label{f:logist-4}
\end{figure}

\begin{figure}[htb!]	
	\centering
	\mbox{
		\subfigure[]{\includegraphics[width=0.3\textwidth]{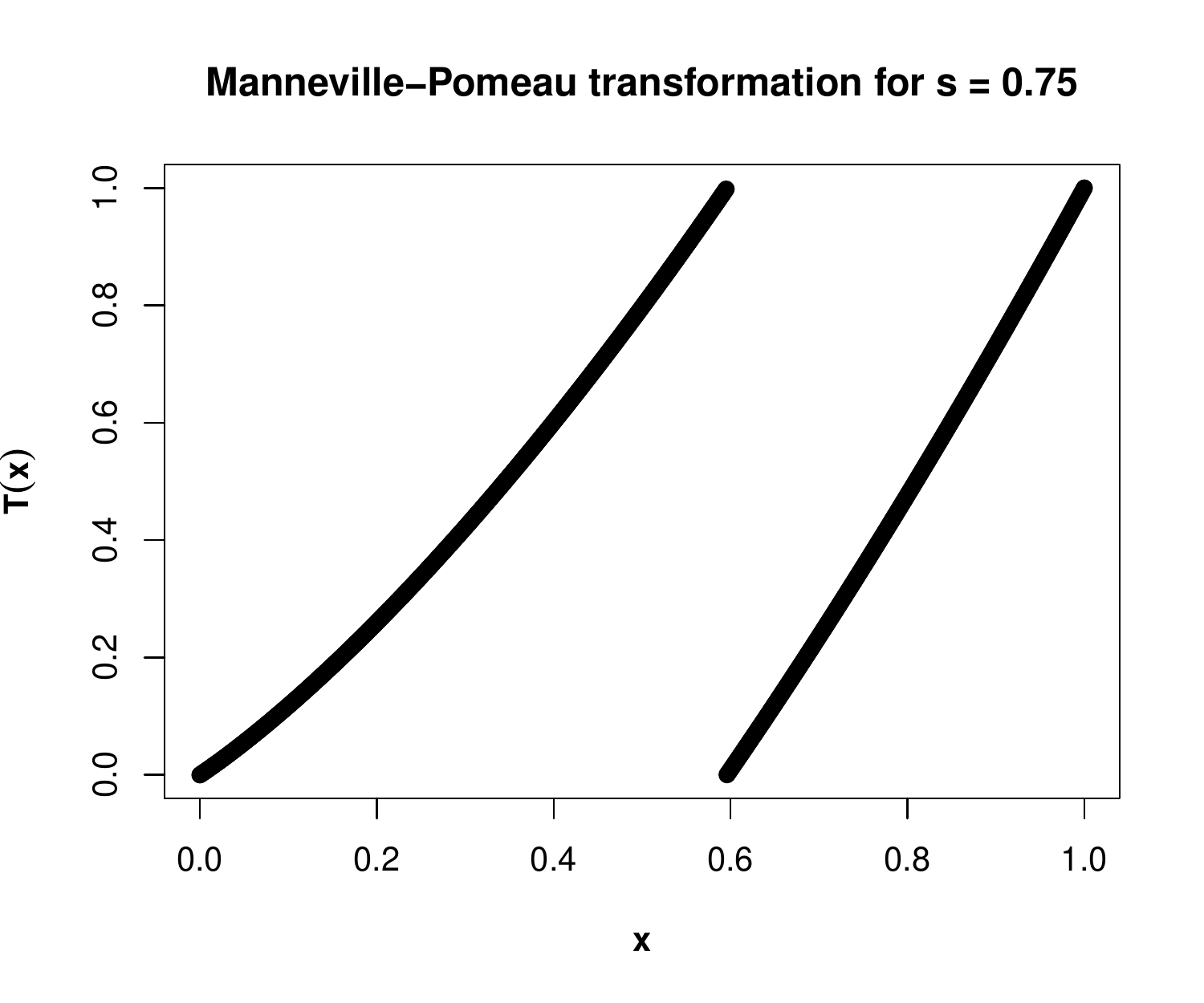}}
		\hskip.1cm
		\subfigure[]{\includegraphics[width=0.3\textwidth]{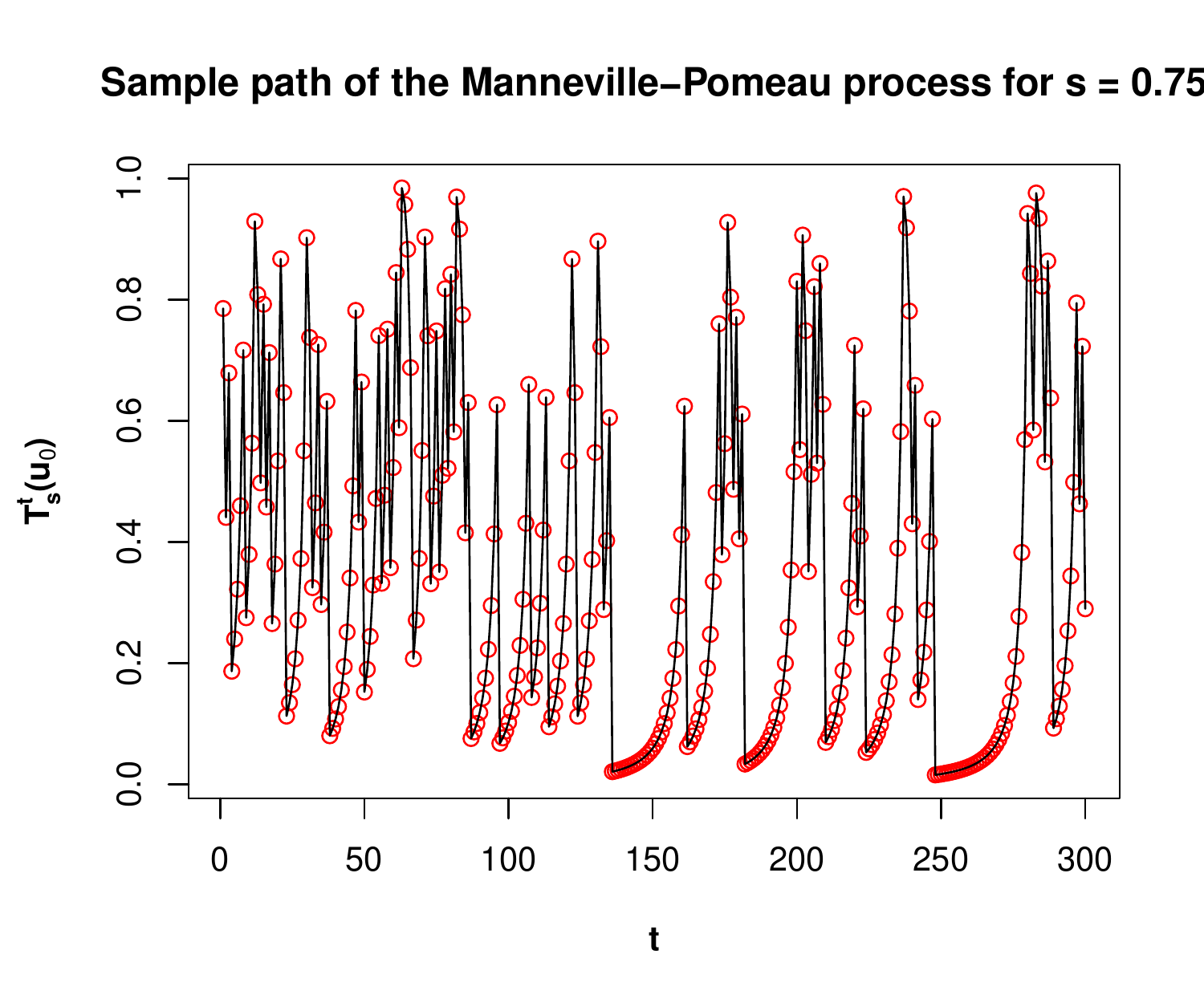}}
		\subfigure[]{\includegraphics[width=0.3\textwidth]{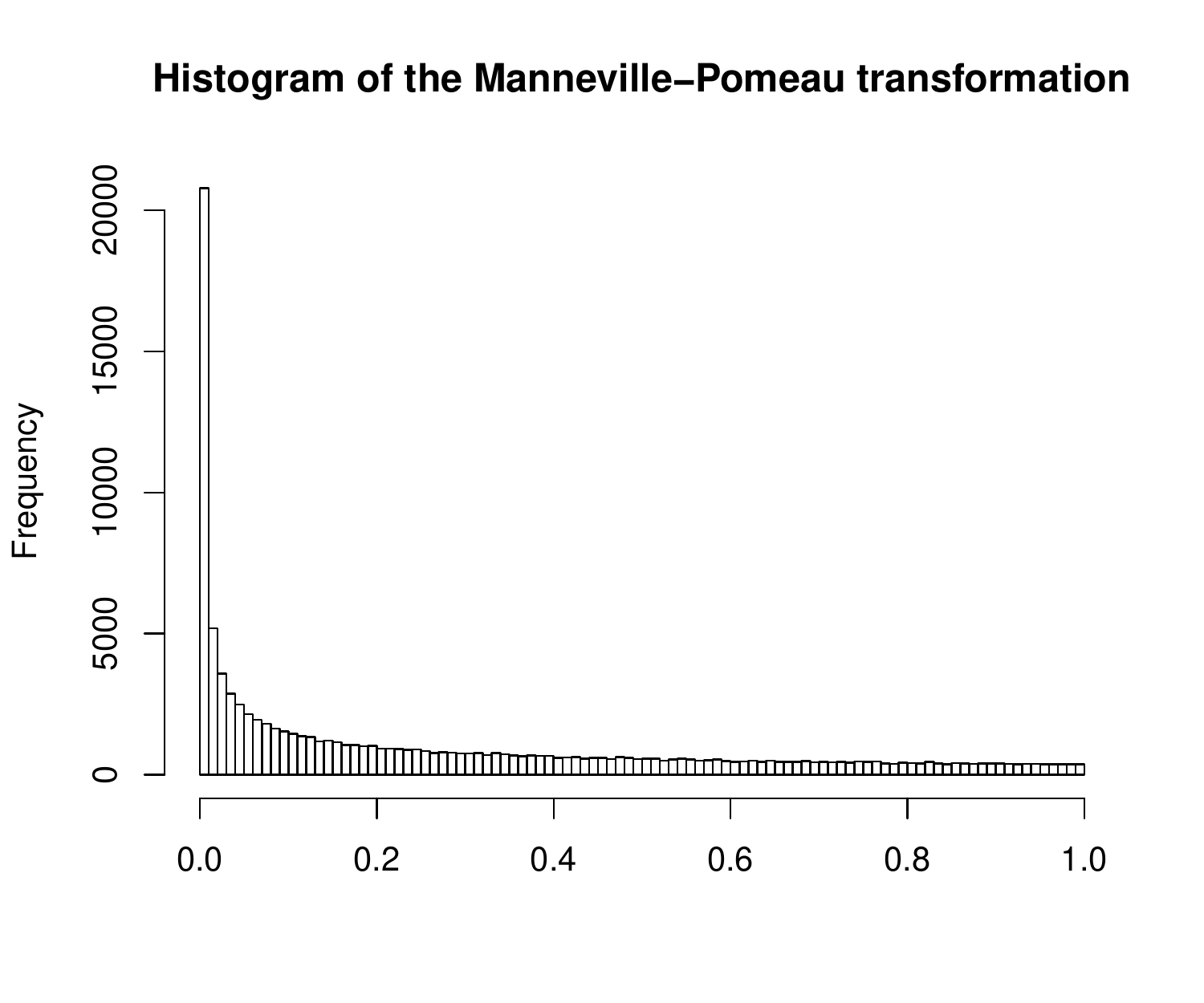}}
		\hskip.1cm
	}
	\caption{(a) A plot of the Manneville-Pomeau transformation \protect\eqref{e:MP} for $s=0.75$.  (b) Sample path of this map for $u_0=\pi/4$ showing laminar behavior near zero. (c) Histogram of the first 100,000 iterates of the map.}\label{f:MP}
\end{figure}

We will now present an example of a family of dynamical systems which will be used in our application in Section
\ref{se:application}. {Others examples can be found in the supplementary material. }  For $s>0$, the Manneville-Pomeau transformation $T_s:[0,1]\rightarrow [0,1]$, is given by
\begin{equation}\label{e:MP}
	T_s(x)=(x+x^{1+s})(\mathrm{mod}\,1).
\end{equation}
For $s\in(0,1)$, there exists an absolutely continuous $T_s$-invariant probability measure \citep{Thaler1980}, which can be seen in Figure \ref{f:MP}(c). For $s\geq1$ there exists an absolutely continuous invariant measure which is only $\sigma$-finite (not a probability measure). Figure \ref{f:MP}(a) show the Manneville-Pomeau transformation for $s=0.75$.
The Manneville-Pomeau transformation presents a property referred to as transition to turbulence through intermittency \citep{Eck1981}.
The Manneville-Pomeau transformation has an  indifferent fixed point  at $0$ and, hence, it is not uniformly expanding. The chaotic processes associated to $T_s$ are often called Manneville-Pomeau processes which present a very slow correlation decay when $s\in(0.5,1)$, characteristic of long range dependent processes and it is commonly viewed as an alternative model for long range dependence outside the classical duet of Fractional Brownian Motion and ARFIMA processes. This slow decay is mainly due to the presence of laminar behavior near zero, which can be seen in Figure \ref{f:MP}(b).

\subsection{General results on the $\beta$ARC model}\label{2.3}
{   We will now present some general results about the $\beta$ARC model. Some of them concern the pure $\beta$ARC model, some other involve the presence of the dynamics alongside random covariates (but no AR structure) while others deal with the full $\beta$ARC model. }
We begin with Proposition \ref{cov}. This result, for the pure $\beta$ARC model, shows that the covariance structure  of the dynamics $\{\mu_t\}_{t\geq1}$ ultimately determines the unconditional covariance structure of the process $\{Y_t\}_{t\geq1}$. To the best of our knowledge, similar results are not available for any other competing GARMA-like processes in the literature.
\begin{prop}\label{cov}
	Let $\{Y_t\}_{t\geq1}$ be a pure $\beta$ARC  process following \eqref{e:mutp0}. Then, for all $t,h>0$,
	\begin{enumerate}
		\item $\E(Y_t)=\E(\mu_t)$.
		\item $\var(Y_t)=\var(\mu_t)+\dfrac1{1+\nu}\E\big(\mu_t(1-\mu_t)\big)$.
		\item $\mathrm{Cov}(Y_t,Y_{t+h})=\mathrm{Cov}(\mu_t,\mu_{t+h})$.
	\end{enumerate}
\end{prop}
\proof Item (a) follows from $\E(Y_t)=\E\big(\E(Y_t|\F_{t-1})\big)=\E(\mu_t)$, while (b) follows from the identity $\var(Y_t)=\var\big(\E(Y_t|\F_{t-1})\big)+\E\big(\var(Y_t|\F_{t-1})\big)$ and (a).
\noindent As for (c), for any $t,h>0$,  from item (a) we obtain
\begin{align}\label{e:auxcov}
	\mathrm{Cov}(Y_t,Y_{t+h})=\E(Y_tY_{t+h})-\E(Y_t)\E(Y_{t+h})=\E(Y_tY_{t+h})-\E(\mu_t)\E(\mu_{t+h}).
\end{align}
Now, notice that $\mu_t$ is $\F_{1}$-measurable, for all $t>0$, so that
\begin{align}\label{e:auxcov2}
	\E(Y_tY_{t+h})&=\E\big(\E(Y_tY_{t+h}|\F_{t+h-1})\big)=\E\big(Y_t\E(Y_{t+h}|\F_{t+h-1}))=
	\E(Y_t\mu_{t+h})\nonumber\\
	&=\E\big(\E(Y_t\mu_{t+h}|\F_{t-1})\big)=\E\big(\mu_{t+h}\E(Y_{t}|\F_{t-1})\big)=
	\E(\mu_t\mu_{t+h}),
\end{align}
and the result follows upon replacing \eqref{e:auxcov2} into \eqref{e:auxcov}.\qed

Conditions for stationarity of {other} dynamical models for time series following a GARMA approach are traditionally very hard to obtain and remain an open subject for most traditional models, such as $\beta$ARMA \citep{Rocha2009}, $\beta$ARFIMA \citep{Pumi2019} and KARMA models \citep{Bayers} {except under trivial scenarios}. {In the next result we show that $\beta$ARC models  are stationary in a very broad specification and under easily verifiable conditions. }

\begin{thm}\label{thm1}
	Let $\{Y_t\}_{t\geq 1}$ be a $\beta$ARC model with $\nu>0$ and
	\[\eta_t=g(\mu_t)=\alpha+\bs x_t^\prime\bs\beta+h\big(T^{t-1}(U_0)\big),\]
	where $\{\bs x_t\}_{t\geq 1}$ is a set of random covariates, $g$ and $h$  are twice continuously differentiable, one to one link functions, {and $U_0$ is a random variable such that $T^t(U_0)\in(0,1)$, for all $t \geq 0$, with probability 1}. Then $\{(Y_t,\mu_t)\}_{t\geq 1}$ is jointly stationary if and only if $\{ \mu_t\}_{t\geq 1}$ is stationary.
\end{thm}
\proof
Suppose that $\{\mu_t\}_{t\geq1}$ is stationary.  For any arbitrary positive integer $k$, let $t_1,\cdots,t_k$ be distinct time points, $\bs t = (t_1, \cdots, t_k)$, $\bs t + h = (t_1+h, \cdots, t_k+h)$,  $\bs{Y_t}=(Y_{t_1},\cdots,Y_{t_k})$ and $\bs{\mu_t} = (\mu_{t_1},\cdots,\mu_{t_k})$.  Using Riemman-Stieltjes integration we have
\[
F_{Y_t|\mu_t}(y|z) = P(Y_t \leq y |\mu_t = z) = \int_{0}^ydF_{Y_t|\mu_t}(u|z), \quad \forall y,z \in(0,1),
\]
\[
F_{\bs{\mu_t}}( v_1, \dots, v_k)  = P(\mu_{t_1}\leq v_1, \dots, \mu_{t_k}\leq v_k) =  \int_{0}^{v_1}\dots \int_{0}^{v_k}dF_{\bs{\mu_t}}(z_1,\dots,z_k)
\]
and
\begin{align*}
	F_{\bs{Y_t}, \bs{\mu_t}}(u_1, \dots, u_k, v_1, \dots, v_k) & = P(Y_{t_1} \leq u_1, \dots, Y_{t_k} \leq u_k, \mu_{t_1}\leq v_1, \dots, \mu_{t_k}\leq v_k)\\
	& = \int_{0}^{u_1}\dots \int_{0}^{u_k}\int_{0}^{v_1}\dots \int_{0}^{v_k}dF_{\bs{Y_t}, \bs{\mu_t}}(y_1,\dots,y_k,z_1,\dots,z_k),
\end{align*}
for all $u_1, \dots, u_k, v_1, \dots, v_k \in (0,1)$, where $dF_{Y_t,\mu_t}$,  $dF_{\bs{\mu_t}}$ and $dF_{\bs{Y_t}, \bs{\mu_t}}$ are the integrands of the Riemann-Stieltjes integrals.  Observe that,   for all $t > 0$, given  $\mu_{t} = z$, the random variable $Y_{t}$ depends, neither on the past information $\{Y_{s}, \mu_{s}\}_{s < t}$, nor on the future $\mu_s$, $s > t$, so that
\[
dF_{\bs{Y_t}, \bs{\mu_t}}(y_1,\dots,y_k,z_1,\dots,z_k) = dF_{\bs{\mu}_{\bs{t}}}(z_1,\dots,z_k) \prod_{j=1}^kdF_{Y_{t_j}|\mu_{t_j}}(y_j|z_j).
\]
It is easy to see that $dF_{Y_t|\mu_t}(y|z) = f(y;z,v|\F_{t-1})dy$, where $f$ is the conditional density defined by \eqref{e:density}  and that $dF_{Y_t|\mu_t}(y|z) = dF_{Y_1|\mu_1}(y|z) = dF_{Y_{t+h}|\mu_{t+h}}(y|z) $, for all $t,h>0$, so that, from the stationarity of $\{\mu_t\}_{t\geq 1}$, it follows that,  for all $A_i, B_i \subset (0,1)$, $i = 1,\dots,k$,
\begin{align*}
P(Y_{t_1} \in A_1,\dots, &Y_{t_k} \in A_k, \mu_{t_1} \in B_1, \dots, \mu_{t_k} \in  B_k)  = \\
&=\int_{A_1}\dots \int_{A_k}\int_{B_1}\dots \int_{B_k}dF_{\bs{Y_t}, \bs{\mu_t}}(y_1,\dots,y_k,z_1,\dots,z_k)\\
 &= \int_{A_1}\dots \int_{A_k}\int_{B_1}\dots \int_{B_k}dF_{\bs{\mu}_{\bs{t}}}(z_1,\dots,z_k) \prod_{j=1}^kdF_{Y_{t_j}|\mu_{t_j}}(y_j|z_j)\\
 &= \int_{A_1}\dots \int_{A_k}\int_{B_1}\dots \int_{B_k}dF_{\bs{\mu}_{\bs{t}+h}}(z_1,\dots,z_k) \prod_{j=1}^kdF_{Y_{t_j+h}|\mu_{t_j+h}}(y_j|z_j)\\
&= P(Y_{t_1+h} \in A_1,\dots, Y_{t_{k}+h} \in A_k, \mu_{t_1+h} \in B_1, \dots, \mu_{t_k+h} \in  B_k) .
\end{align*}
This implies that $\{(Y_t,\mu_t)\}_{t\geq 1}$ is jointly stationary.  The converse is obvious.
\qed


\begin{coro}\label{corozero}
	Under the conditions of Theorem \ref{thm1}, if $\{ \mu_t\}_{t\geq 1}$ is stationary, then so is $\{Y_t\}_{t\geq 1}$.
\end{coro}
\proof Observe that, if $\{\mu_t\}_{t\geq 1}$ is stationarity then from Theorem \ref{thm1} $\{(Y_t,\mu_t)\}_{t\geq 1}$ is jointly stationary and hence, $\{Y_t\}_{t\geq 1}$ is stationary.
\qed

\begin{coro}
	Let $T_{\bs \theta}$ be a dynamical system with ACIM given by $\lambda_T$ and let $\{Y_t\}_{t\geq 1}$ be a pure chaotic $\beta$ARC model with $\nu>0$ where $\mu_t=T_{\bs \theta}^{t-1}(U_0)$ and $U_0$ is chosen accordingly to $\lambda_T$. Then
	$\{Y_t\}_{t\geq 1}$ is stationary and the common marginal distribution $F_{Y_t}$ is absolutely continuous with respect to the Lebesgue measure, with unconditional density  given by
	\[
	f_{Y_t}(y)=\int_0^1 f_{Y_t|\mu_t}(y|z)\lambda_T(dz),
	\]
	where $f_{Y_t|\mu_t}(y|z) = f(y;z,\nu|\F_{t-1})$ is the conditional density of $Y_t$ given $\mu_t$, defined by \eqref{e:density}.
\end{coro}

\proof
The stationarity of  $\{Y_t\}_{t\geq 1}$
follows immediately from Corollary \ref{corozero}, as $\mu_t$ is clearly stationary in this case.
Now, let $f_{Y_t,\mu_t}$ denote the joint density of $\big(Y_t, \mu_t\big)$, so we have
\[
f_{Y_t}(y)=\int_0^1 f_{Y_t,\mu_t}(y,z)dz = \int_0^1 f_{Y_t|\mu_t}(y|z)\lambda_T(dz),
\]
and the proof is complete \qed.

\begin{coro}
	Let $T_{\bs \theta}$ be a dynamical system with ACIM given by $\lambda_T$, and let $\{\bs x_t\}_{t\geq 1}$ be a set of random covariates. Suppose $\{Y_t\}_{t\geq 1}$ is a $\beta$ARC model with $\nu>0$ where
	\[\eta_t=g(\mu_t)=\alpha+\bs x_t^\prime\bs\beta+h\big(T^{t-1}(U_0)\big),\]
	for two twice continuously differentiable, one to one link functions $g$ and $h$.
	Suppose $U_0$ is chosen according to $\lambda_T$.
	Then if $\{\bs x_t\}_{t\geq 1}$ is stationary, so is $\{Y_t\}_{t\geq 1}$.
\end{coro}
\proof Since $g$ and $h$ are both measurable functions, $\{\mu_t\}_{t\geq 1}$ is stationary if and only if  $\{\bs x_t\}_{t\geq 1}$ is stationary and the result follows immediately from Theorem \ref{thm1}. \hfill \qed

\begin{rmk}
	The proof of Theorem \ref{thm1} is also valid under the full specification \eqref{e:model}. However, verification of the hypothesis under \eqref{e:model} is difficult since it is not presented in an autoregressive fashion as we write $\eta_t$ in terms of past values of $Y_t$ and $\bs x_t$, which depends on the past of $\eta_t$ in a non-trivial way. In this scenario it is challenging to obtain stationarity conditions for $\{\eta_t\}_{t\geq 1}$ under the full specification \eqref{e:model}. This and the recursive nature of $\mu_t$ for similar GARMA-like models, such as the $\beta$ARMA, $\beta$ARFIMA and KARMA, make obtaining stationarity conditions a non-trivial problem for these models.
\end{rmk}

\begin{figure}[h!]
	\centering
	\mbox{
		\subfigure[]{\includegraphics[width=0.45\textwidth]{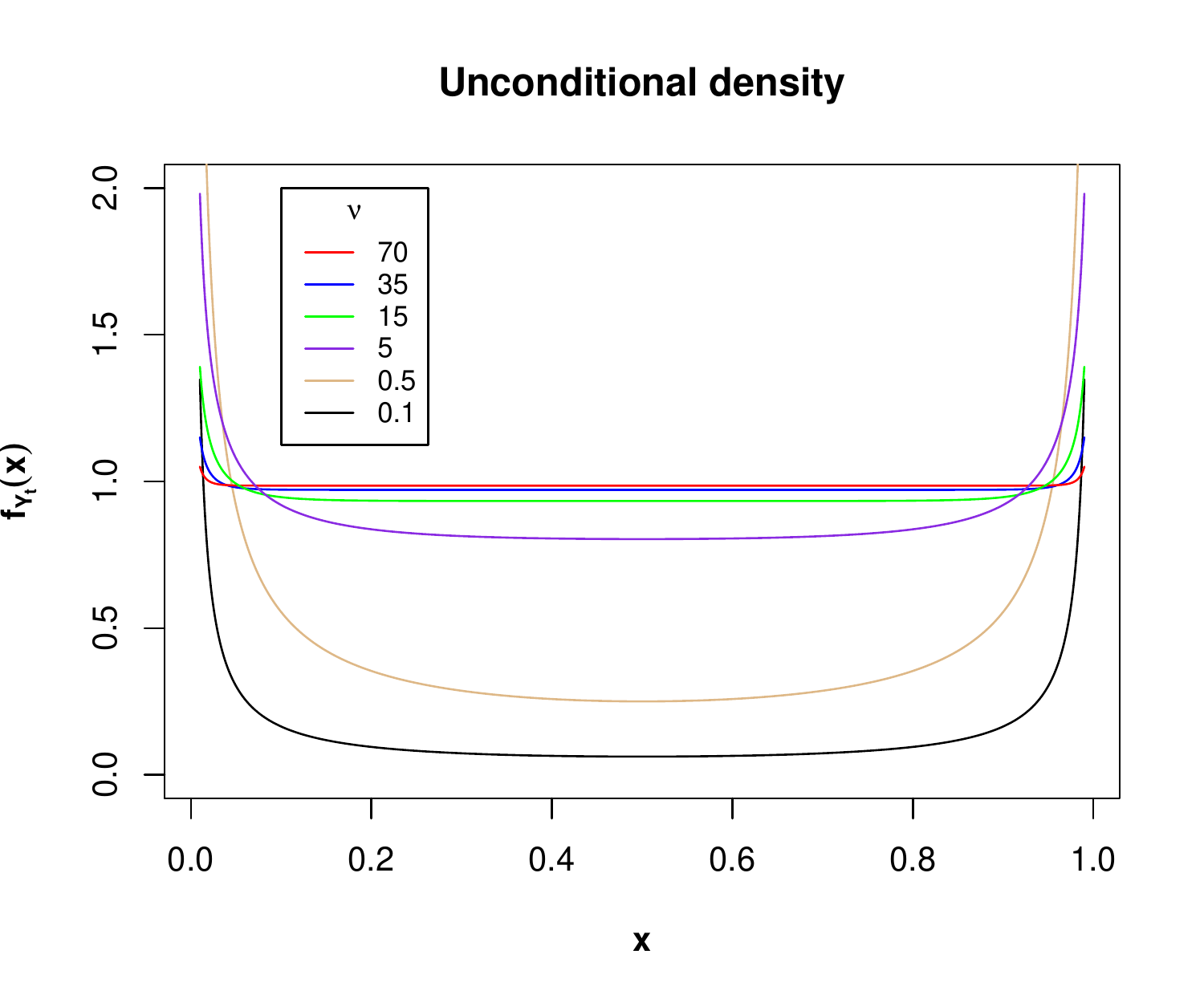}}		\hskip.5cm
		\subfigure[]{\includegraphics[width=0.45\textwidth]{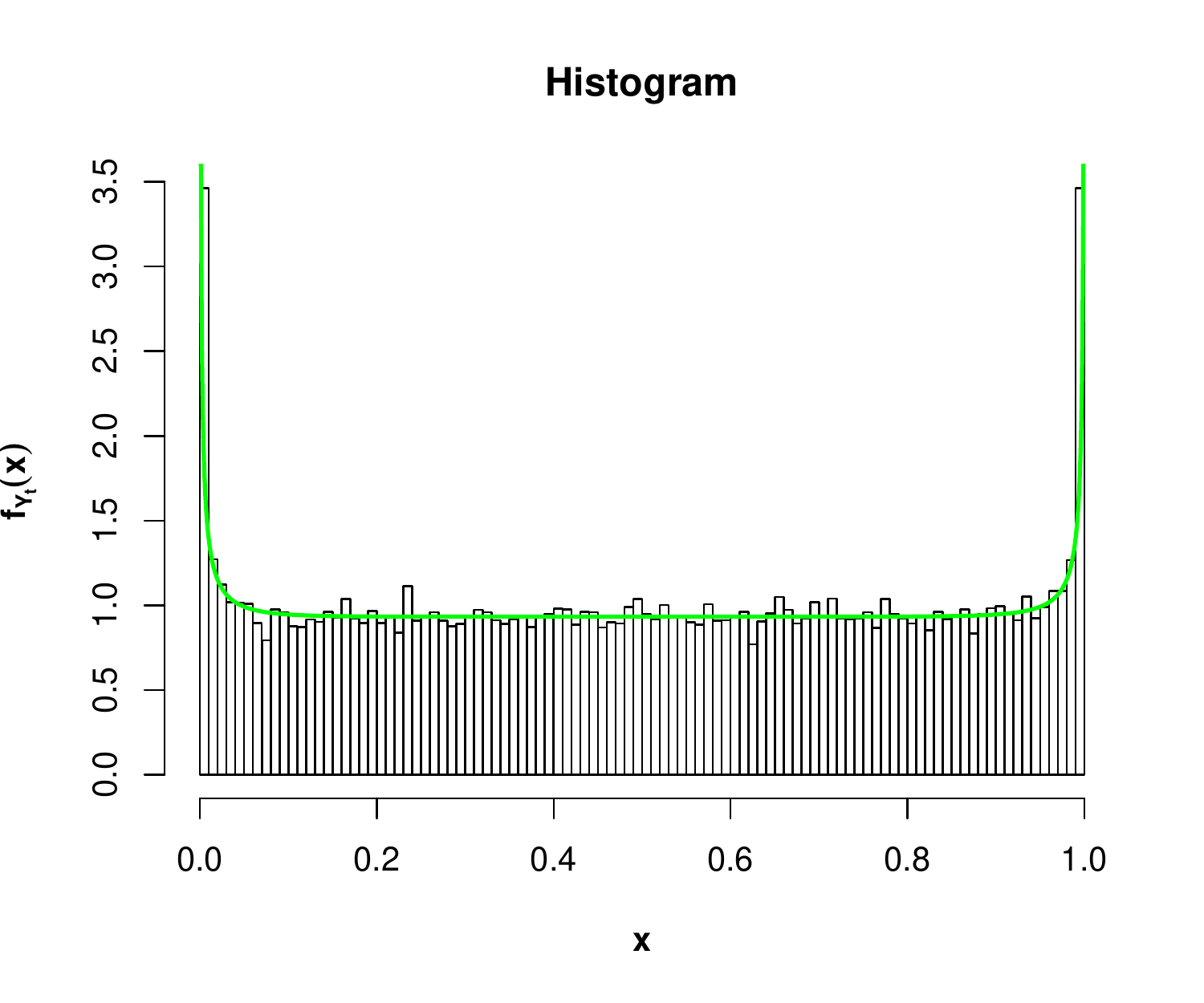}}
	}
	\caption{(a) The unconditional density of the pure $\beta$ARC model with $T_k(x)=(kx) \mathrm{mod}(1)$ and (b) histogram of an associated sample of size $n=30,000$ starting at $u_0=\pi/4$, with $\nu=15$ and $k=3$, showing the associated unconditional density (green).}
	\label{exe}
\end{figure}

{Before moving on, let us analyze an example that will motivate the next result.
	We shall analyze  the stationarity of the $\beta$ARC model with $T_k(x)=(kx) \mathrm{mod}(1)$} for an integer $k>0$. In this case, the Lebesgue measure in $[0,1]$ is $T_k$ invariant and the unconditional distribution of $Y_t$ is given by
\[f_{Y_t}(x)=\frac{\Gamma(\nu)(1-x)^{\nu-1}}{x}\int_0^1\left[\frac{x}{1-x}\right]^{\nu z}\frac{1}{\Gamma(\nu z)\Gamma\big(\nu(1-z)\big)} dz.\]
The behavior of $f_{Y_t}$ depends on the magnitude of $\nu$. In Figure \ref{exe} we show the behavior for several values of $\nu$.


\FloatBarrier

\begin{figure}
	\centering
	\mbox{
		\subfigure[]{\includegraphics[width=0.45\textwidth]{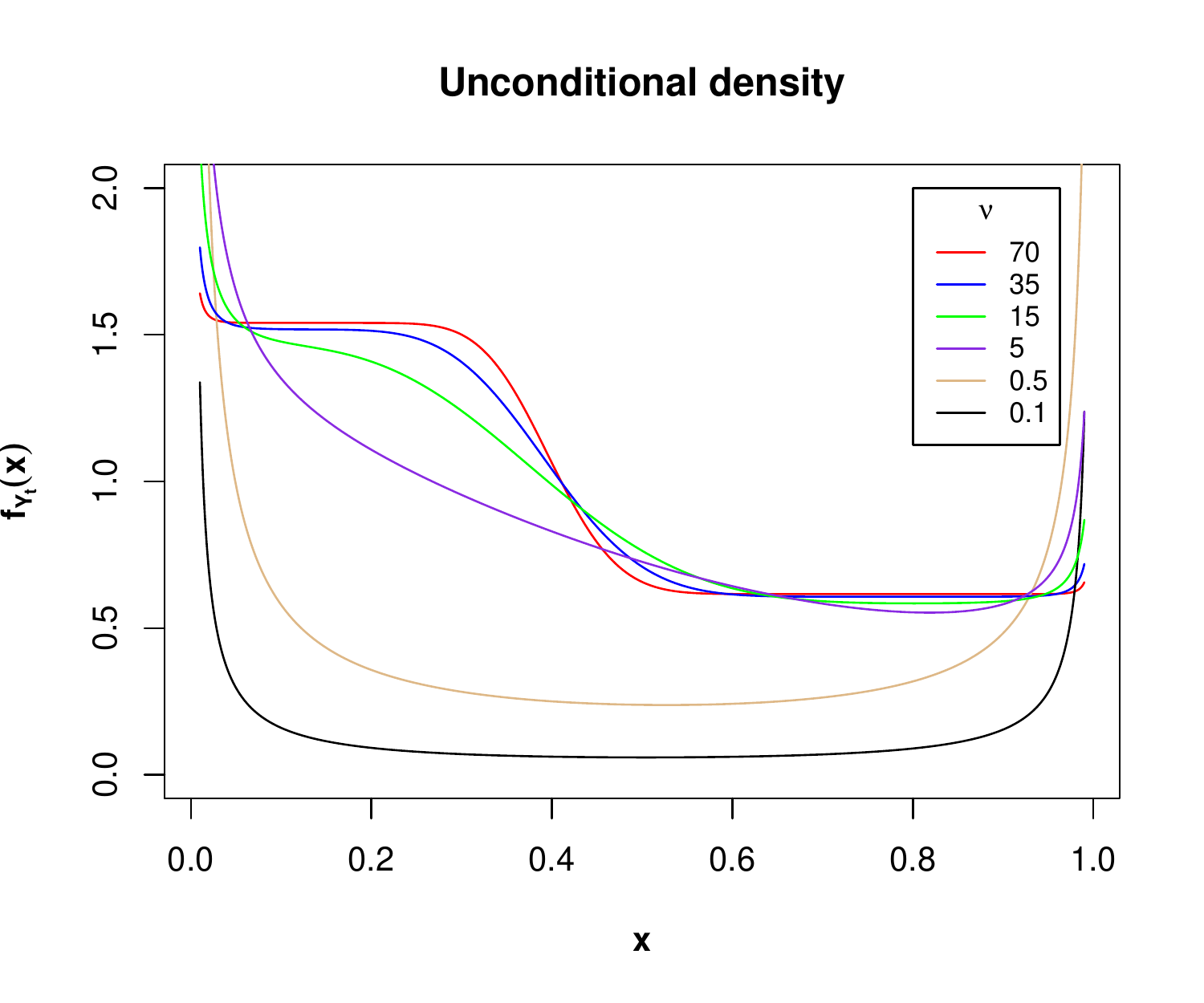}}		\hskip.5cm
		\subfigure[]{\includegraphics[width=0.45\textwidth]{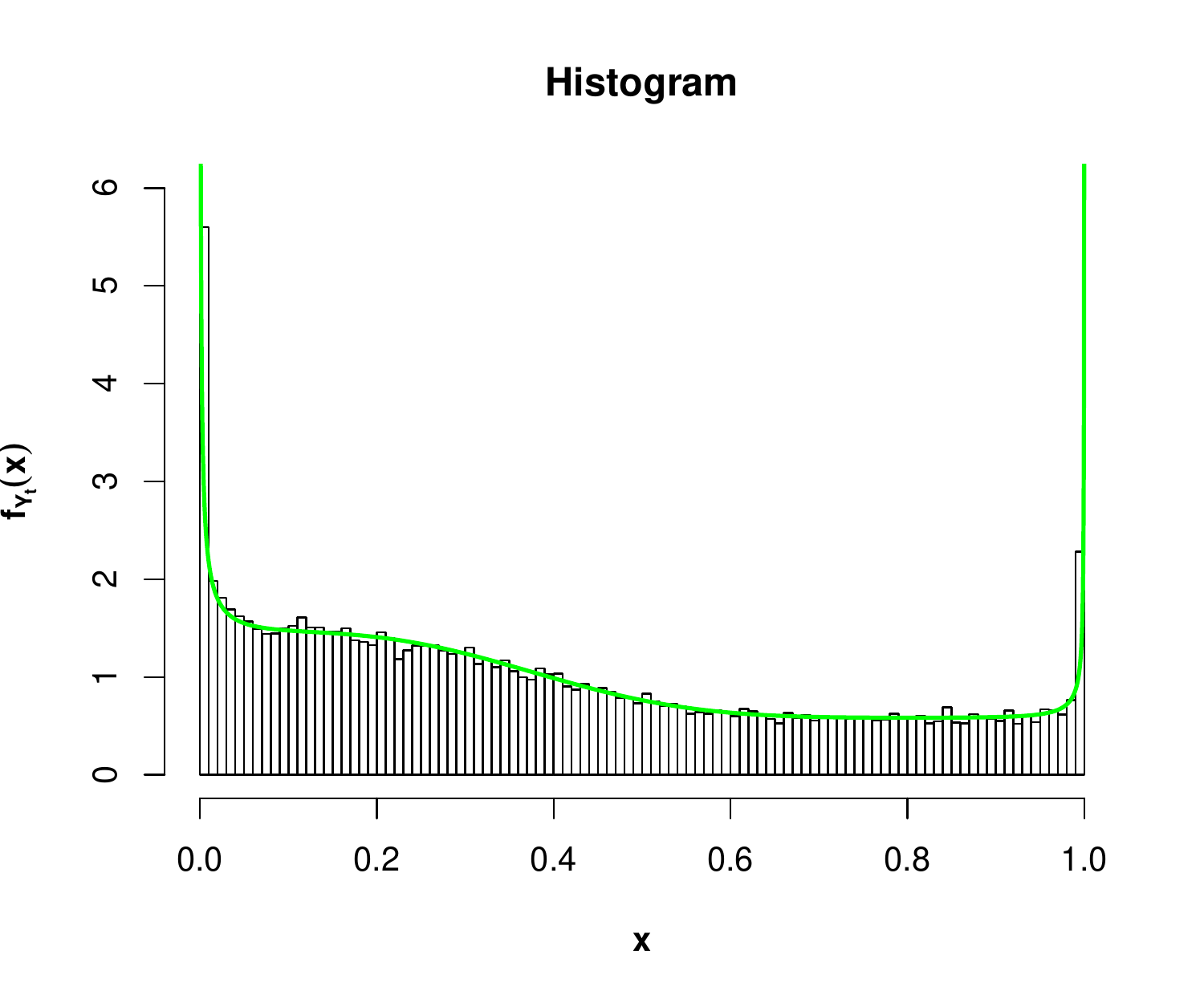}}
	}
	\caption{(a) The unconditional density of the pure $\beta$ARC model for \protect\eqref{e:dissertRafael} and (b) histogram of an associated sample of size $n=30,000$ starting at $u_0=\pi/4$, with $\nu=15$ and $\theta=0.4$ showing the associated unconditional density (green).}
	\label{exe2}
\end{figure}

{Let
	\small
	\begin{equation}\label{e:dissertRafael}
		T_{\theta}(x)=\begin{cases}
			\frac{x}{\theta} & \mbox{if} \;\; 0 \leq x < \theta, \\
			\frac{\theta (x- \theta)}{1-\theta} & \mbox{if} \;\; \theta \leq x \leq 1.
		\end{cases}
	\end{equation}
	\normalsize
	More details regarding this map can be found in the Supplementary material (Map 2).}
Now consider the pure $\beta$ARC model coupled with \eqref{e:dissertRafael}. For any $\theta\in(0,1)$, the unconditional distribution of $Y_t$ is given by
\begin{align*}
	f_{Y_t}(x)=\frac{\Gamma(\nu)(1-x)^{\nu-1}}{x}\bigg(\int_0^x \left[\frac{x}{1-x}\right]^{\nu z}&\frac{1}{(2-z)\Gamma(\nu z)\Gamma\big(\nu(1-z)\big)} dz+\\
	+&	\int_x^1\left[\frac{x}{1-x}\right]^{\nu z}\frac{1}{z(2-z)\Gamma(\nu z)\Gamma\big(\nu(1-z)\big)} dz\bigg).
\end{align*}
In Figure \ref{exe2} we show the behavior of $f_{Y_t}$ for several values of $\nu$.

{The next result shows, as the last example suggests, that the larger the precision parameter is}, the closer the $\beta$ARC model resembles its conditional mean.
{Observe that the result does not require stationarity to hold.}

\begin{thm}\label{resemble}
	Let $\{Y_t\}_{t\geq 1}$ be a $\beta$ARC process.
	Then, for each fixed $t>0$,
	\[Y_t\overset{d}{{\longrightarrow}}\mu_t, \quad \mbox{as } \nu\to\infty.\]
\end{thm}
\proof
First, for fixed $t>0$, observe that $\var(Y_t|\F_{t-1})=\frac{\mu_t(1-\mu_t)}{1+\nu}\rightarrow 0$ as $\nu\to\infty$.
Now we can also use the fact that $\E(Y_t|\F_{t-1})=\mu_t$ and Chebysheff's inequality to conclude that $Y_t$ conditionally converges in probability to $\mu_t$, which implies convergence in distribution.	
Therefore, for any $0<c<1$ which is a continuity point of $F_{\mu_t}$, we have
\[ P( Y_t \leq c | \mu_t=z) \longrightarrow
\begin{cases}
1 \mbox{ if } c>z,\\ 0 \mbox{ if } c<z,
\end{cases}
\]
when $\nu\to\infty$, which implies
\[P( Y_t \leq c) = \int_0^1  P( Y_t \leq c | \mu_t=z) dF_{\mu_t}(z)
\rightarrow \int_0^c 1\, dF_{\mu_t}(z)=  P(\mu_t \leq c)
\] when $\nu\to\infty$, by the Lebesgue dominated convergence theorem.\qed

In the case of the pure $\beta$ARC model, if the map $T_{\bs \theta}$ has an ACIM $\lambda_T$ and $U_0$ is distributed according to $\lambda_T$, the distribution of $\mu_t$ is given by $\lambda_{T}$. Therefore, Theorem \ref{resemble} and Birkhoff's Theorem suggests that the histogram is a valuable tool in choosing the family of maps $T_{\bs \theta}$ to be used to model a given time series.

The next theorem presents a simple condition under which the strong law of large numbers holds for $\beta$ARC process.
In particular, for a stationary $\beta$ARC process, the strong law of large numbers for $\{Y_t\}_{t\geq 1}$ is related to the covariance structure of the dynamical system $\{\mu_t\}_{t\geq 1}$.

\begin{thm}\label{t:as}
	Let $\{Y_t\}_{t\geq1}$ be a $\beta$ARC process  and $\varphi:[0,1]\to\R$ be a measurable function such that $\E(\varphi(Y_t)^2)<\infty$, for all $t>0$. If
	\begin{equation}\label{cv}
		\sum_{k=1}^\infty \frac{\sup_{t\geq 1}\Big\{\big|\mathrm{Cov}\big(\varphi(Y_t),\varphi(Y_{t+k})\big)\big|\Big\}}{k^q}<\infty,\quad\mbox{for some $0\leq q<1$},
	\end{equation}
	and
	\begin{equation}\label{var}
		\sum_{k=1}^\infty \frac{\var\big(\varphi(Y_k)\big)\ln(k)^2}{k^2}<\infty
	\end{equation}
	then
	\begin{equation}\label{as}
		\lim_{n\to \infty} \frac{1}{n} \sum_{\ell=0}^{n-1} \big[\varphi(Y_\ell) - \E\big(\varphi(Y_\ell)\big)\big]=0,
		\quad \mbox{a.s.}
	\end{equation}
\end{thm}
\proof {Observe that, conditions \eqref{cv} and \eqref{var} translate into conditions (3.2) and (3.1) in theorem 1 in \cite{hu2008}, respectively. Hence, the result in the mentioned theorem hold which translates into  \eqref{as}.}\qed

\begin{rmk}
	Observe that if $\{Y_t\}_{t\geq1}$ is stationary, so is $\{\varphi(Y_t)\}_{t>0}$, hence \eqref{var} is always satisfied and condition \eqref{cv} becomes
	\[\sum_{k=1}^\infty \frac{\big|\mathrm{Cov}\big(\varphi(Y_t),\varphi(Y_{t+k})\big)\big|}{k^q}<\infty,\quad\mbox{for some $0\leq q<1$}.\]
	Moreover, in this case the conclusion is a Birkhoff-type theorem since \eqref{as} becomes
	\[\lim_{n\to \infty} \frac{1}{n} \sum_{\ell=0}^{n-1} \varphi(Y_\ell) = \int \varphi(z) dF_{Y_t}(z) \quad \mbox{a.s.}\]
\end{rmk}

\begin{rmk}
	An interesting corollary to Theorem \ref{t:as} is obtained by taking $\varphi$ as the identity function. In view of Proposition \ref{cov} and Corollary \ref{corozero}, for a pure $\beta$ARC associated to a dynamical system presenting ACIM $\lambda_T$, {with $U_0\sim\lambda_T$}, a sufficient condition for \eqref{cv} to hold is
	\begin{equation}\label{covmu}
		\sum_{k=1}^\infty \frac{|\mathrm{Cov}(\mu_t,\mu_{t+k})|}{k^q}<\infty,\quad\mbox{for some $0\leq q<1$}.
	\end{equation}
	This result is very convenient since a vast literature concerning the covariance structure of dynamical systems is available.  For instance, it is well known that if the dynamical system is hyperbolic, then the covariance decays exponentially fast \citep[see][]{baladi, HK} and the condition \eqref{covmu} holds for all $q\in[0,1)$. Furthermore, if the system presents long range dependence in the sense that $\mathrm{Cov}(\mu_t,\mu_{t+k})\sim L(k)k^{-b}$, for $0<b<1$, for some slowly varying function $L$, then condition \eqref{covmu} holds, for all $1-b<q<1$. This is the case, for instance, for the Manneville-Pomeau map (Map 4) when $s\in(0.5,1)$. Finally, in this context, \eqref{covmu} is a sufficient condition for a strong law of large number for $Y_t$ to hold.
\end{rmk}


\section{Partial Maximum Likelihood Inference}\label{inf}

Parameter inference in the proposed model can be done via partial maximum likelihood estimation (PMLE). Let $\{(y_t, \bs x_t^\prime)^\prime\}_{t=1}^n$ be a sample from a $\beta$ARC$(p)$ model following \eqref{e:density} and \eqref{e:model} for a given transformation $T_{\bs\theta}$ depending on an identifiable vector of parameters $\bs\theta=(\theta_1,\cdots,\theta_r)^\prime\in\Omega_T\subseteq\R^r$ and $h$ a suitable link function. We shall assume that $u_0\in(0,1)$ is known and such that $T^t(u_0)\notin \{0,1\}$ for all $t$.  Let $\bs\gamma:=(\nu,\alpha,\bs\beta^\prime,\bs\phi^\prime,\bs\theta^\prime)^\prime\in\Omega\subseteq(0,\infty)\times\R^{p+l+1}\times\Omega_T$ be the $(l+p+r+2)$-dimensional vector of parameter related to the model, where $\Omega$ denotes the parameter space. Upon writing
\begin{align*}
	\ell_t(\bs\gamma)=\log\big(f(y_t;\bs\gamma|\F_{t-1})\big)=\log\big(&\Gamma(\nu)\big)-\log\big(\Gamma(\mu_t\nu)\big)-\log\big(\Gamma\big(\nu(1-\mu_t)\big)\big)+\\
&+(\mu_t\nu-1)\log(y_t)+\big(\nu(1-\mu_t)-1\big)\log(1-y_t),
\end{align*}
the log-likelihood associated to model \eqref{e:density} and \eqref{e:model} is given by
\[\ell(\bs\gamma):=\sum_{t=1}^n \ell_t(\bs\gamma).\]
The partial maximum likelihood estimator is then defined as
\begin{align}\label{e:min}
	\widehat{\bs\gamma}=\underset{\bs\gamma\in\Omega}{\mathrm{argmax}}\big\{\ell(\bs\gamma)\big\}.
\end{align}
To obtain the PMLE we need to solve the optimization problem \eqref{e:min}, which can be done upon finding the score function and solving a non-linear system, by using, for instance, the BFGS optimization algorithm. Alternatively, the optimization problem can also be solved by using other methods such as Nelder-Mead.

Since $\eta_t$ is generally a non-linear function of $\bs\theta$, the asymptotic theory of the PMLE in the context of $\beta$ARC models requires some non-trivial adaptations of the existing theory for GARMA-like models \citep[presented, for instance, in][]{Fokianos2004}. {A rigorous large sample theory for the PMLE in the context of $\beta$ARC models is subject of a future paper.} In the next section (and in the supplementary material accompanying the paper), we shall study the finite sample performance of the PMLE in the context of $\beta$ARC models.


\section{Monte Carlo Simulation}\label{sim}

In this section we present a short Monte Carlo simulation study to analyze the finite sample performance of the PMLE in the context of $\beta$ARC models. For the sake of brevity, we shall only consider a single scenario. A more extensive Monte Carlo simulation study considering several different scenarios is presented in the Supplementary material accompanying this paper.

We consider parameter estimation via PMLE for a pure chaotic $\beta$ARC model with map $T_k(x)=(kx)(\mathrm{mod}\,1)$ for $k\in\{3,5,7\}$, considered known, and three different starting points $u_0\in \{0.2 + \pi/100,0.5 + \pi/100,0.8 + \pi/100\}$. We present the results for $\nu=40$ (other cases are presented in the supplementary material). We generate samples $\{y_t\}_{t=1}^n$ for $n\in\{100,500,1000\}$ by setting
\[ \mu_t := T_k^{t-1}(u_0) \quad \mbox{and} \quad y_t \sim \mbox{Beta}\big(\nu\mu_t, \nu(1-\mu_t)\big) \]
For all scenarios we perform  $1,000$ replications. To obtain the PMLE we solve the optimization problem \eqref{e:min}.
The maximization of the objective function was performed by considering the so-called Nelder-Mead algorithm implemented in Fortran by Alan Miller\footnote{available at \color{blue}{https://jblevins.org/mirror/amiller}} and adapted by the authors to handle parameter constraints using the ideas implemented in the matlab function \texttt{fminsearchbnd}\footnote{see \color{blue}{www.mathworks.com/matlabcentral/fileexchange/8277-fminsearchbnd}}. To start the optimization algorithm we calculate the likelihood function for $\nu\in\{5, 50, 100\}$ and select the one with higher likelihood value as starting point.

All computer codes were written by the authors. The most demanding task of parameter estimation was implemented in FORTRAN, while the other tasks were implemented in R \citep{R2018} version 3.6.1. The necessary shared libraries were also compiled in R version 3.6.1.

\subsubsection*{Results}

\begin{figure}[h!]
	\centering
	\includegraphics[width=0.9\textwidth]{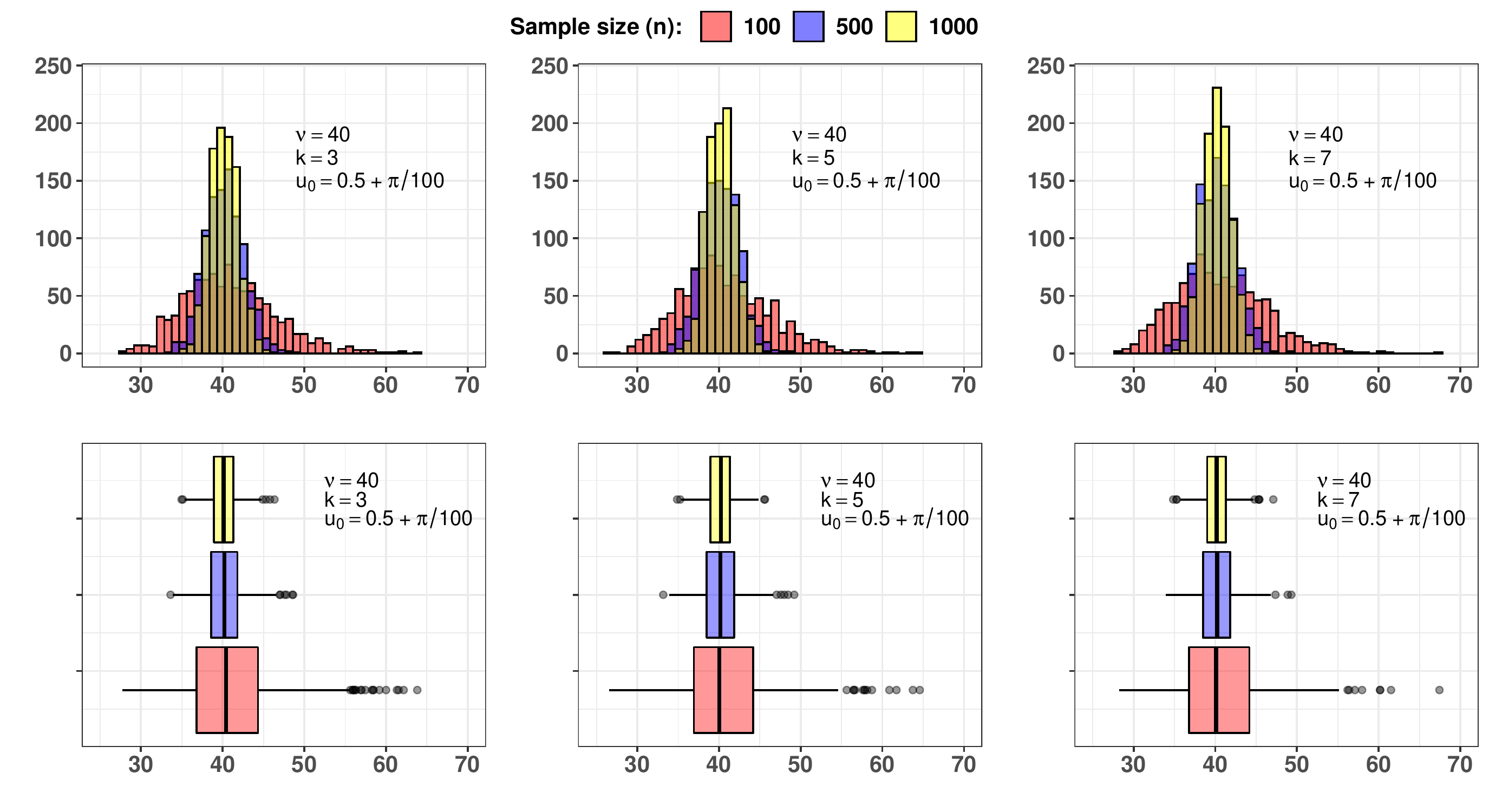}
	\caption{Boxplots of the estimated values, from 1,000 replications, for the parameter $\nu=40$ considering the map $T_k(x)=(kx)(\mathrm{mod}\,1)$ with  $k\in\{3,5,7\}$, $n\in \{100,500,1000\}$ and starting points $u_0\in\{0.2+\pi/100,0.5+\pi/100,0.8+\pi/100\}$.}\label{hkx2}
\end{figure}

Table \ref{t1} presents the simulation results. Highlighted in blue and red are the best and worst scenarios in each case, respectively. We observe that as $n$ increases, the bias and standard deviation of the estimated values decrease. From the results we found no relation between $u_0$ and $k$ with the estimated value of $\nu$. Figure \ref{hkx2} presents the histograms and boxplots of the results for $u_0=0.5+\pi/100$ (the other cases are analogous and can be found in the supplementary material). The histograms suggests that the PMLE in the context of $\beta$ARC models satisfy a central limit theorem. Indeed, applying a Shapiro-Wilk test to the results presented in the top left plot ($k=3$), for $n$ equals 100, 500 and 1,000 the test yields p-values equal to $0.0000$, 0.1077 and 0.2744, respectively.

\begin{table}[h!]
	\renewcommand{\arraystretch}{1.1}
\setlength{\tabcolsep}{3pt}
	\centering
	\caption{Simulation Results for parameter $\nu$ considering the map $T_k(x)=(kx)(\mathrm{mod}\,1)$ with
		$k \in\{3,5,7\}$ and sample size $n\in \{100, 500, 1000\}$: the mean estimated
		value of $\nu$ over 1,000 replications ($\bar{\nu}$), for $\nu=40$,
		the standard deviation of the estimates ($sd_\nu$) and the mean absolute percentage
		error (MAPE).}\label{t1}\vspace{0.4\baselineskip}
	\begin{tabular}{c|ccc|ccc|ccc}
		\hline
		&	 \multicolumn{3}{c|}{$n = 100$} &  \multicolumn{3}{c|}{$n = 500$} &   \multicolumn{3}{c}{$n = 1,000$} \\
		\cline{2-10}
		\multicolumn{1}{c|}{$u_0$} & \multicolumn{1}{c|}{$\bar{\nu}$}	&	\multicolumn{1}{c|}{$sd_{\nu}$}		& \multicolumn{1}{c|}{MAPE}	 &	
		\multicolumn{1}{c|}{$\bar{\nu}$}	&	\multicolumn{1}{c|}{$sd_{\nu}$}		& \multicolumn{1}{c|}{MAPE}	&	
		\multicolumn{1}{c|}{$\bar{\nu}$}	&	\multicolumn{1}{c|}{$sd_{\nu}$}		& \multicolumn{1}{c}{MAPE}\\
		\hline
		\multicolumn{10}{c}{$k = 3 $}\\
		\hline
		$ 0.2 +\frac{\pi}{100}$ &  40.78  &  \color{blue}{\bf 5.4388}  &  \color{blue}{\bf 10.75} &   \color{blue}{\bf 40.18}  &   \color{blue}{\bf 2.3437}  &  4.73 &   \color{blue}{\bf 40.14}  &  1.6885  &  3.42 \\
		$ 0.5 +\frac{\pi}{100}$ &  40.92  &  5.6838  &  11.19 &  40.23  &  2.3867  &  4.78 &  40.15  &  1.7157  &  3.46 \\
		$ 0.8 +\frac{\pi}{100}$ &  40.76  &  5.5986  &  11.00 &  40.30  &  \color{red}{\bf 2.4250} &  \color{red}{\bf 4.90} &  40.19  &  \color{blue}{\bf 1.6813}  &  3.40 \\
		\hline
		\multicolumn{10}{c}{$k = 5 $}\\
		\hline
		$ 0.2 +\frac{\pi}{100}$ &  40.78  &  5.5718  &  10.89 &  \color{red}{\bf 40.47} &  2.3862  &  4.83 &  \color{red}{\bf 40.37}  &  1.6937  &  3.51 \\
		$ 0.5 +\frac{\pi}{100}$ &  40.72  &  5.6021  &  10.94 &  \color{blue}{\bf 40.18} &  2.3819  &  4.77 &  40.24  &  \color{red}{\bf 1.7350} &\color{red}{\bf 3.52} \\
		$ 0.8 +\frac{\pi}{100}$ &  40.94  &  5.4653  &  10.86 &  \color{blue}{\bf 40.18}  &  2.3451  &  \color{blue}{\bf 4.64} &  40.17  &  1.6912  &  3.35 \\
		\hline
		\multicolumn{10}{c}{$k = 7 $}\\
		\hline
		$ 0.2 +\frac{\pi}{100}$ &  \color{red}{\bf 40.99} &  \color{red}{\bf 5.6991}  &  \color{red}{\bf 11.31} &  40.24  &  2.3998  &  4.78 &  40.16  &  1.7052  &  3.46 \\
		$ 0.5 +\frac{\pi}{100}$ &  \color{blue}{\bf 40.64}  &  5.5486  &  11.06 &  40.23  &  2.3665  &  4.77 &  40.21  &  1.7173  &  3.44 \\
		$ 0.8 +\frac{\pi}{100}$ &  40.82  &  5.4884  &  10.95 &  40.20  &  2.3940  &  4.78 &  40.15  &  1.7040  &  \color{blue}{\bf 3.39} \\
		\hline
	\end{tabular}
\end{table}


\section{Real data Application}\label{se:application}

In this section we illustrate the usefulness of the $\beta$ARC model in modeling real data. The variable of interest is the proportion of stocked hydroelectric energy in South Brazil. The data are monthly averages from January 2001 to April 2017 and can be freely downloaded from ONS's (the Brazilian national operator of the electrical system) website ({\color{blue} http://www.ons.org.br}). For comparison with other models, 6 months of data, from November 2016 to April 2017, were reserved for out-of-sample forecasting, yielding a sample of size $n=190$ for fitting purposes. This data was first considered in \cite{Scherer} where the authors fit a $\beta$ARMA(1,1) model to the data and compare its forecasting capabilities with 4 other models: the KARMA(1,1) of \cite{Bayers}, Gaussian ARMA(1,1) and AR(2) models and also with the Holt exponential smoothing algorithm. Our goal is to fit the proposed $\beta$ARC model models and compare to the results reported in \cite{Scherer}.

\begin{figure}[h!]
	\centering
	{\includegraphics[width=0.7\textwidth]{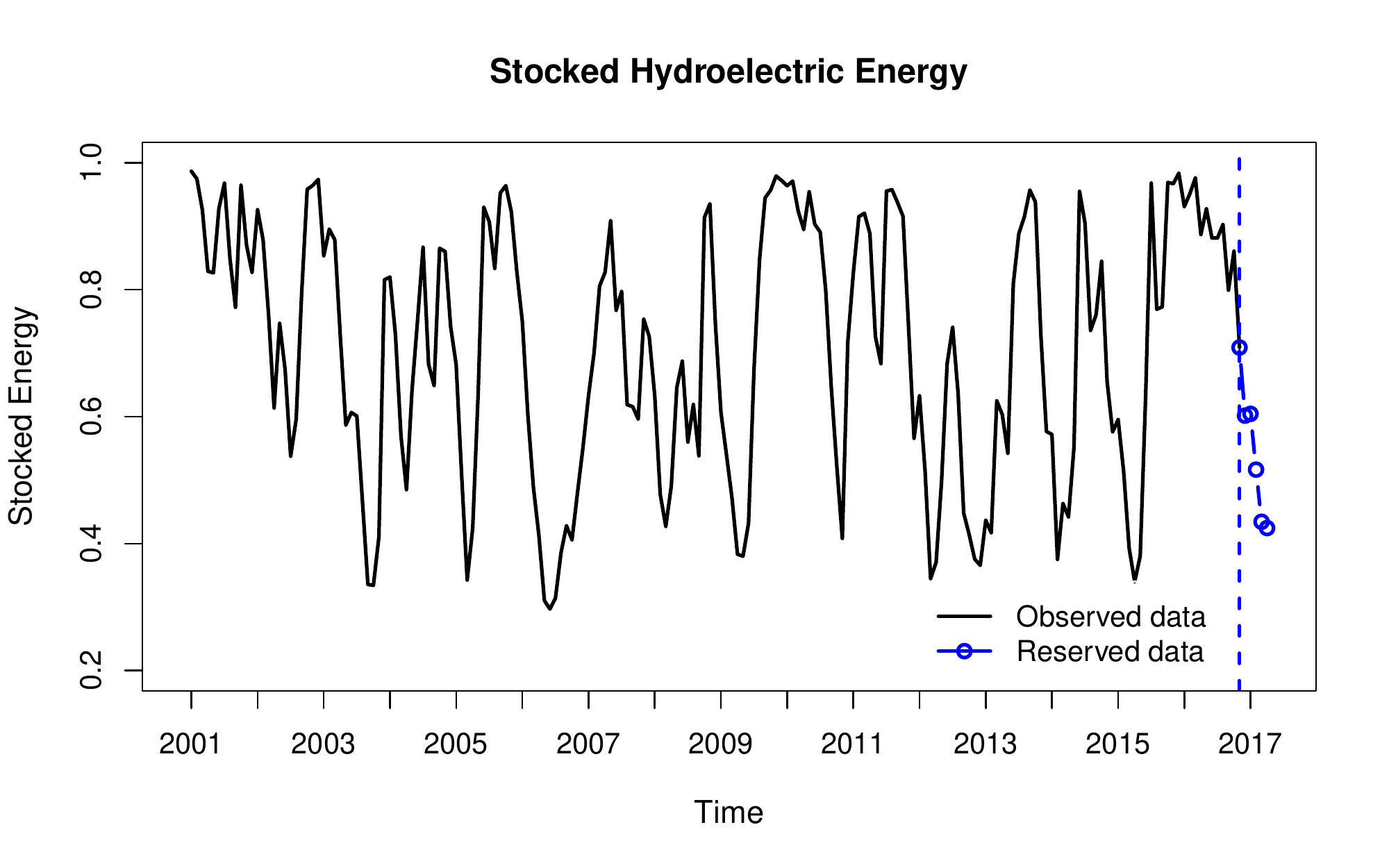}}
	\caption{ Time series plot of the proportion of stocked energy, showing the data used to fit the model (black) and the reserved data (blue). }\label{serie}
\end{figure}

Figure \ref{serie} brings the time series plot of the data. In order to fit a $\beta$ARC model, we apply the Manneville-Pomeau transformation \eqref{e:MP} with $h$ as the identity link. For $g$ in $\eta_t$ we take the cloglog link given by $g(x)=\log(-\log(1-x))$. To obtain the PMLE estimates based on the log-likelihood we first employ a L-BFGS-B optimization with numerical derivatives and then a Nelder-Mead optimization algorithm starting at the values obtained from the L-BFGS-B. This approach showed better results in practice. As for the $p$-values, they are obtained from Wald's $z$ test, based on the numerical hessian (a rigorous asymptotic theory for the PMLE in the context of $\beta$ARC model is under development and shall be presented in another paper).

One delicate computational problem is defining which $u_0$ to apply. Observe that although $s$ is identifiable, $u_0$ is not as the Manneville-Pomeau transformation presents two full branches. However, the sample path of the transformation is identifiable (except for $u_0$), hence, the specific value of $u_0$ brings no useful information in practice, but it is needed to start the PMLE. We overcome this problem with a simple strategy: optimization was performed based on a grid  of 900 initial points for $u_0$ starting at $\pi/1000$ and ending on $1-\pi/1000$. The whole process takes less than 2 minutes in any average computer running Windows 10.

For model selection, we consider only models for which all the coefficients were significant and that the residuals (defined as $y_t-\mu_t$) did not present any serial correlation, condition tested using the Ljung-Box test considering $m=20$ lags. Among all models satisfying these conditions, we chose the one with the smallest in-sample mean absolute prediction error (MAPE-IN), which we shall call Model 1, and also the one with the highest likelihood (Model 2). Under both metrics, a simple $\beta$ARC(1) model satisfied the aforementioned conditions.

\begin{figure}[h!]
	\centering
	\subfigure[]{\includegraphics[width=0.6\textwidth]{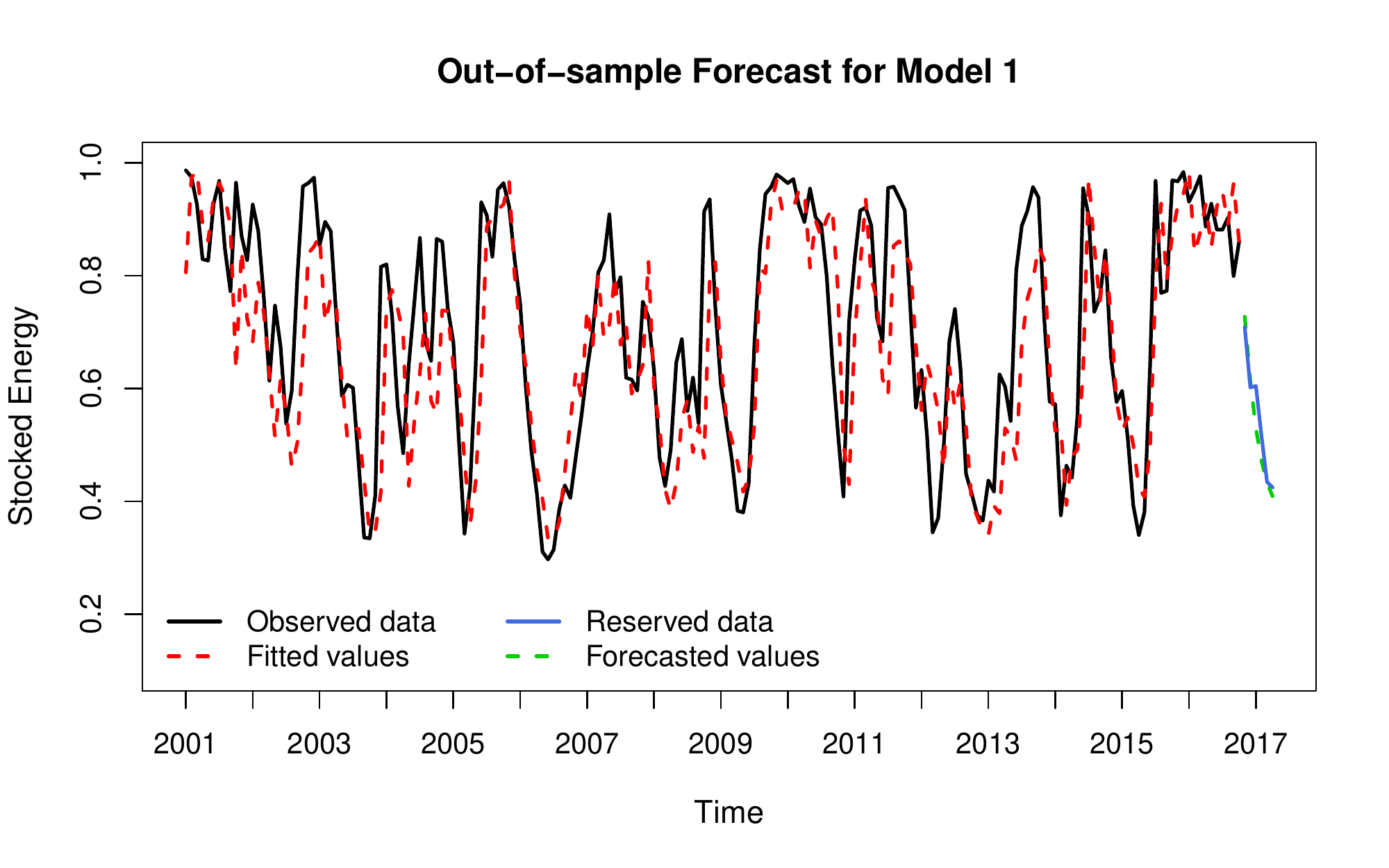}}
	\subfigure[]{\includegraphics[width=0.6\textwidth]{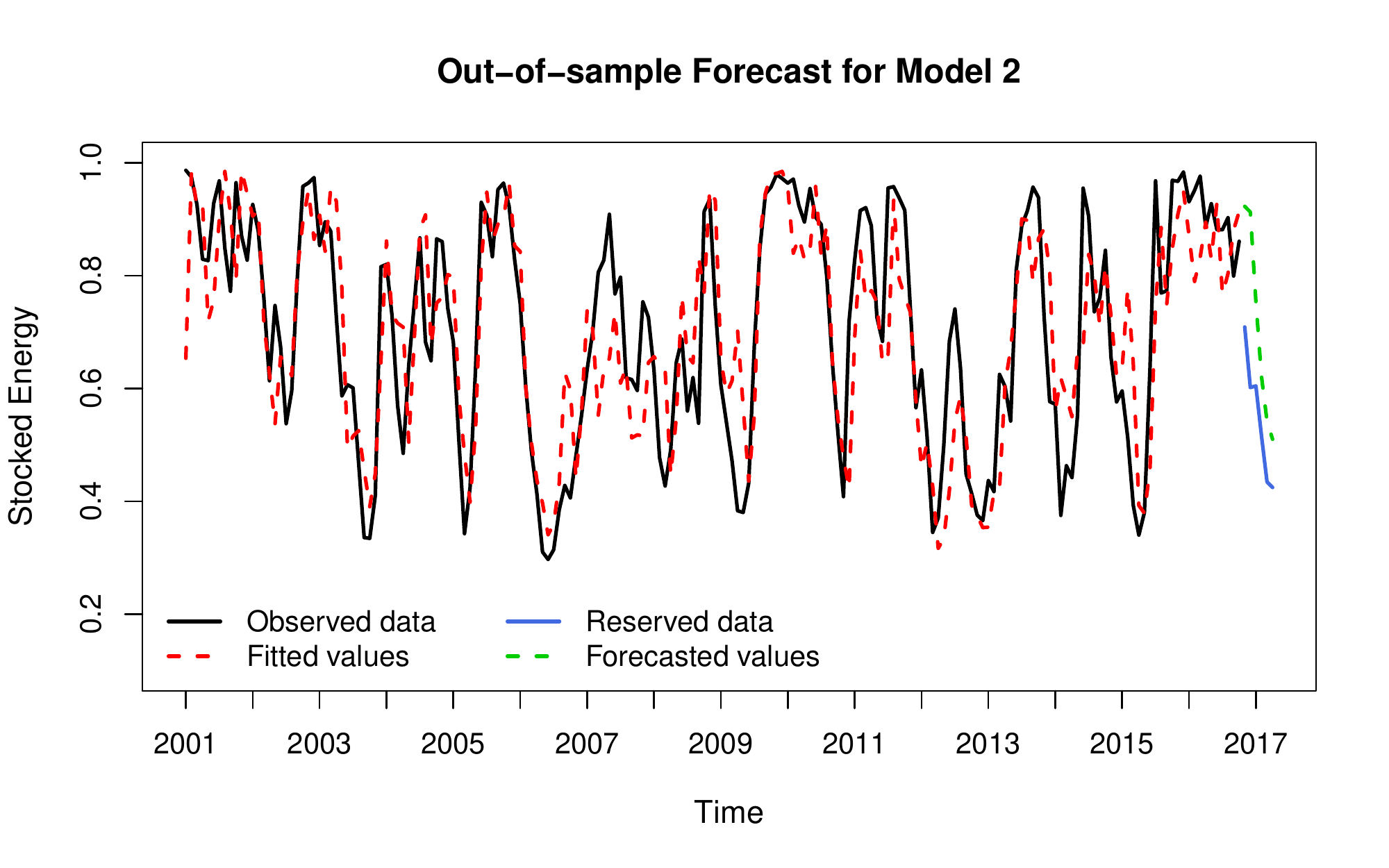}}
	\caption{The observed time series and the in-sample and out-of-sample forecasted values for the fitted $\beta$ARC models 1 (a) and 2 (b).}\label{las}
\end{figure}

Table \ref{t:fitted} presents the fitted $\beta$ARC(1) models while Table \ref{t:forecast}  presents the corresponding in-sample and out-of-sample accuracy measures. Presented are the mean absolute percentage error (MAPE), the mean percentage error (MPE), the average error (ME), the mean absolute error (MAE) and root mean square error (RMSE) for the in-sample and out-of-sample results. We observe that Model 2, obtained via likelihood, presents borderline better in sample accuracy measures than Model 1, except for the MAPE. In terms of out-of sample performance, however, Model 2 outperforms Model 1 in all measures by a large margin. In Figure \ref{las}(b) and (c), we present the observed time series along with in-sample and out-of-sample forecast for models 1 and 2.
\FloatBarrier

\begin{table}   [h!]
\caption{Fitted $\beta$ARC(1) models for the stocked energy data.}\vspace{.3cm}\label{t:fitted}
\centering
		\begin{tabular}{crc|crc}
			\hline
			\multicolumn  {3}{c|}{Model 1: smallest MAPE-IN} & \multicolumn  {3}{c}{Model 2: highest likelihood} \\
			\hline
			&Estimate      & $p$-value    &             &Estimate      & $p$-value   \\
			\hline
			$\alpha$     &  -0.3170      &  0.0000      &$\alpha$     & -0.3653      &  0.0000      \\
			$\phi_1$     &   0.7634      &  0.0000      &$\phi_1$     &  0.7107      &  0.0000      \\
			$s$    		 &   0.8165      &  --          &$s$   		  &  0.3706      &  --          \\
			$\nu$        &   6.3634      &  --          &$\nu$        & 10.5798      &  --          \\
			\hline
			\multicolumn{3}{c|}{$u_0=0.810052910479796$} &                \multicolumn{3}{c}{$u_0=0.423177621111067$}                                           \\
			\multicolumn{3}{c|}{Log-likelihood: $120.01$} &               \multicolumn{3}{c}{Log-likelihood: $134.70$}                                         \\
			\multicolumn{3}{c|}{AIC: $-232.03$ \qquad BIC: $-219.04$} &   \multicolumn{3}{c}{AIC: $-261.40$ \qquad BIC: $-248.41$}                            \\
			\hline
		\end{tabular}
\end{table}

\begin{table}[h!]
\setlength{\tabcolsep}{3pt}
\caption{In and out-of-sample forecasting measures for the two fitted $\beta$ARC(1) models presented in Table \protect\ref{t:fitted}. (M1 and M2 here stand for Model 1 and Model 2.)}\label{t:forecast}\vspace{.4cm}
\centering
		\begin{tabular}{c|ccccc|ccccc}
			\hline
			&\multicolumn{5}{c|}{In-sample accuracy measures}& \multicolumn{5}{c}{Out-of-sample accuracy measures}\\
			\hline
			& MAPE      & MPE        & ME       & MAE       & RMSE      & MAPE      & MPE       & ME        & MAE       & RMSE       \\
			\hline
			M1 & 14.16\%  &  2.34\% & 0.0337   & 0.0957    & 0.1291    & 4.92\%  & 3.40\%    & 0.0177    & 0.0277    & 0.0372     \\
			M2 & 14.67\% & -1.89\% & 0.0069   & 0.0937    & 0.1217    & 28.63\%  & -28.63\%   & -0.1620   & 0.1620    & 0.1801    \\
			\hline
		\end{tabular}
\end{table}

As mentioned before, \cite{Scherer} also considered the same data and modeled it using 5 different models ($\beta$ARMA(1,1), KARMA(1,1), Gaussian ARMA(1,1), Gaussian AR(2) models and the Holt exponential smoothing algorithm). Among these models, the authors report that the $\beta$ARMA(1,1) presented the smallest AIC (-307.9635) and also the best out-of-sample forecasting performance with an MAE of 0.1839 for the same data considered here (other forecasting accuracy measures were not reported). We observe that both fitted $\beta$ARC(1) outperform the fitted $\beta$ARMA(1,1) model in terms of out-of sample performance. Model 1, in special, present out-of-sample MAE of only 0.0277, considerably smaller than the $\beta$ARMA's MAE.


\section{Conclusion}\label{conc}

Here we introduced the Beta Autoregressive Chaotic ($\beta$ARC) processes, a class of dynamic models for time series taking values on the unit interval. The model follows similar structure of other GARMA-like models \citep[in the sense of][]{Benjamin2003}. The random component of the process was modeled through a beta distribution, conditioned on the past information, while the conditional mean was specified allowing the presence of covariates (random and/or deterministic) and an extra additive term  defined by the iteration of a map $T$ defined on $[0,1]$, inspired on the theory of chaotic processes and dynamical systems. This additive term is able to model a wide variety of behaviors in the processes' conditional mean, including short and long range dependence, attracting and/or repelling fixed or periodic points, presence or absence of absolutely continuous invariant measure, among others, allowing for a much broader and flexible dependence structure compared to competitive GARMA-type models presented in the literature.

In the $\beta$ARC model, the extra additive term's definition borrows ideas from dynamical systems. For this reason, a review on the main definitions concerning one dimensional dynamical systems was presented in order to describe the wide variety of behaviors that $T$ can present. Among the main features of the underlying transformation we focused on the existence of  attracting and/or repelling fixed or periodic points and the  presence or absence of absolutely continuous invariant measure. We also discussed how the characteristics of the chaotic process are reflected into the observed time series. In particular, we showed that, as the precision parameter $\nu$ increases, the closer the sample path resembles the conditional mean's dynamics. We also presented some examples where the systematic component can accommodate short or long range dependence, periodic behavior and/or laminar phases.

We also presented some theoretical results which are new in the literature in the sense that are not known for any other GARMA-like process. For instance, we derived the covariance structure of the $\beta$ARC models and obtained sufficient conditions for stationarity, law of large numbers and a Birkhoff-type result to hold. In particular, we showed that, in the absence of an autoregressive component, if $T$ has an absolute continuous $T$-invariant measure and the covariate process is stationary, then the $\beta$ARC processes is stationary.

A short Monte Carlo simulation study to assess the finite sample performance of the PMLE in the context of pure chaotic $\beta$ARC models was presented. The simulation results show small bias and standard deviations which, as expected, decrease as $n$ increases. Histograms of the simulated results also suggest that the PMLE is asymptotically normally distributed in the context of the simulation. A much broader Monte Carlo simulation study is presented in the supplementary material that accompanies the paper.

Finally, an application of the proposed methodology to real data was presented. The variable of interest is the proportion of stored hydroelectrical energy in Southern Brazil from January 2001 through October 2016. Overall the model was capable of fitting the data very well, outperforming competing standard methods in terms of out-of-sample forecasting accuracy.


\subsection*{Acknowledgments}
T.S. Prass gratefully acknowledges the support of FAPERGS (ARD 01/2017, Processo 17/2551-0000826-0).

\bibliographystyle{rss}
\bibliography{prj_chaos}


\newpage


\section*{Supplementary material}

# A Dynamic Model for Double Bounded Time Series With Chaotic Driven Conditional Averages - Supplementary Material

Guilherme Pumi, Taiane Schaedler Prass and Rafael Rigão Souza

Universidade Federal do Rio Grande do Sul

October 5, 2019

This is a supplementary material for the paper *A Dynamic Model for Double Bounded Time Series With Chaotic Driven Conditional Averages*. Here we present figures and tables describing in details the simulation results.

## 1 Statistical Properties of Dynamical Systems

The proposed $\beta$ARC models strongly rely on the dynamic $T_{\boldsymbol{\theta}}$. For instance, the knowledge of the covariance decay of the dynamical system will translate vis-à-vis to the covariance decay of the process itself (see proposition 2.1 in the main paper). The stationarity of the $\beta$ARC process also depends, in some cases, on the stationarity of the dynamical system. The same can be said about the law of large numbers for the $\beta$ARC model. For this reason, in order to explore the richness of behaviors that can be observed under $\beta$ARC models, in this section we introduce some standard definitions and results in discrete dynamical systems, as well as some examples of one dimensional dynamic systems defined on the interval $[0, 1]$. We begin by some topological oriented definitions (i.e. from standard one dimensional dynamical systems) followed by measure preserving ones (from ergodic theory). There are several general and very comprehensive references for ergodic theory, for instance, Hasselblatt and Katok (1996) and Walters (1982), while for topological aspects we refer the reader to Devaney (2003); Robinson (1998).

Given a transformation $T : [0,1] \to [0,1]$ and $x_0 \in (0,1)$, we define the $k$-th iterate of $T$ evaluated on $x_0$ as the $k$-fold composition $x_k := T^k(x_0) = T\big(T^{k-1}(x_0)\big)$. The sequence $\{x_0, x_1, x_2, \cdots\}$ is called the orbit (or sample path) of $T$. A point $x$ is called a fixed point of $T$ if $T(x) = x$ and it is called a periodic point of $T$ with period $s$ (where $s$ is a positive integer) if $T^s(x) = x$ and $T^k(x) \neq x$, for all $k < s$. Any continuous transformation $T : [a,b] \to [a,b]$ has at least one fixed point. Fixed and periodic points can be very different in its nature, as we will see in the following. If $T$ is differentiable at a fixed point $x$ we say that $x$ is attracting if $|T'(x)| < 1$, repelling if $|T'(x)| > 1$, and indifferent (or neutral) if $|T'(x)| = 1$. Indifferent fixed points are important in the $\beta$ARC model because they can be related to slow correlation decay and long range dependence (see map 4 in what follows). Orbits beginning (or passing) near an attracting fixed points converge to this fixed point, while points near a repelling fixed point are pushed away after some iterations of the map, and the subsequent behavior depends on the global properties of the function. Sometimes the sequence of iterates of a repelling fixed point can come back to a neighborhood of this point after some iterations. It is easy to generalize the former definitions concerning fixed point to periodic points: we say that a periodic point $x$ of period $s$ is attracting (respectively repelling or indifferent) if $|(T^s)'(x)| < 1$ (resp. $|(T^s)'(x)| > 1$ or $|(T^s)'(x)| = 1$). In order to illustrate some of the former concepts, Figure 1 shows the graph of the map $T(x) = \theta x(1-x)$, for $\theta = 2.5$ and $3.5$ and the first three elements of the orbit of $x_0 = 0.3$ and $0.7$, respectively. The diagonal $y = x$ is included to help drawing the orbit of

$x_0$. Figure 1(a) shows that the orbit approaches the attracting fixed point $p = 0.6$. The fact that $|T'(p)| < 1$ makes $p$ an attracting fixed point and the negative sign of the derivative at $p$ makes the orbit oscillate around $p$, while converging to $p$. In contrast, Figure 1(b) pictures the behavior of the map near the repelling fixed point.

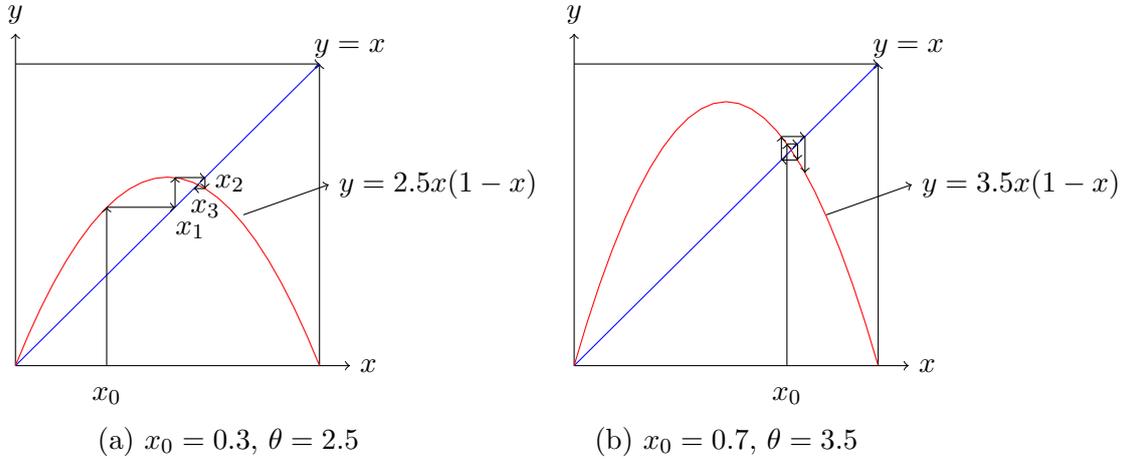

(a) $x_0 = 0.3$, $\theta = 2.5$      (b) $x_0 = 0.7$, $\theta = 3.5$

Figure 1: The first three elements of the orbit of $x_0$ for the map $T(x) = \theta x(1-x)$. In (a) $p = 0.6$ is an attracting fixed point while in (b) $p = 5/7$ is a repelling fixed point.

Now we move our attention to ergodic theory: let $T : [0, 1] \to [0, 1]$ be a Borel measurable transformation. $\lambda_T$ is called a $T$-invariant probability measure (or invariant measure for short) if $\lambda_T$ is a probability measure defined on the Borel sets of $[0, 1]$ which satisfies $\lambda_T\big(T^{-1}(A)\big) = \lambda_T(A)$ for all measurable set $A \subset [0, 1]$. Any continuous function defined on the closed interval $[0, 1]$ (in fact any continuous function defined on a compact set) has at least one invariant measure. If such invariant measure is absolutely continuous with respect to the Lebesgue measure, then we call it an ACIM. ACIM plays an important role in the context of $\beta$ARC models, hence, its existence deserves special attention. Alerted by the fact that there are several cases where a map has no ACIM (for example, if the map has attracting periodic orbits), it is important to consider conditions under which the existence of an ACIM is assured. One of such conditions involves *hyperbolicity*, which, in the one-dimensional case, is equivalent to uniform expansivity: we say that $T$ is uniformly expanding if $T$ is continuously differentiable and there exists $\rho > 1$ such that $|T'(x)| \geq \rho$, for all $x \in (0, 1)$. This kind of maps are the (one-dimensional) *hyperbolic maps*. Any fixed or periodic point of an uniformly expanding map is repelling. Also, if we require the derivative of such maps to be Hölder-continuous, then they present an unique ACIM, which gives positive mass to any open subset and is also ergodic (see Boyarsky and Gora, 1997; Hasselblatt and Katok, 1996; de Melo and Van Strien, 1993, and references therein).

A central result in ergodic theory is the Birkhoff Ergodic Theorem which uses the concept of ergodic measure. An invariant measure $\lambda_T$ is called ergodic if the only measurable sets that are invariant for $T$ are sets of full or zero measure, i.e., if $T^{-1}(A) = A$ implies $\lambda_T(A) = 0$ or $\lambda_T(A) = 1$ (ergodicity implies that it is not possible to split the dynamics into two invariant sets with both having nonzero measure). It can be shown that the set of invariant measures for a given dynamical system is a convex set and the ergodic measures are its extremal points. Therefore, whenever the set of invariant measures for a given map in non-empty, there exists at least one ergodic measure for this map. Birkhoff Ergodic Theorem states that, if $\lambda_T$ is an ergodic invariant measure for $T$, then, for any $\lambda_T$-integrable function $f : [0, 1] \to \mathbb{R}$ and for $\lambda_T$-almost all $x \in [0, 1]$, we have

$$\lim_{n \to \infty} \frac{1}{n} \sum_{k=0}^{n-1} f\big(T^k(x)\big) = \int_0^1 f d\lambda_T.$$

In particular, by taking $f(x) = I_A(x)$, the indicator of $A$, Birkhoff's theorem implies the convergence of the histogram (that is, the sample density) to the associated density (Ding and Zhou, 2009). This application shows why ergodic theory is usually described as the statistical analysis of dynamical systems.

Before passing to some examples of maps we will use in the $\beta$ARC model, we present, for completeness, the traditional definition of a chaotic map. A map is called chaotic if it satisfies the following three conditions: (i) it is sensitive to initial conditions, which means a small change in $x_0$ can lead to a completely different path (orbit); (ii) it is topologically transitive (i.e. there exists $x_0$ such that the orbit of the map beginning at $x_0$ is dense in the domain of the map) and (iii) the set of periodic orbits is dense (any open subset of the domain of the map contains periodic points - if such open subset is small, than usually periodic orbits that intersect the open subset have big periods). The first condition is the most interesting one, and it is usually easy to verify. For instance, the presence of attracting fixed or periodic point implies the map is not chaotic, as none of the three conditions are possible under this hypothesis (we remark the fact that the dense periodic orbits required by the third item above are repelling ones).

Figure 2 exemplifies the usual behavior of a chaotic map: even though the sample paths presented start only 0.01 apart of each other, after less than 10 steps the sample behavior of the processes are completely different.

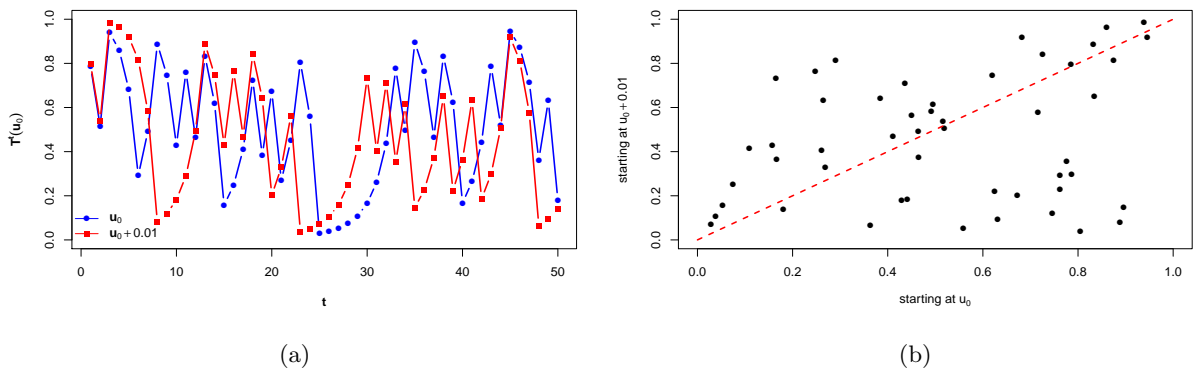

Figure 2: (a) Chaotic behavior of the Mannevile-Pomeau transformation for parameter $s = 0.3$ starting at $u_0 = \pi/4$ (blue) and $u_1 = u_0 + 0.01$ (red) and the associated scatter plot (b).

In the sequel we present 4 examples of maps displaying several interesting behaviors. Maps 1 and 2 are examples of chaotic maps, as well as map 4, while map 3 can be chaotic or not, depending on the parameter $\theta$. We observe, however, that a deep understanding of the three conditions that defines a chaotic map is not required to work with $\beta$ARC model, so that we refrain for further discussing the subject. More details can be found in Devaney (2003).

**Map 1.** One of the simplest dynamical systems is the map defined by

$$T_k(x) = (kx)(\operatorname{mod} 1),$$

where $k$ is a natural number greater or equal to 2. Figure 3(a) presents this function for $k = 3$. It is the basic model of hyperbolic, uniformly expanding map, presenting a simple ACIM: the Lebesgue probability measure. As a consequence of Birkhoff's theorem, there exist a set of total Lebesgue measure in $(0, 1)$ such that, for any point $u_0$ in this set, the orbit $\{T_k^t(u_0)\}_{t \geq 0}$ will be dense on the unit interval, and the histogram of a size $n$ sample from this orbit will converge to the constant function, as $n$ increases. Although there exists periodic points of any period for this maps, all periodic points are repelling and the set of such points is countable, and therefore has zero Lebesgue measure.

3

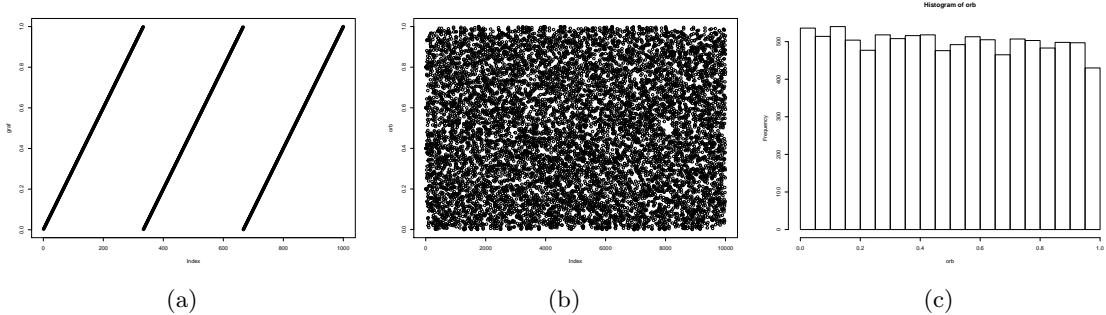

(a)                  (b)                  (c)

Figure 3: (a) The map $T(x) = 3x \pmod 1$ (b) Sample path of this map for $u_0 = \pi/10$: the absence of attracting fixed or periodic points makes possible the chaotic behavior of this map (c) Histogram of the first 10,000 iterates of the map.

**Map 2.** The second example is the map

$$T_\theta(x) = \begin{cases} \frac{x}{\theta} & \text{if } 0 \leq x < \theta, \\ \frac{\theta(x-\theta)}{1-\theta} & \text{if } \theta \leq x \leq 1. \end{cases} \qquad (1)$$

In Figure 4(a) we present a plot from this transformation for $\theta = 0.4$. Depending on the parameter $\theta$ this map can be expanding or not, but its second iterate is always expanding. Therefore, this map also have an ACIM that can be explicitly calculated (see Lopes et al., 1996) and is given by the invariant density $f_\theta(x) = (2-\theta)^{-1}\theta^{-I(x<\theta)}$ which can be recovered in the histogram of the orbit, for large samples (see Figure 4(b) and (c)).

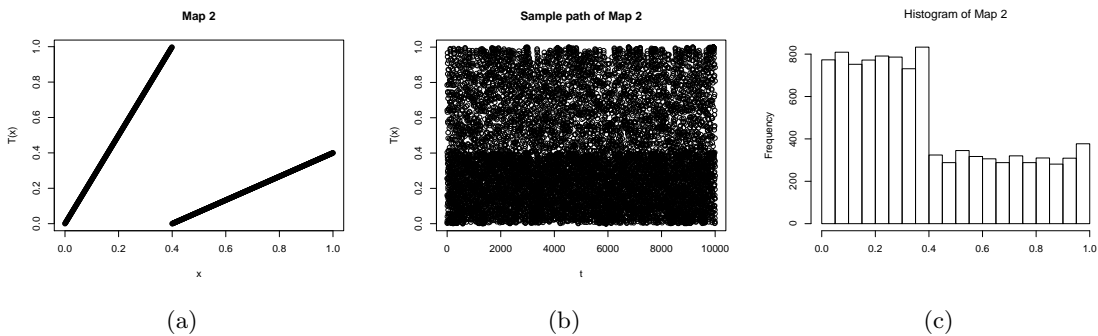

(a)                  (b)                  (c)

Figure 4: (a) A plot of the map given by (1) for $\theta = 0.4$. (b) Sample path of this map for $u_0 = \pi/4$. (c) Histogram of the first 10,000 iterates of the map.

**Map 3.** For $0 \leq \theta \leq 4$, the Logistic map[1] is given by

$$T_\theta(x) = \theta x(1-x).$$

This map is not expanding: in fact it has derivative zero on $x = 1/2$. Therefore, its analysis is more complicated than the preceding examples: we can observe different behaviors for the Logistic map, depending on the parameter $\theta$. For $\theta < 3$ the Logistic map has an attracting fixed point, as can be seen in Figure 1. The most interesting behaviors occurs for $\theta > 3$. As it can be seen in Figure 5, the fixed points are both repellors, as the derivative has modulus greater than 1. For $\theta = 10/3$ the logistic map has an attracting periodic orbit of period 2, shown in Figure 6(a). For $\theta = 4$ the Logistic map has no attracting periodic points, which allows very complicated sample paths (Figure 6(b)), and has an ACIM which can be seen in Figure 6(c).

---
[1] also called *quadratic family*

Figure 7(a)-(i) present some sample paths of the Logistic Map for $\theta$ going from 3.55 to 3.95. For more information on the Logistic map, see Devaney (2003) and Robinson (1998).

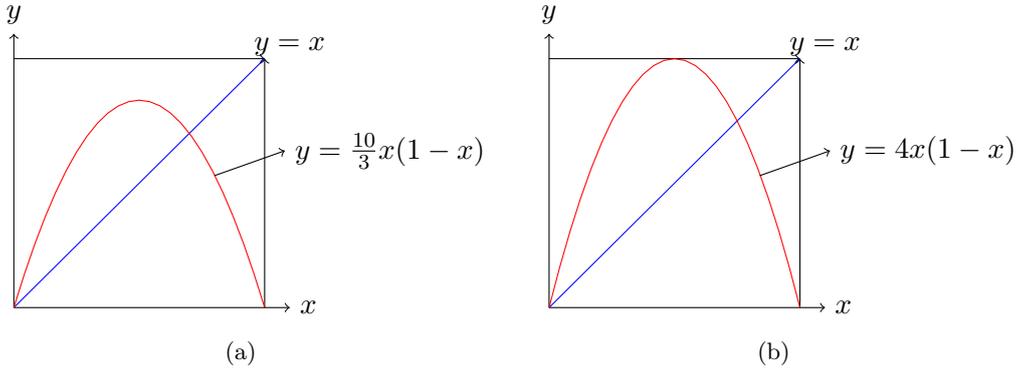

Figure 5: (a) The logistic map for $\theta = 10/3$. (b) The logistic map for $\theta = 4$.

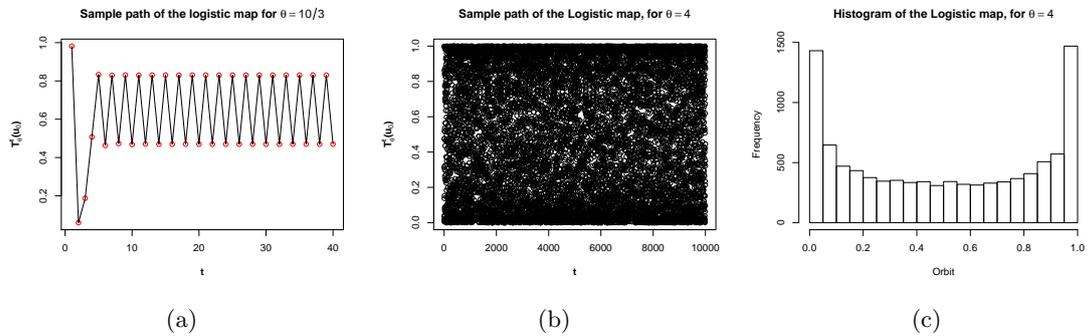

Figure 6: (a) Plot of a sample path associated to the logistic map for $\theta = 10/3$ starting at $u_0 = \pi/3.2$ showing an attracting periodic orbit of period 2; (b) sample path of the logistic map for $\theta = 4$; (c) Histogram of the sample path in (b).

**Map 4.** Now we introduce a family of maps presenting and indifferent fixed point on $x = 0$, that may or may not have an absolute continuous invariant probability measure depending on the parameter chosen. For $s > 0$, the Manneville-Pomeau transformation $T_s : [0, 1] \to [0, 1]$, is given by

$$T_s(x) = (x + x^{1+s})(\text{mod } 1). \tag{2}$$

For $s \in (0, 1)$, there exists an absolutely continuous $T_s$-invariant probability measure (Thaler, 1980), which can be seen in Figure 8(c). For $s \geq 1$ there exists an absolutely continuous invariant measure which is only $\sigma$-finite (not a probability measure). Figure 8(a) show the Manneville-Pomeau transformation for $s = 0.75$. The Manneville-Pomeau transformation presents a property referred to as transition to turbulence through intermittency (Eckmann, 1981). It also has an indifferent fixed point at 0 and, hence, it is not uniformly expanding. The chaotic processes associated to $T_s$ are often called Manneville-Pomeau processes which present a very slow correlation decay when $s \in (0.5, 1)$, characteristic of long range dependent processes and it is commonly viewed as an alternative model for long range dependence outside the classical duet of Fractional Brownian Motion and ARFIMA processes. This slow decay is mainly due to the presence of laminar behavior near zero, which can be seen in Figure 8(b).

5

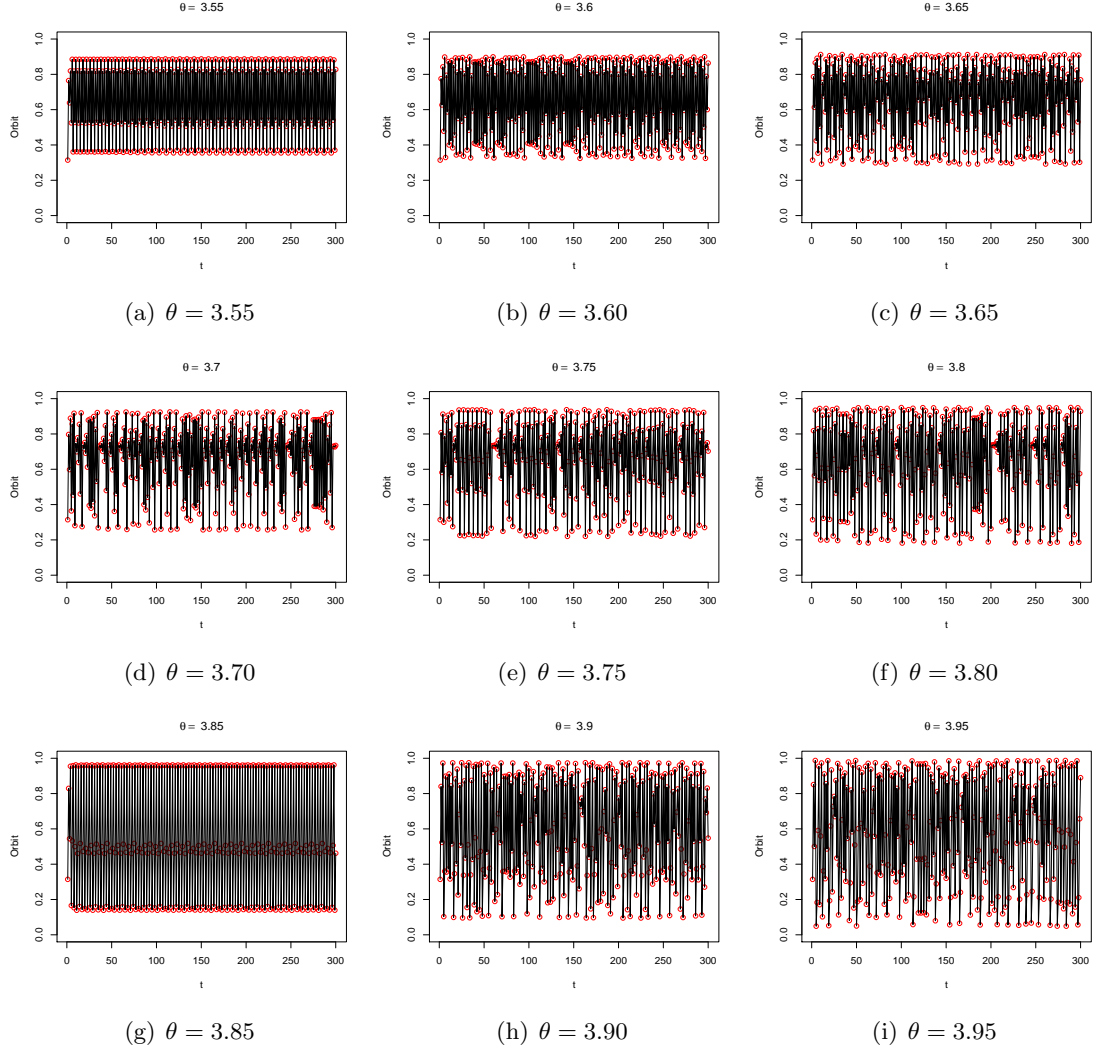

Figure 7: Sample paths of the logistic map for $\theta$ from 3.55 to 3.95 with $u_0 = \pi/10$.

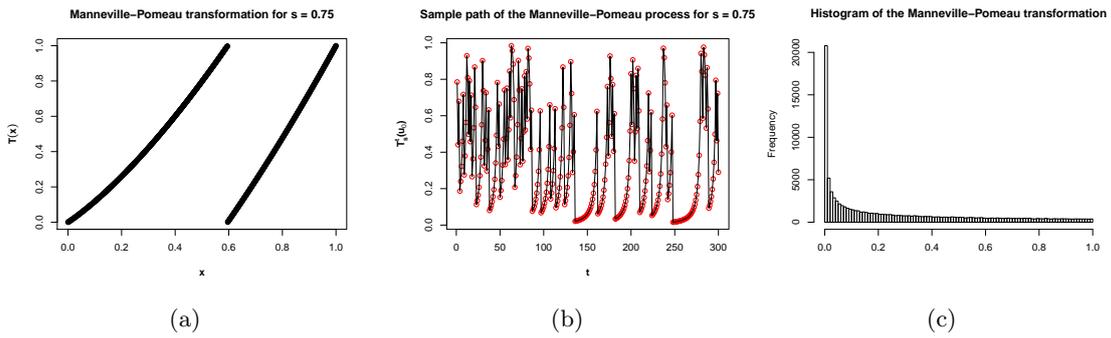

Figure 8: (a) A plot of the Manneville-Pomeau transformation (2) for $s = 0.75$. (b) Sample path of this map for $u_0 = \pi/4$ showing laminar behavior near zero. (c) Histogram of the first 100,000 iterates of the map.

## 2 Monte Carlo Simulation - Pure chaotic $\beta$ARC case

A Monte Carlo simulation study was performed to analyze the finite sample performance of the PMLE on the context of pure chaotic $\beta$ARC models considering two different type of maps:

map 1 and map 3 (logistic map). For map 1 we assume $k$ fixed and known while for map 3 we assume both, $\theta$ known and $\theta$ unknown. Different values of $\nu$ were considered to generate models with small and large conditional variance. We also investigate the influence of $k$, $\theta$ and the sample size $n$. All computer codes were written by the authors. The most demanding task of parameter estimation was implemented in FORTRAN, while the other tasks were implemented in R (R Core Team, 2018) version 3.6.1. The necessary shared libraries were also compiled in R version 3.6.1.

**Data Generating Process**

To generate the samples from the pure chaotic $\beta$ARC process the following was set.

1) For map 1 we consider $k \in \{3, 5, 7\}$ and for map 3 we consider $\theta \in \{1, 2, 3, 3.3, 3.6, 3.99\}$.

2) Three different values of $\nu$ were considered, namely, $\nu \in \{10, 40, 120\}$.

3) Three different values of $u_0$ were considered, namely, $u_0 \in \{0.2 + \pi/100, 0.5 + \pi/100, 0.8 + \pi/100\}$.

4) The time series $\{y_t\}_{t=1}^n$ was generated as follows

$$\mu_t := T^{t-1}(u_0) \quad \text{and} \quad y_t \sim \text{Beta}(\nu\mu_t, \nu(1-\mu_t))$$

were $T$ denotes map 1 or 3 and. Observe that $y_t$ will follow the beta distribution parameterized as in (2) in the paper.

5) All time series were generated with sample size $n = 1,000$.

6) For all scenarios we perform $1,000$ replications.

**Parameter Estimation**

To obtain the PMLE we numerically solve the optimization problem

$$\widehat{\boldsymbol{\gamma}} = \operatorname*{argmax}_{\boldsymbol{\gamma} \in \Omega} \{\ell(\boldsymbol{\gamma})\}.$$

The following was set.

1) The maximization of the objective function was performed by considering the so-called Nelder-Mead algorithm implemented in Fortran by Alan Miller[2] and adapted by the authors to handle parameter constraints using the ideas implemented in the matlab function `fminsearchbnd`[3].

2) In all scenarios, $u_0$ is assumed known.

3) To start the optimization algorithm we calculate the log-likelihood function in a grid of initial points and select as starting point the one with higher log-likelihood value. For $\nu$ we calculate the likelihood in the set $\nu \in \{5, 50, 100\}$. For the logistic map, when $\theta$ was assumed unknown, we calculate the likelihood in the set $\theta \in \{1.1, 2.1, 3.1, 3.5, 3.8\}$.

4) For all scenarios, estimation was performed by considering samples $\{y_t\}_{t=1}^n$ with $n \in \{100, 500, 1000\}$ from the simulated time series.

---

[2]available at https://jblevins.org/mirror/amiller
[3]see www.mathworks.com/matlabcentral/fileexchange/8277-fminsearchbnd

## 2.1 Simulation Results for the pure chaotic $\beta$ARC case

Overall, the simulation results show that as $n$ increases, the estimator's bias and variance decrease so that the estimated values are close to the true ones. As expected, we found no relation between $u_0$ and the estimator's performance. Below are some more specific comments.

1) Results for map 1 are presented in Table 1. For $n$ fixed we found no relation between $k$, $\nu$ and $u_0$ and the estimation performance. As expected, the standard deviation decreases as $n$ increases. From $n = 100$ to $n = 500$, the average reduction in bias is about $1.5\%$ and from $n = 500$ to $n = 1,000$ is about $0.14\%$ on average, considering all $\nu$ and $u_0$. As for the average reduction in standard deviation, from $n = 100$ to $n = 500$ and from $n = 500$ to $n = 1,000$ we have reductions of about $57.1\%$ and $28.5\%$, respectively. The average reduction in MAPE, from $n = 100$ to $n = 500$ and from $n = 500$ to $n = 1,000$ is about $56.2\%$ and $28.0\%$, respectively. Figures 9, 10 and 11 present the histograms of the estimated values of $\nu$ for each $u_0$, $k$ and sample sizes $n$. It is interesting to observe the evolution of the histograms as the sample size increase. To exemplify, in the top left panel of Figure 9, we show the histograms for $\nu = 10$, $k = 3$ and $u_0 = 0.2 + \pi/100$ and $n \in \{100, 500, 1000\}$. Applying a Shapiro-Wilk test for $n$ equal 100, 500 and 1,000 we obtain p-values equal to $1.5 \times 10^{-9}$, $0.0037$ and $0.0874$, respectively, suggesting that the PMLE satisfy a central limit theorem in the context of this simulation. The other cases are similar. Figure 12, 13 and 14 present boxplots of the simulation results. We found that there is little difference among the histograms and boxplot for different $u_0$, except for the scale (related to $\nu$). For each fixed $n$ and $\nu$, the histograms and boxplots for different values of $u_0$ and $k$ are similar.

2) Results for the logistic map assuming $\theta$ known are presented in Table 2. We found no relation between $\theta$, $\nu$ and $u_0$ and the estimation performance, for fixed $n$. Except in 8 cases (out of 162), the bias always decreases as $n$ increases. The average reduction in standard deviation, from $n = 100$ to $n = 500$ and from $n = 500$ to $n = 1,000$ are about $56.7\%$ and $28.2\%$, respectively. The average reduction in MAPE, from $n = 100$ to $n = 500$ and from $n = 500$ to $n = 1,000$ are about $55.7\%$ and $28.0\%$, respectively. Notice that the reduction in standard deviation and MAPE for the logistic map is very close to the ones reported for map 1. Figures 15, 16 and 17 present the histograms of the estimated values of $\nu$, $\theta$ and $n$ for $u_0 = 0.5 + \pi/100$. For other $u_0$ the results are similar and are not shown. Boxplots are presented in Figures 18, 19 and 20. The histograms and boxplot present similar behavior as those for map 1.

3) Notice that the reduction in standard deviation and MAPE in items 1 and 2 are close to what would be expected if the PMLE were asymptotically normally distributed.

4) For the logistic map with $\theta$ unknown, the simulation results are presented in Table 3 (for $\nu$) and 4 (for $\theta$) and overall the estimation of $\nu$ and $\theta$ are very good. Again we found no relation between $\theta$, $\nu$ and $u_0$ and the estimation performance, for fixed $n$. The estimation of $\theta$ is much more precise than that of $\nu$, with very small standard deviation and bias. For $\theta = 1$, the estimation of $\nu$ always presented the highest standard deviation. Boxplots and histograms are presented together in Figures 21 (for $\nu = 10$), 22 (for $\nu = 40$) and 23 (for $\nu = 120$). We show the results only for $u_0 = 0.5 + \pi/100$ as the other cases are similar. The boxplot for the different $\nu$ are very similar to each other. Especially for $\theta = 1$, the scatter plots are highly affected by some outliers. Despite the outliers in the estimation, the scatter plots and histograms suggests that, in this case, the PMLE for $(\nu, \theta)$ is asymptotically jointly normally distributed.

## Simulation Results for Map 1

Table 1: Simulation Results for parameter $\nu$ considering Map 1 with $k \in \{3, 5, 7\}$ and sample size $n \in \{100, 500, 1000\}$: the mean estimated value of $\nu$ over 1,000 replications ($\bar{\nu}$), for $\nu \in \{10, 40, 120\}$, the standard deviation of the estimates ($sd_\nu$) and the mean absolute percentage error (MAPE).

| $k$ | $u_0$ | $n = 100$ | | | $n = 500$ | | | $n = 1,000$ | | |
|---|---|---|---|---|---|---|---|---|---|---|
| | | $\bar{\nu}$ | $sd_\nu$ | MAPE | $\bar{\nu}$ | $sd_\nu$ | MAPE | $\bar{\nu}$ | $sd_\nu$ | MAPE |
| | | | | | $\nu = 10$ | | | | | |
| 3 | $0.2 + \frac{\pi}{100}$ | 10.30 | 1.3247 | **10.62** | **10.07** | **0.5595** | **4.44** | **10.04** | **0.4008** | **3.21** |
| | $0.5 + \frac{\pi}{100}$ | 10.21 | **1.2975** | 10.28 | 10.12 | 0.5707 | 4.63 | 10.08 | 0.4084 | 3.27 |
| | $0.8 + \frac{\pi}{100}$ | 10.26 | 1.2999 | **10.13** | 10.14 | 0.5709 | 4.63 | 10.10 | 0.4085 | 3.36 |
| 5 | $0.2 + \frac{\pi}{100}$ | 10.19 | 1.3252 | 10.35 | **10.15** | 0.5715 | **4.71** | **10.17** | **0.4159** | **3.57** |
| | $0.5 + \frac{\pi}{100}$ | 10.18 | 1.3004 | 10.27 | 10.10 | 0.5689 | 4.61 | 10.15 | 0.4116 | 3.49 |
| | $0.8 + \frac{\pi}{100}$ | **10.43** | 1.3129 | 10.56 | **10.15** | 0.5712 | 4.66 | 10.13 | 0.4029 | 3.39 |
| 7 | $0.2 + \frac{\pi}{100}$ | **10.16** | **1.3330** | 10.57 | 10.11 | 0.5624 | 4.53 | 10.08 | 0.4059 | 3.31 |
| | $0.5 + \frac{\pi}{100}$ | 10.22 | 1.3052 | 10.37 | 10.14 | 0.5713 | 4.69 | 10.12 | 0.4028 | 3.31 |
| | $0.8 + \frac{\pi}{100}$ | 10.35 | 1.3224 | 10.46 | 10.11 | **0.5782** | 4.65 | 10.09 | 0.4111 | 3.35 |
| | | | | | $\nu = 40$ | | | | | |
| 3 | $0.2 + \frac{\pi}{100}$ | 40.78 | **5.4388** | **10.75** | **40.18** | **2.3437** | 4.73 | **40.14** | 1.6885 | 3.42 |
| | $0.5 + \frac{\pi}{100}$ | 40.92 | 5.6838 | 11.19 | 40.23 | 2.3867 | 4.78 | 40.15 | 1.7157 | 3.46 |
| | $0.8 + \frac{\pi}{100}$ | 40.76 | 5.5986 | 11.00 | 40.30 | **2.4250** | **4.90** | 40.19 | **1.6813** | 3.40 |
| 5 | $0.2 + \frac{\pi}{100}$ | 40.78 | 5.5718 | 10.89 | **40.47** | 2.3862 | 4.83 | **40.37** | 1.6937 | 3.51 |
| | $0.5 + \frac{\pi}{100}$ | 40.72 | 5.6021 | 10.94 | **40.18** | 2.3819 | 4.77 | 40.24 | **1.7350** | **3.52** |
| | $0.8 + \frac{\pi}{100}$ | 40.94 | 5.4653 | 10.86 | **40.18** | 2.3451 | **4.64** | 40.17 | 1.6912 | 3.35 |
| 7 | $0.2 + \frac{\pi}{100}$ | **40.99** | **5.6991** | **11.31** | 40.24 | 2.3998 | 4.78 | 40.16 | 1.7052 | 3.46 |
| | $0.5 + \frac{\pi}{100}$ | **40.64** | 5.5486 | 11.06 | 40.23 | 2.3665 | 4.77 | 40.21 | 1.7173 | 3.44 |
| | $0.8 + \frac{\pi}{100}$ | 40.82 | 5.4884 | 10.95 | 40.20 | 2.3940 | 4.78 | 40.15 | 1.7040 | **3.39** |
| | | | | | $\nu = 120$ | | | | | |
| 3 | $0.2 + \frac{\pi}{100}$ | 122.64 | 17.1119 | 11.23 | 120.56 | 7.3838 | 4.96 | 120.35 | 5.1672 | 3.47 |
| | $0.5 + \frac{\pi}{100}$ | 122.55 | 17.3253 | 11.47 | **120.47** | 7.3247 | 4.91 | 120.33 | 5.2495 | 3.49 |
| | $0.8 + \frac{\pi}{100}$ | 122.59 | **17.6640** | **11.50** | 120.63 | 7.3055 | 4.90 | 120.38 | 5.3209 | 3.58 |
| 5 | $0.2 + \frac{\pi}{100}$ | 122.66 | 17.0043 | 11.26 | **121.24** | **7.4303** | **5.02** | **120.73** | 5.2244 | 3.54 |
| | $0.5 + \frac{\pi}{100}$ | **122.82** | 17.2089 | 11.35 | 120.56 | 7.3794 | 5.01 | 120.57 | **5.3540** | **3.62** |
| | $0.8 + \frac{\pi}{100}$ | 122.57 | 17.0216 | 11.29 | 120.54 | 7.2335 | 4.85 | 120.41 | 5.1375 | 3.45 |
| 7 | $0.2 + \frac{\pi}{100}$ | 122.85 | 17.3146 | 11.33 | 120.55 | **7.3723** | 4.95 | 120.37 | 5.2457 | 3.51 |
| | $0.5 + \frac{\pi}{100}$ | **122.43** | 17.2454 | 11.27 | 120.53 | 7.3290 | 4.90 | 120.43 | 5.2120 | 3.49 |
| | $0.8 + \frac{\pi}{100}$ | 122.66 | **16.8971** | **11.05** | **120.47** | 7.2517 | **4.83** | **120.28** | **5.1218** | **3.42** |

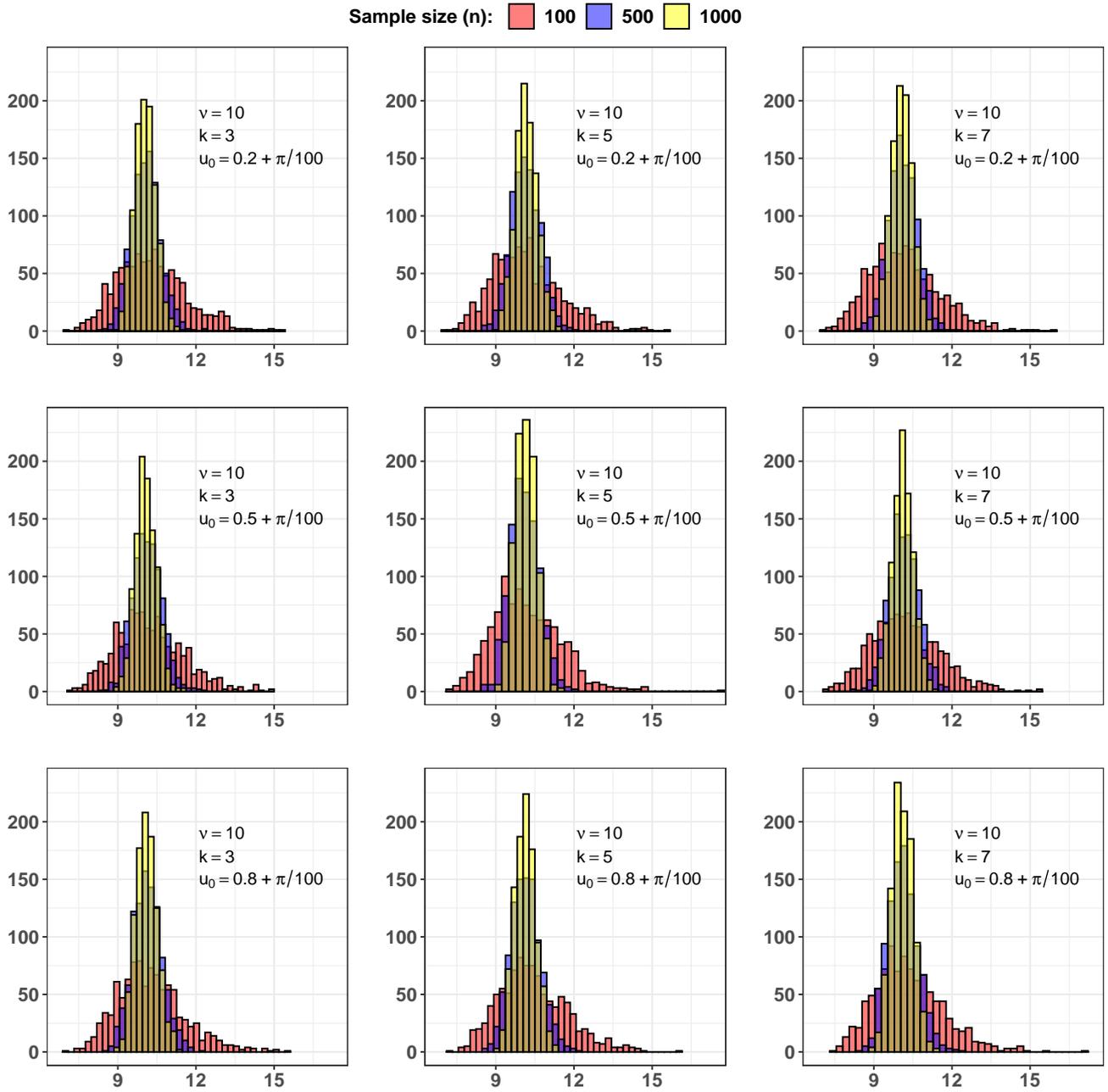

Figure 9: Histogram of the estimated values, from 1,000 replications, for the parameter $\nu = 10$ considering Map 1 with $k \in \{3, 5, 7\}$ (first, second and third column of each subfigure, respectively), $n \in \{100, 500, 1000\}$ and starting points $u_0 \in \{0.2+\pi/100, 0.5+\pi/100, 0.8+\pi/100\}$ (first, second and third row of each subfigure, respectively).

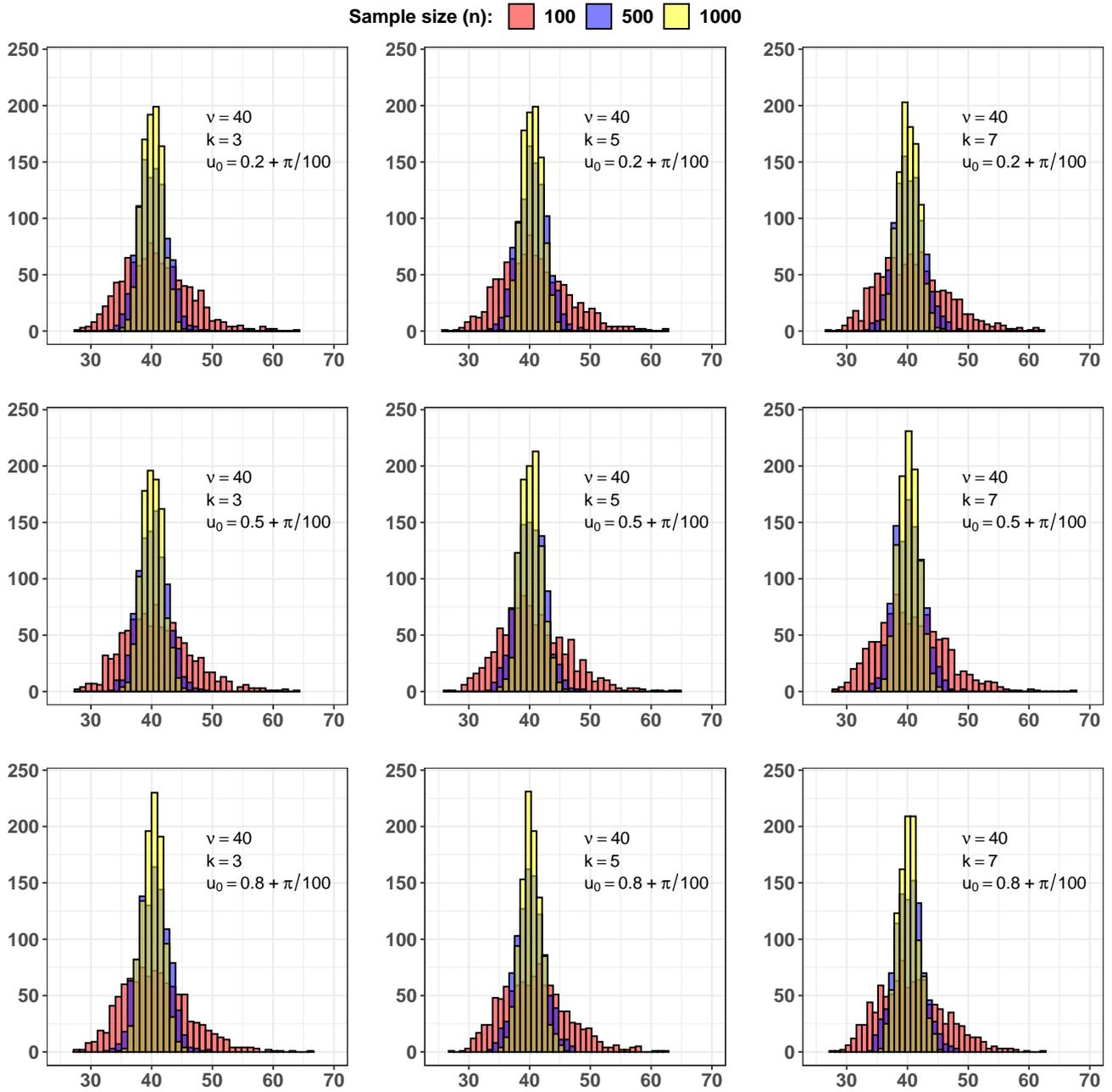

Figure 10: Histogram of the estimated values, from 1,000 replications, for the parameter $\nu = 40$ considering Map 1 with $k \in \{3, 5, 7\}$ (first, second and third column of each subfigure, respectively), $n \in \{100, 500, 1000\}$ and starting points $u_0 \in \{0.2+\pi/100, 0.5+\pi/100, 0.8+\pi/100\}$ (first, second and third row of each subfigure, respectively).

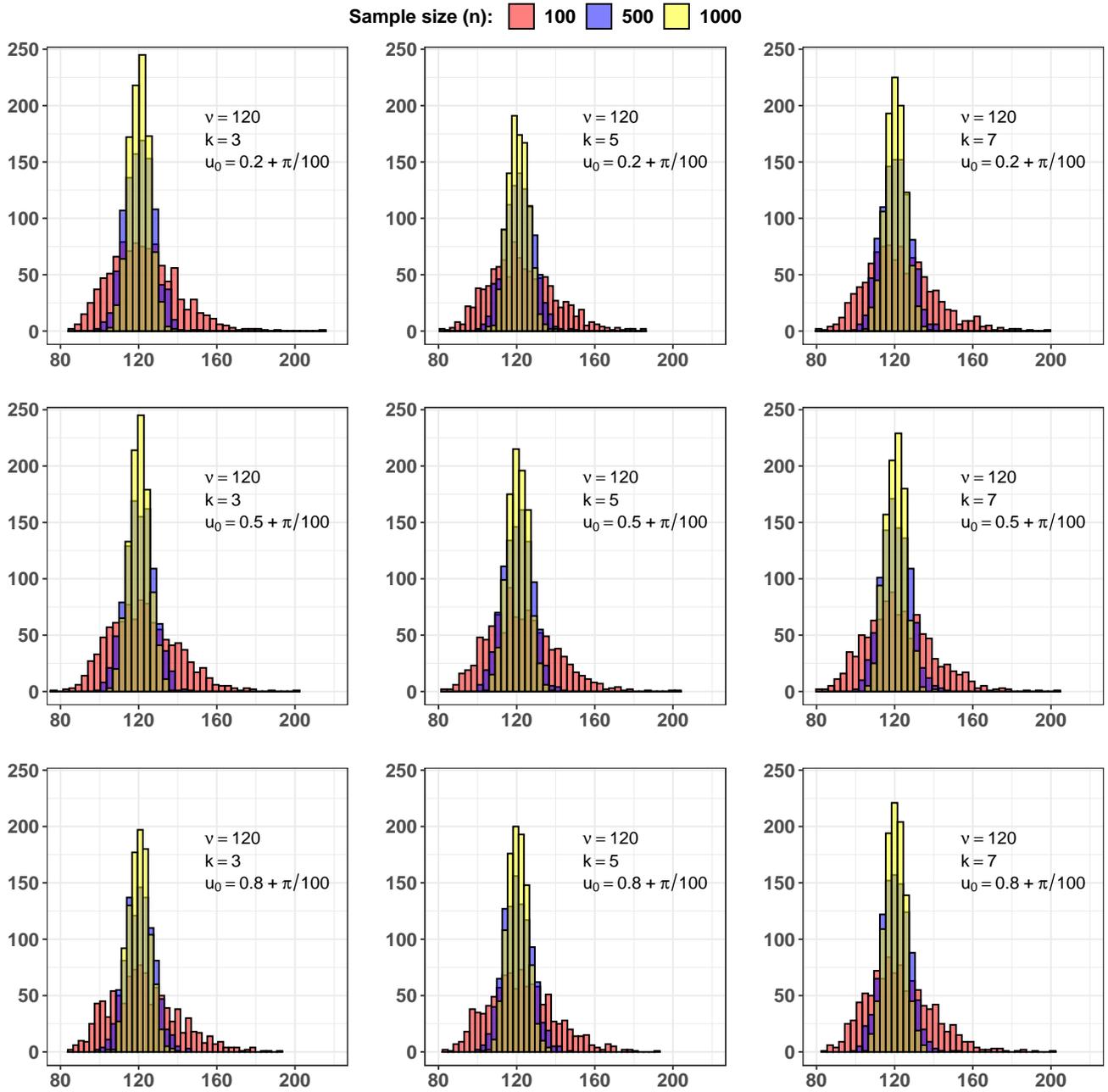

Figure 11: Histogram of the estimated values, from 1,000 replications, for the parameter $\nu = 120$ considering Map 1 with $k \in \{3, 5, 7\}$ (first, second and third column of each subfigure, respectively), $n \in \{100, 500, 1000\}$ and starting points $u_0 \in \{0.2+\pi/100, 0.5+\pi/100, 0.8+\pi/100\}$ (first, second and third row of each subfigure, respectively).

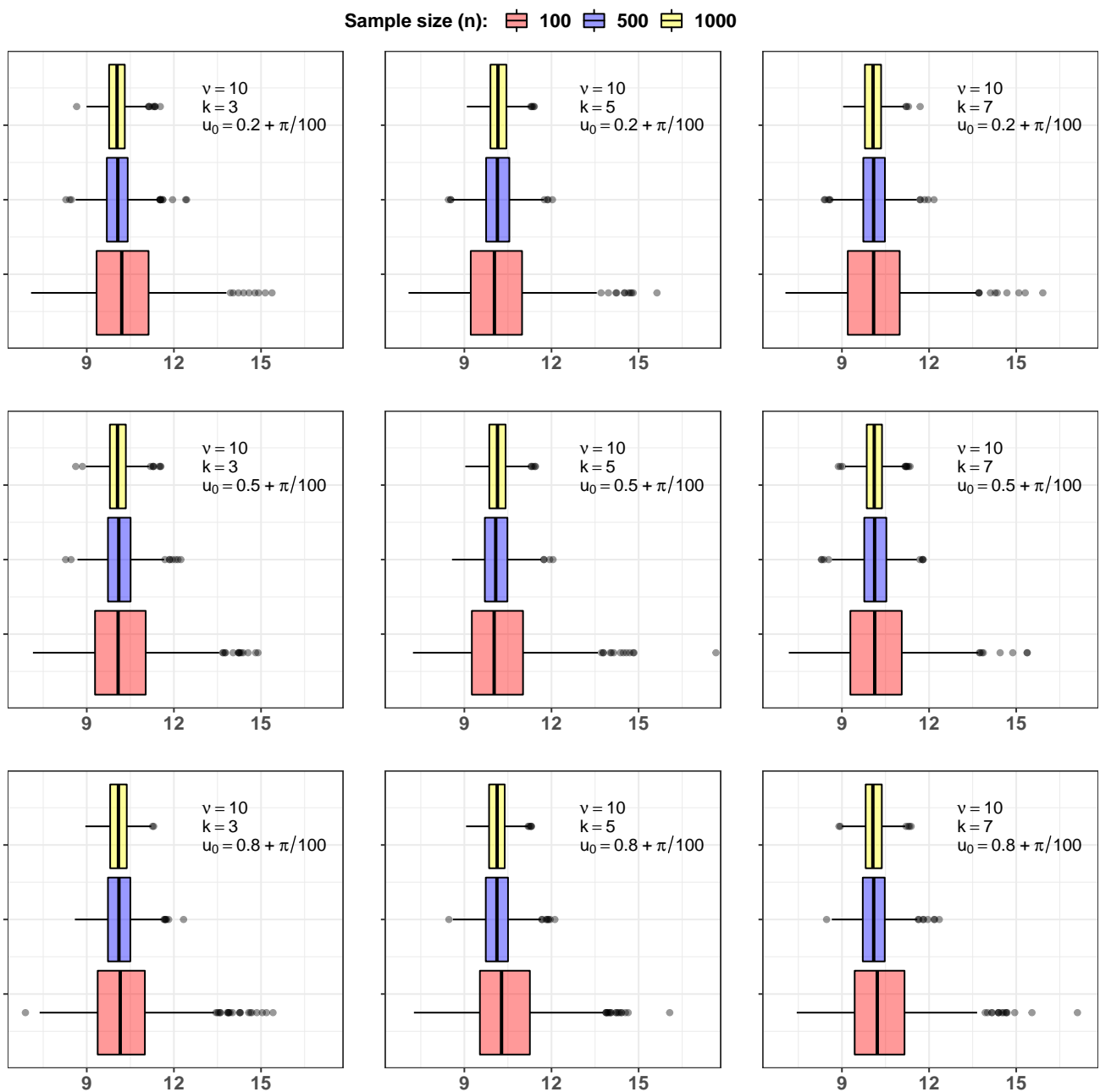

Figure 12: Boxplots of the estimated values, from 1,000 replications, for the parameter $\nu = 10$ considering Map 1 with $k \in \{3, 5, 7\}$ (first, second and third column of each subfigure, respectively), $n \in \{100, 500, 1000\}$ and starting points $u_0 \in \{0.2+\pi/100, 0.5+\pi/100, 0.8+\pi/100\}$ (first, second and third row of each subfigure, respectively).

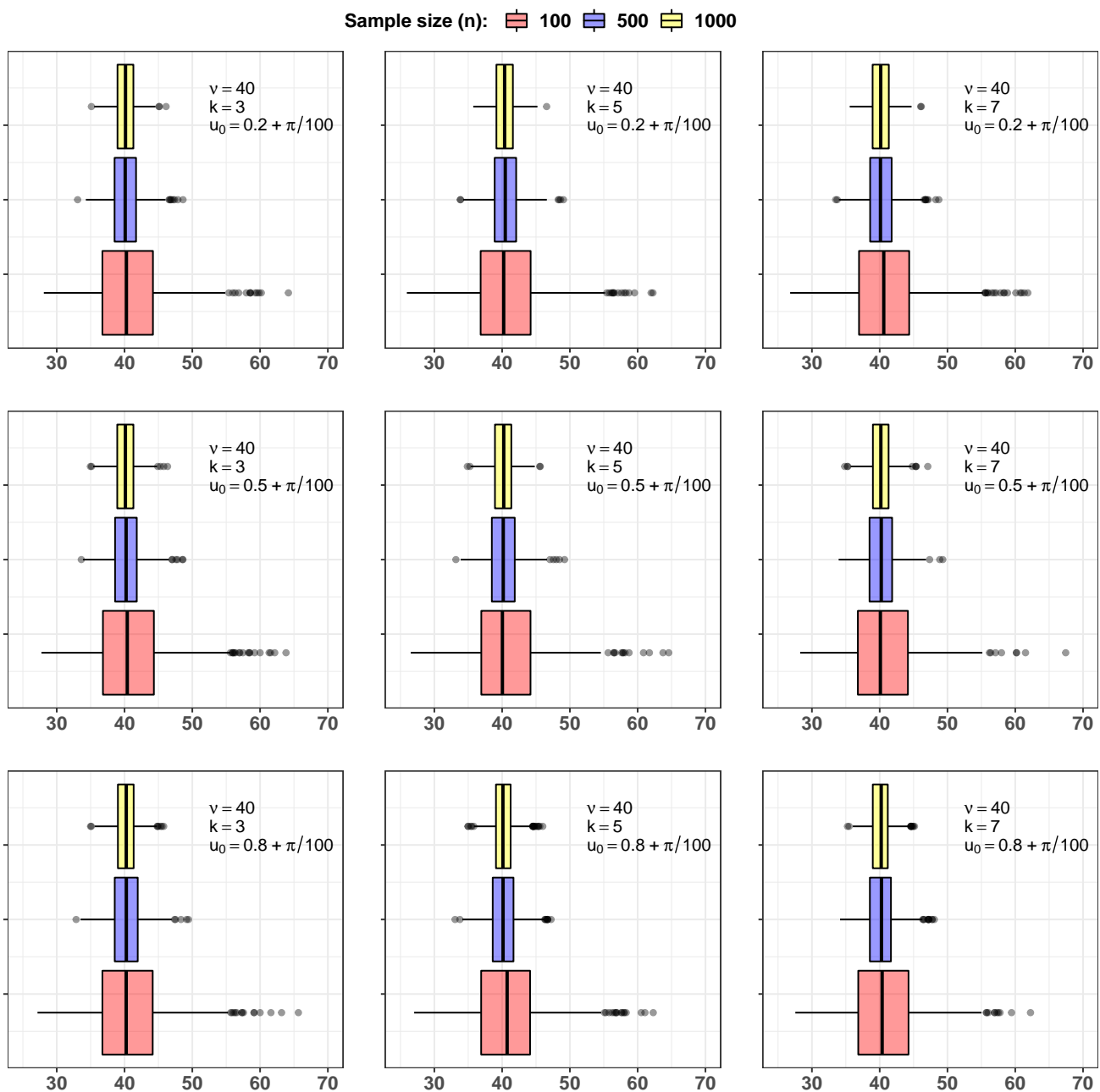

Figure 13: Boxplots of the estimated values, from 1,000 replications, for the parameter $\nu = 40$ considering Map 1 with $k \in \{3, 5, 7\}$ (first, second and third column of each subfigure, respectively), $n \in \{100, 500, 1000\}$ and starting points $u_0 \in \{0.2+\pi/100, 0.5+\pi/100, 0.8+\pi/100\}$ (first, second and third row of each subfigure, respectively).

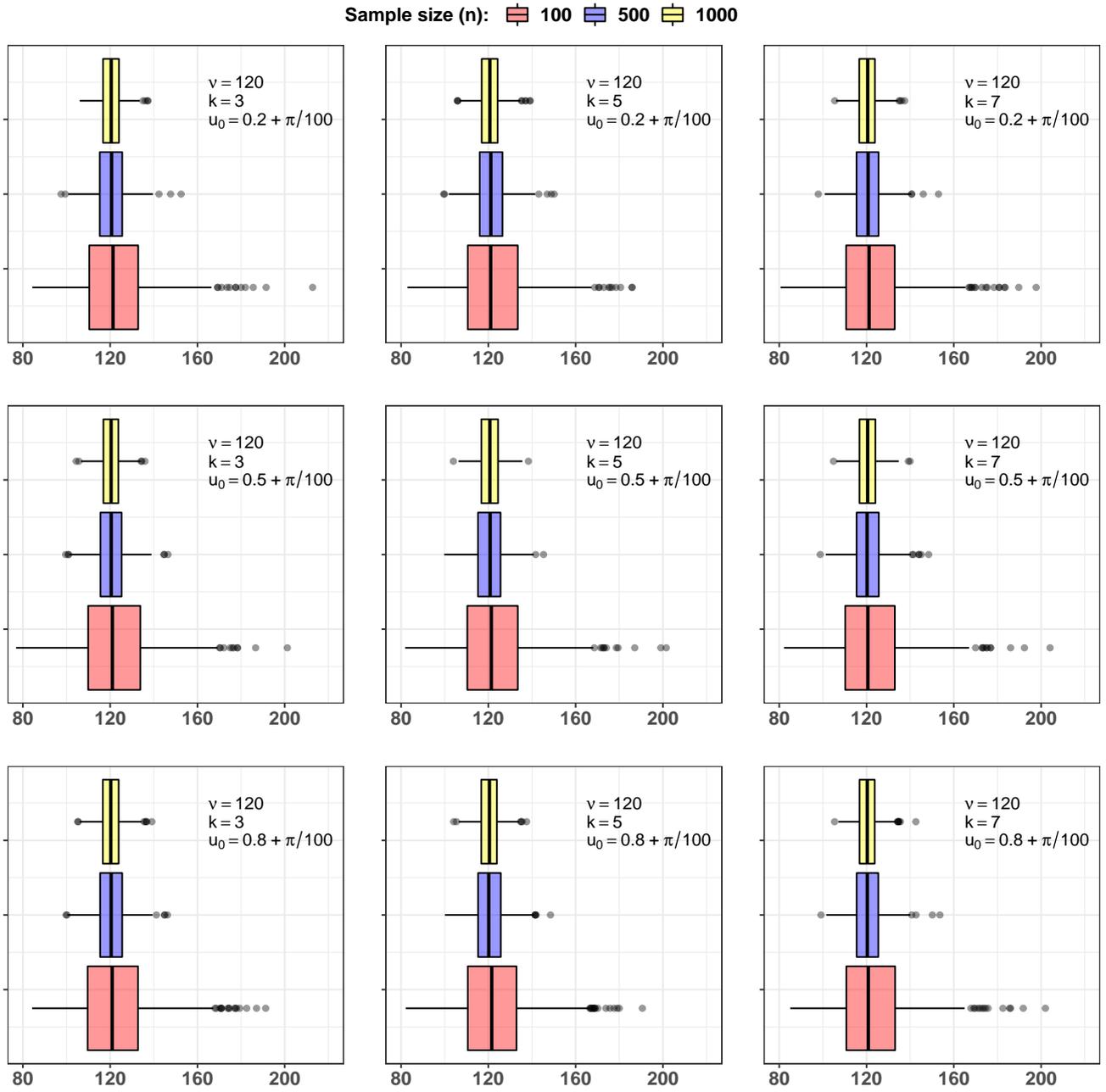

Figure 14: Boxplots of the estimated values, from 1,000 replications, for the parameter $\nu = 120$ considering Map 1 with $k \in \{3, 5, 7\}$ (first, second and third column of each subfigure, respectively), $n \in \{100, 500, 1000\}$ and starting points $u_0 \in \{0.2+\pi/100, 0.5+\pi/100, 0.8+\pi/100\}$ (first, second and third row of each subfigure, respectively).

# Simulation Results for Map 3 - known $\theta$

Table 2: Simulation results for parameter $\nu$ considering the logistic map with $\theta \in \{1, 2, 3, 3.3, 3.6, 3.99\}$ and sample size $n \in \{100, 500, 1000\}$: the mean estimated value of $\nu$ over 1,000 replications ($\bar{\nu}$), for $\nu \in \{10, 40, 120\}$, the standard deviation of the estimates ($sd_\nu$) and the mean absolute percentage error (MAPE).

| $\theta$ | $u_0$ | $n = 100$ | | | $n = 500$ | | | $n = 1,000$ | | |
|---|---|---|---|---|---|---|---|---|---|---|
| | | $\bar{\nu}$ | $sd_\nu$ | MAPE | $\bar{\nu}$ | $sd_\nu$ | MAPE | $\bar{\nu}$ | $sd_\nu$ | MAPE |
| | | | | | $\nu = 10$ | | | | | |
| | $0.2 + \frac{\pi}{100}$ | 10.32 | 2.0676 | 16.10 | **9.62** | 0.8258 | 7.43 | **9.99** | 0.6104 | 4.88 |
| 1 | $0.5 + \frac{\pi}{100}$ | 10.34 | 2.0728 | **16.18** | **9.62** | 0.8207 | **7.45** | **9.99** | **0.6126** | **4.90** |
| | $0.8 + \frac{\pi}{100}$ | 10.32 | **2.0793** | 16.17 | **9.62** | **0.8264** | 7.41 | 10.00 | 0.6103 | 4.87 |
| | $0.2 + \frac{\pi}{100}$ | 10.19 | 1.3792 | 10.78 | 10.03 | 0.5952 | 4.76 | 10.02 | 0.4223 | 3.37 |
| 2 | $0.5 + \frac{\pi}{100}$ | 10.19 | 1.3727 | 10.73 | 10.03 | 0.5943 | 4.76 | 10.02 | 0.4214 | 3.36 |
| | $0.8 + \frac{\pi}{100}$ | 10.19 | 1.3809 | 10.75 | 10.03 | 0.5943 | 4.76 | 10.02 | 0.4222 | 3.37 |
| | $0.2 + \frac{\pi}{100}$ | 10.22 | 1.3935 | 10.96 | 10.04 | 0.5864 | 4.71 | 10.02 | 0.4238 | 3.40 |
| 3 | $0.5 + \frac{\pi}{100}$ | 10.20 | 1.3723 | 10.76 | 10.04 | 0.5816 | 4.69 | 10.02 | 0.4226 | 3.38 |
| | $0.8 + \frac{\pi}{100}$ | 10.22 | 1.4019 | 11.00 | 10.04 | 0.5866 | 4.72 | 10.02 | 0.4236 | 3.40 |
| | $0.2 + \frac{\pi}{100}$ | 10.16 | 1.3306 | 10.55 | **10.02** | 0.5881 | 4.66 | 10.02 | 0.4122 | **3.28** |
| 3.3 | $0.5 + \frac{\pi}{100}$ | 10.16 | 1.3318 | 10.49 | **10.02** | 0.5919 | 4.69 | 10.02 | 0.4151 | 3.30 |
| | $0.8 + \frac{\pi}{100}$ | 10.17 | 1.3259 | 10.52 | 10.03 | 0.5888 | 4.66 | 10.02 | 0.4125 | **3.28** |
| | $0.2 + \frac{\pi}{100}$ | 10.15 | 1.2784 | **10.05** | 10.04 | 0.5664 | **4.51** | 10.02 | 0.4096 | **3.28** |
| 3.6 | $0.5 + \frac{\pi}{100}$ | **10.14** | 1.3099 | 10.22 | **10.02** | 0.5713 | 4.54 | 10.02 | 0.4125 | 3.30 |
| | $0.8 + \frac{\pi}{100}$ | 10.17 | 1.3212 | 10.39 | 10.04 | 0.5737 | 4.55 | 10.02 | 0.4123 | 3.32 |
| | $0.2 + \frac{\pi}{100}$ | 10.29 | 1.2899 | 10.19 | 10.22 | 0.5894 | 4.92 | 10.23 | 0.3973 | 3.67 |
| 3.99 | $0.5 + \frac{\pi}{100}$ | **10.50** | **1.2423** | 10.31 | 10.34 | 0.5606 | 5.21 | **10.26** | 0.3994 | 3.81 |
| | $0.8 + \frac{\pi}{100}$ | 10.33 | 1.2516 | 10.08 | 10.22 | **0.5386** | 4.64 | 10.23 | **0.3840** | 3.54 |
| | | | | | $\nu = 40$ | | | | | |
| | $0.2 + \frac{\pi}{100}$ | 41.43 | **7.1666** | **13.92** | **40.06** | 3.5357 | 6.91 | **38.96** | 2.4457 | 5.37 |
| 1 | $0.5 + \frac{\pi}{100}$ | **41.46** | 7.1134 | 13.84 | 40.07 | 3.5259 | 6.89 | **38.96** | 2.4277 | 5.31 |
| | $0.8 + \frac{\pi}{100}$ | 41.44 | 7.1591 | 13.91 | **40.06** | 3.5578 | **6.99** | **38.96** | 2.4539 | **5.39** |
| | $0.2 + \frac{\pi}{100}$ | 40.91 | 5.8146 | 11.56 | 40.12 | 2.4245 | 4.87 | 40.11 | 1.7406 | 3.46 |
| 2 | $0.5 + \frac{\pi}{100}$ | 40.91 | 5.7889 | 11.49 | 40.12 | 2.4232 | 4.87 | 40.11 | 1.7396 | 3.46 |
| | $0.8 + \frac{\pi}{100}$ | 40.93 | 5.8329 | 11.56 | 40.13 | 2.4254 | 4.88 | 40.12 | 1.7401 | 3.45 |
| | $0.2 + \frac{\pi}{100}$ | 40.91 | 5.7385 | 11.35 | **40.17** | 2.4135 | 4.94 | 40.09 | 1.7159 | 3.44 |
| 3 | $0.5 + \frac{\pi}{100}$ | 40.88 | 5.6447 | 11.10 | 40.16 | 2.4069 | 4.92 | 40.09 | 1.7150 | 3.43 |
| | $0.8 + \frac{\pi}{100}$ | 40.90 | 5.7280 | 11.33 | 40.17 | 2.4191 | 4.95 | 40.09 | 1.7190 | 3.45 |
| | $0.2 + \frac{\pi}{100}$ | 40.71 | 5.6888 | 11.36 | 40.11 | 2.3649 | 4.70 | 40.10 | 1.7427 | 3.53 |
| 3.3 | $0.5 + \frac{\pi}{100}$ | 40.71 | 5.6447 | 11.20 | 40.10 | 2.3530 | 4.70 | 40.09 | 1.7353 | 3.50 |
| | $0.8 + \frac{\pi}{100}$ | 40.72 | 5.6974 | 11.39 | 40.11 | 2.3661 | 4.70 | 40.10 | 1.7429 | 3.53 |
| | $0.2 + \frac{\pi}{100}$ | 40.85 | 5.8039 | 11.34 | 40.15 | 2.3855 | 4.76 | 40.11 | 1.7407 | 3.50 |
| 3.6 | $0.5 + \frac{\pi}{100}$ | 40.87 | 5.7128 | 11.23 | 40.15 | 2.3577 | 4.71 | 40.10 | 1.7201 | 3.45 |
| | $0.8 + \frac{\pi}{100}$ | 40.80 | 5.7541 | 11.21 | 40.11 | 2.4050 | 4.80 | 40.09 | 1.7171 | 3.42 |
| | $0.2 + \frac{\pi}{100}$ | 40.80 | 5.4407 | 10.82 | 40.15 | 2.4033 | 4.76 | 40.11 | **1.6925** | 3.39 |
| 3.99 | $0.5 + \frac{\pi}{100}$ | **40.68** | **5.3732** | **10.54** | 40.15 | 2.3751 | 4.71 | **40.08** | 1.6734 | 3.38 |
| | $0.8 + \frac{\pi}{100}$ | 40.76 | 5.4431 | 10.62 | 40.13 | **2.3325** | **4.66** | 40.09 | 1.6801 | **3.36** |

Table 2 – *Continued from previous page*

| | $\nu = 120$ | | | | | | | | | |
|---|---|---|---|---|---|---|---|---|---|---|
| 1 | $0.2 + \frac{\pi}{100}$ | 124.36 | 18.0834 | 12.07 | **121.08** | 9.1308 | 6.08 | 120.38 | 6.9623 | 4.58 |
| | $0.5 + \frac{\pi}{100}$ | **124.41** | 18.1450 | 12.07 | 121.05 | **9.1493** | 6.07 | 120.37 | 6.9619 | 4.57 |
| | $0.8 + \frac{\pi}{100}$ | 124.28 | **18.2572** | **12.15** | 121.06 | 9.1413 | **6.12** | 120.36 | **6.9714** | **4.59** |
| 2 | $0.2 + \frac{\pi}{100}$ | 122.74 | 17.3998 | 11.51 | 120.38 | 7.2552 | 4.86 | 120.37 | 5.2043 | 3.47 |
| | $0.5 + \frac{\pi}{100}$ | 122.75 | 17.3952 | 11.50 | 120.39 | 7.2535 | 4.86 | 120.37 | 5.2049 | 3.47 |
| | $0.8 + \frac{\pi}{100}$ | 122.78 | 17.4615 | 11.54 | 120.39 | 7.2625 | 4.86 | 120.37 | 5.2101 | 3.47 |
| 3 | $0.2 + \frac{\pi}{100}$ | 122.92 | 17.3504 | 11.41 | 120.55 | 7.3037 | 4.94 | **120.40** | 5.2727 | 3.51 |
| | $0.5 + \frac{\pi}{100}$ | 122.99 | 17.2597 | 11.42 | 120.53 | 7.3254 | 4.96 | 120.39 | 5.2812 | 3.50 |
| | $0.8 + \frac{\pi}{100}$ | 122.96 | 17.3952 | 11.44 | 120.55 | 7.3092 | 4.95 | **120.40** | 5.2729 | 3.51 |
| 3.3 | $0.2 + \frac{\pi}{100}$ | 122.62 | 17.3120 | 11.44 | 120.39 | 7.2279 | 4.83 | 120.29 | 5.1907 | 3.45 |
| | $0.5 + \frac{\pi}{100}$ | 122.67 | 17.1222 | 11.33 | 120.39 | 7.1640 | 4.81 | 120.29 | 5.1684 | 3.44 |
| | $0.8 + \frac{\pi}{100}$ | 122.62 | 17.3228 | 11.44 | 120.38 | 7.2265 | 4.83 | 120.29 | 5.1916 | 3.45 |
| 3.6 | $0.2 + \frac{\pi}{100}$ | 122.74 | 17.5434 | 11.70 | 120.41 | 7.2348 | 4.86 | 120.32 | 5.2341 | 3.51 |
| | $0.5 + \frac{\pi}{100}$ | 122.54 | 16.9107 | 11.22 | 120.40 | 7.1674 | 4.77 | 120.33 | 5.1701 | 3.45 |
| | $0.8 + \frac{\pi}{100}$ | 122.84 | 17.4679 | 11.57 | 120.46 | 7.2986 | 4.86 | 120.34 | 5.2397 | 3.50 |
| 3.99 | $0.2 + \frac{\pi}{100}$ | 122.54 | 16.8729 | 10.98 | 120.38 | 7.2459 | 4.77 | 120.27 | **4.9935** | 3.34 |
| | $0.5 + \frac{\pi}{100}$ | 122.51 | **16.7931** | 11.16 | 120.45 | 7.2631 | 4.82 | 120.34 | 5.0686 | **3.40** |
| | $0.8 + \frac{\pi}{100}$ | **122.03** | 16.8531 | **11.04** | **120.31** | **7.0893** | **4.68** | **120.24** | 5.1664 | 3.43 |

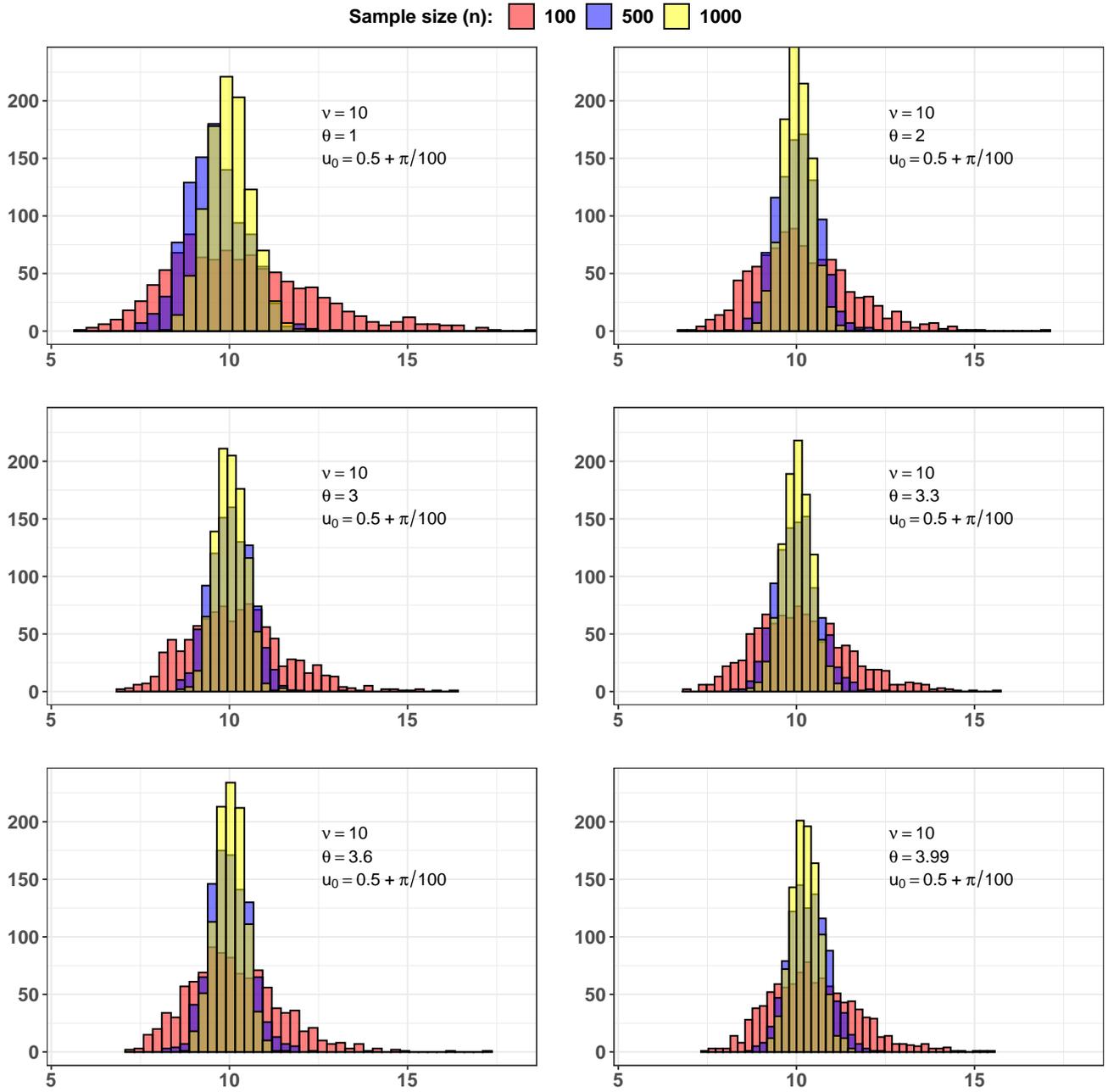

Figure 15: Histogram of the estimated values, from 1,000 replications, for the parameter $\nu = 10$ considering logistic map with $\theta \in \{1, 2, 3, 3.3, 3.6, 3.99\}$, $n \in \{100, 500, 1000\}$ and $u_0 = 0.5 + \pi/100$.

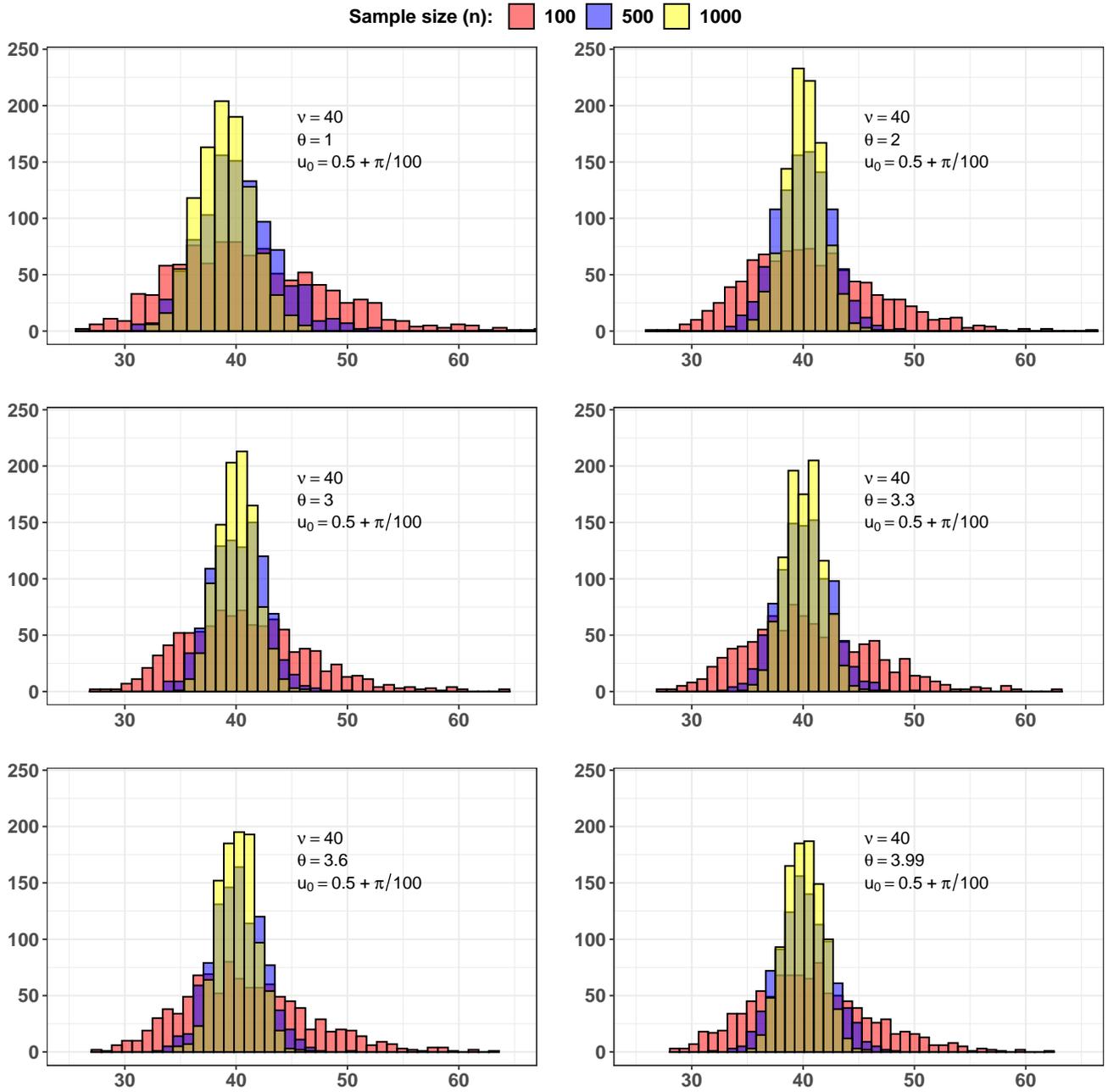

Figure 16: Histogram of the estimated values, from 1,000 replications, for the parameter $\nu = 40$ considering logistic map with $\theta \in \{1, 2, 3, 3.3, 3.6, 3.99\}$, $n \in \{100, 500, 1000\}$ and $u_0 = 0.5 + \pi/100$.

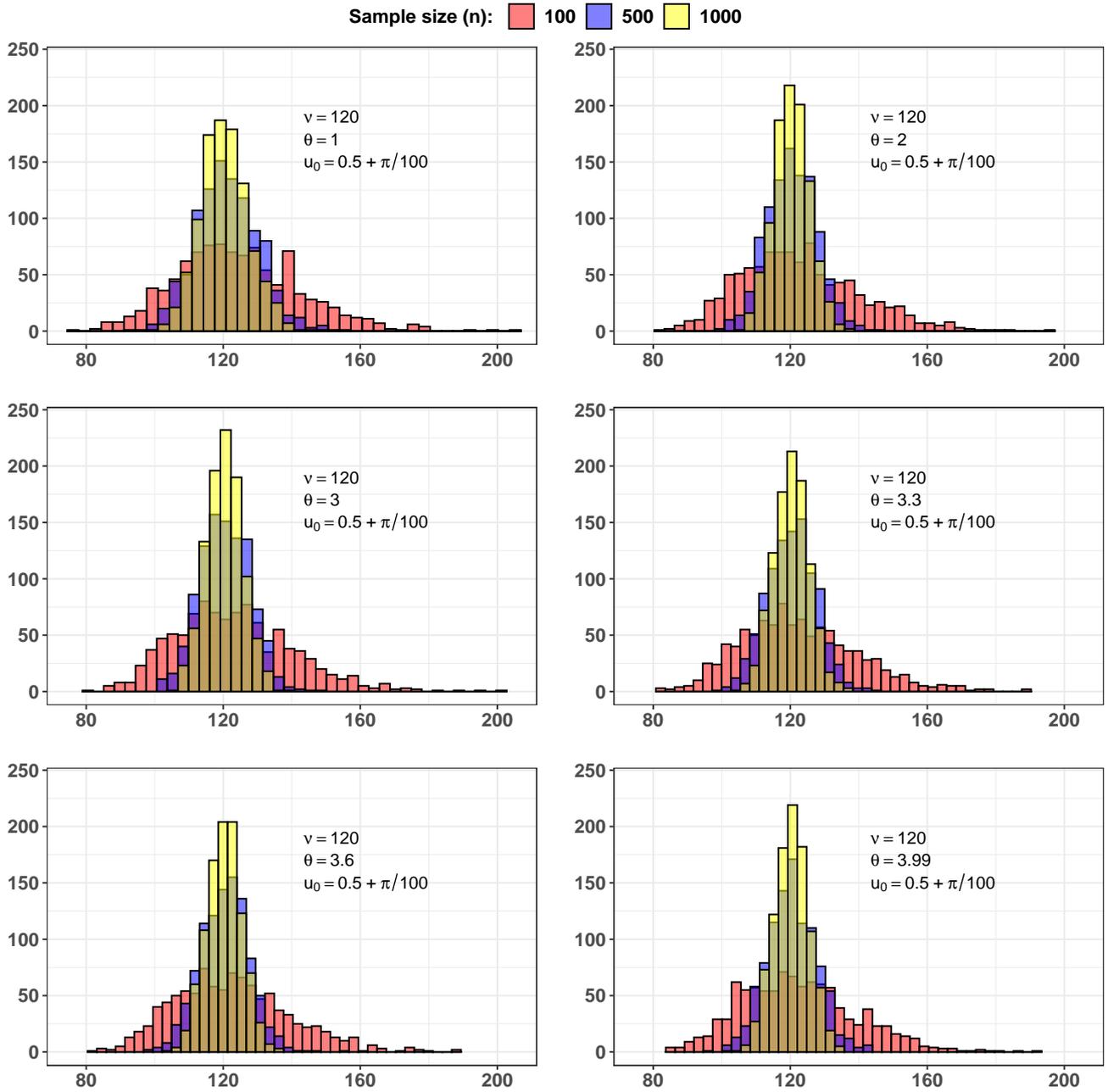

Figure 17: Histogram of the estimated values, from 1,000 replications, for the parameter $\nu = 120$ considering logistic map with $\theta \in \{1, 2, 3, 3.3, 3.6, 3.99\}$, $n \in \{100, 500, 1000\}$ and $u_0 = 0.5 + \pi/100$.

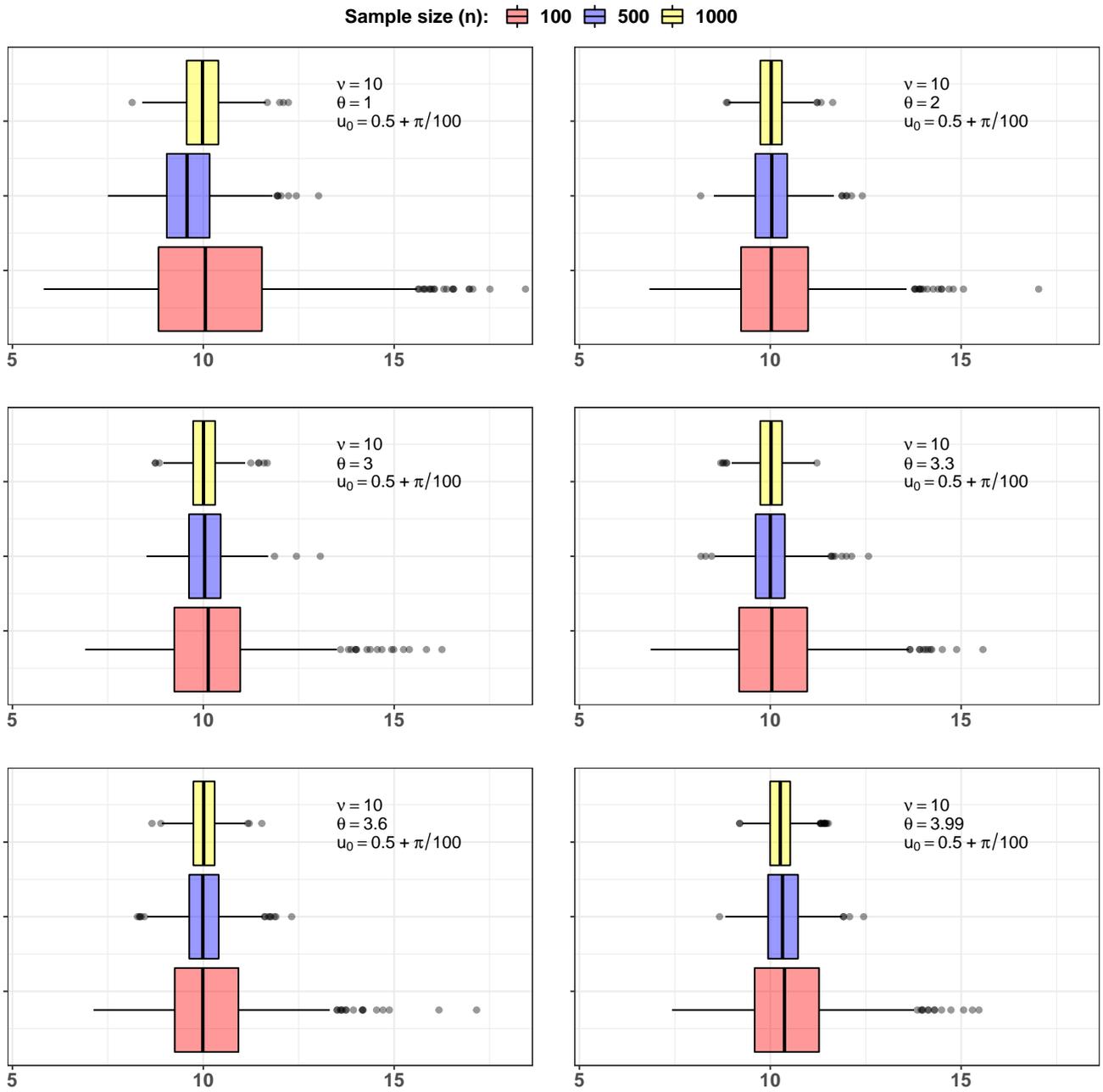

Figure 18: Boxplots of the estimated values, from 1,000 replications, for the parameter $\nu = 10$ considering the logistic map with $\theta \in \{1, 2, 3, 3.3, 3.6, 3.99\}$, $n \in \{100, 500, 1000\}$ and $u_0 = 0.5 + \pi/100$.

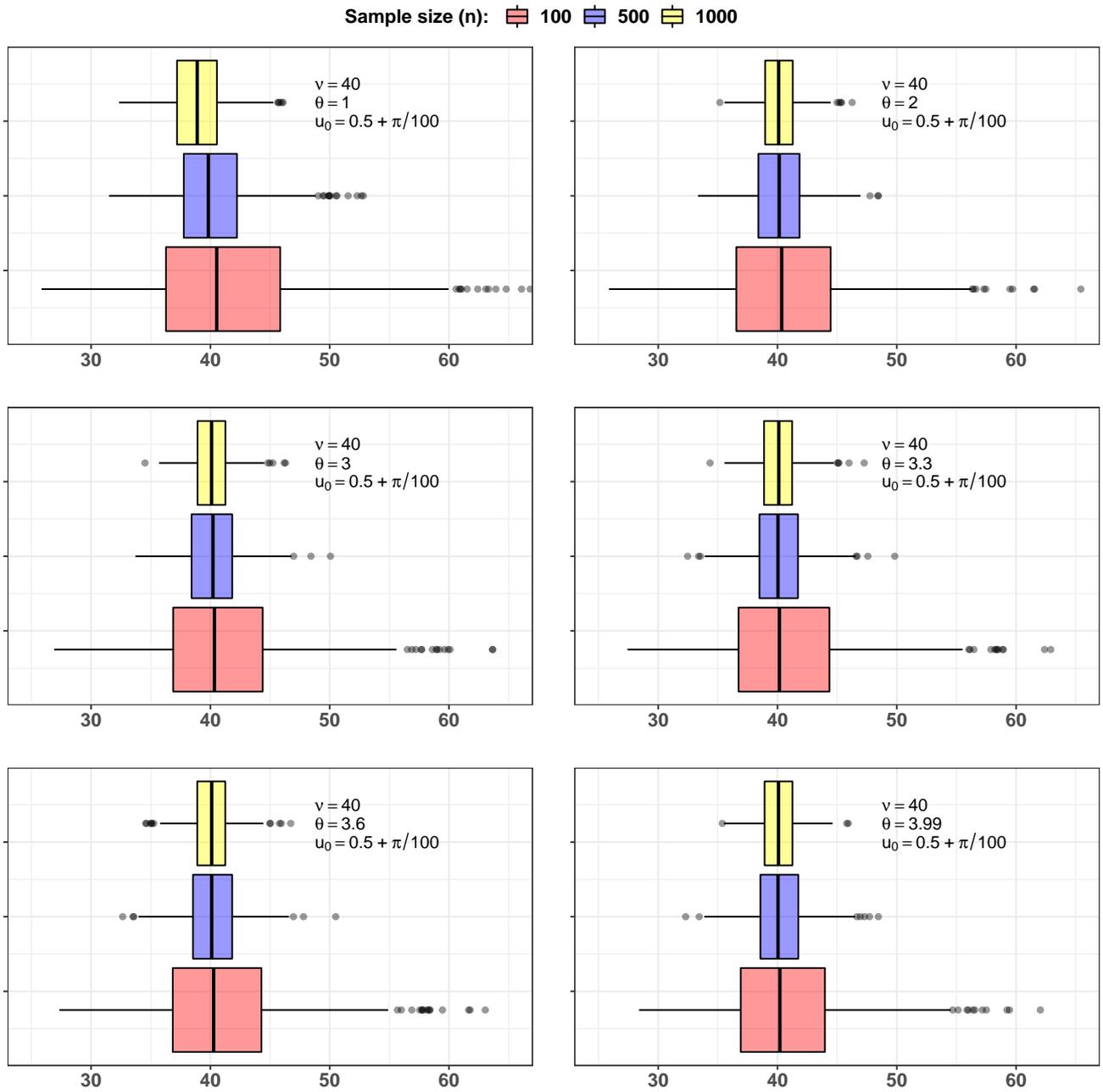

Figure 19: Boxplots of the estimated values, from 1,000 replications, for the parameter $\nu = 40$ considering the logistic map with $\theta \in \{1, 2, 3, 3.3, 3.6, 3.99\}$, $n \in \{100, 500, 1000\}$ and $u_0 = 0.5 + \pi/100$.

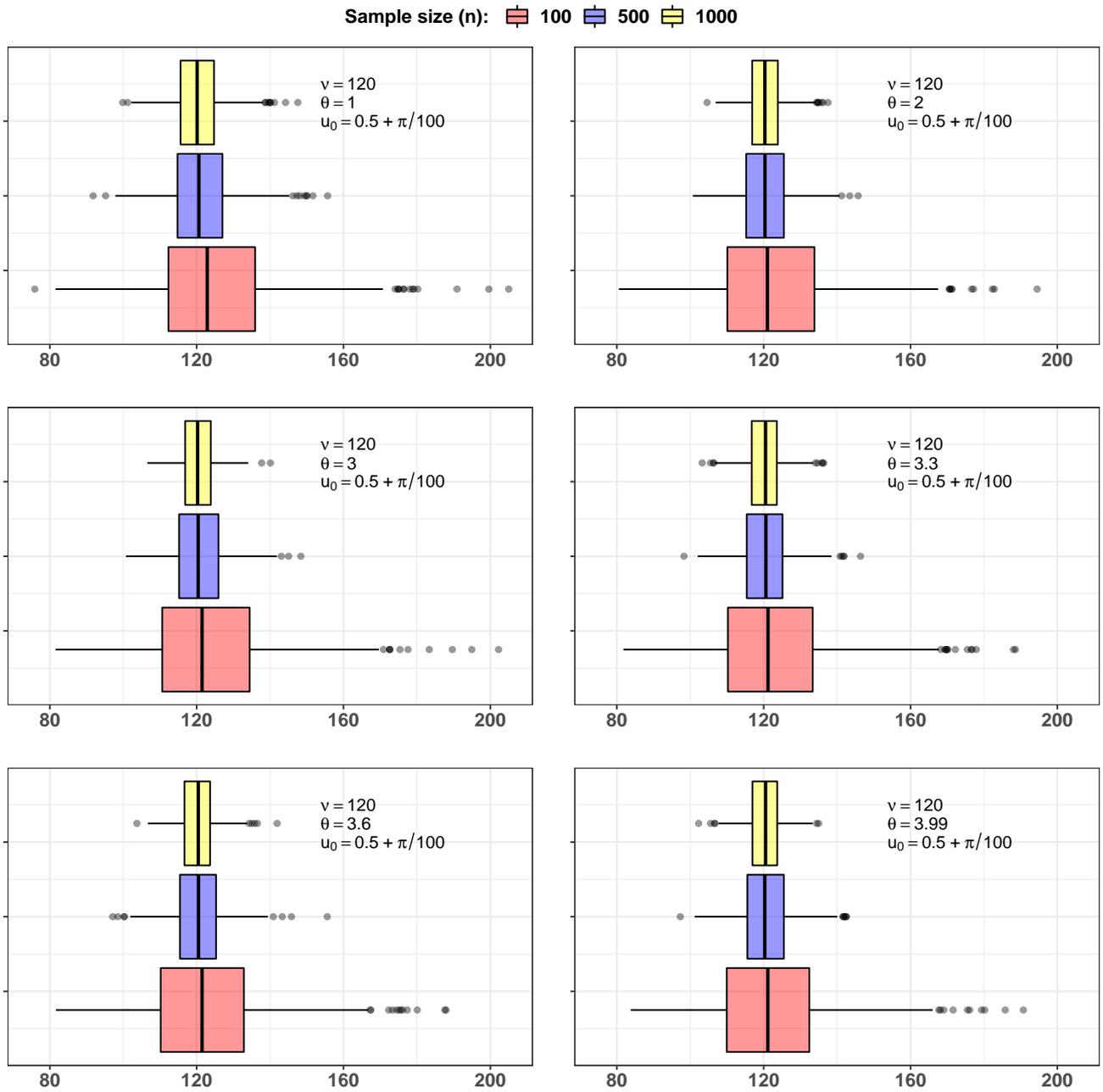

Figure 20: Boxplots of the estimated values, from 1,000 replications, for the parameter $\nu = 120$ considering the logistic map with $\theta \in \{1, 2, 3, 3.3, 3.6, 3.99\}$, $n \in \{100, 500, 1000\}$ and $u_0 = 0.5 + \pi/100$.

# Simulation Results for Map 3 - unknown $\theta$

Table 3: Simulation Results for parameter $\nu$ considering the logistic map with $\theta \in \{1, 2, 3, 3.3\}$ and sample size $n \in \{100, 500, 1000\}$.

| $\theta$ | $u_0$ | $n = 100$ | | | $n = 500$ | | | $n = 1,000$ | | |
|---|---|---|---|---|---|---|---|---|---|---|
| | | $\bar{\nu}$ | $sd_\nu$ | MAPE | $\bar{\nu}$ | $sd_\nu$ | MAPE | $\bar{\nu}$ | $sd_\nu$ | MAPE |
| | | | | | $\nu = 10$ | | | | | |
| 1 | $0.2 + \frac{\pi}{100}$ | 10.32 | 2.0676 | 16.11 | **9.62** | 0.8259 | 7.43 | 9.99 | 0.6104 | 4.88 |
| | $0.5 + \frac{\pi}{100}$ | **10.34** | 2.0727 | **16.18** | **9.62** | **0.8270** | **7.45** | 9.99 | **0.6126** | **4.90** |
| | $0.8 + \frac{\pi}{100}$ | 10.32 | **2.0793** | 16.17 | **9.62** | 0.8264 | 7.41 | **10.00** | 0.6102 | 4.87 |
| 2 | $0.2 + \frac{\pi}{100}$ | 10.29 | 1.3976 | 10.94 | **10.05** | 0.5973 | 4.78 | **10.03** | 0.4230 | 3.38 |
| | $0.5 + \frac{\pi}{100}$ | 10.29 | **1.3914** | 10.88 | **10.05** | 0.5967 | 4.79 | **10.03** | 0.4222 | **3.37** |
| | $0.8 + \frac{\pi}{100}$ | 10.29 | 1.3998 | 10.94 | **10.05** | 0.5965 | 4.78 | **10.03** | 0.4230 | 3.38 |
| 3 | $0.2 + \frac{\pi}{100}$ | 10.30 | 1.4264 | 11.24 | **10.05** | 0.5895 | 4.74 | 10.02 | 0.4238 | 3.40 |
| | $0.5 + \frac{\pi}{100}$ | **10.28** | 1.4032 | 11.02 | **10.05** | **0.5833** | **4.70** | **10.03** | **0.4216** | **3.37** |
| | $0.8 + \frac{\pi}{100}$ | 10.30 | 1.4307 | 11.24 | **10.05** | 0.5898 | 4.75 | **10.03** | 0.4233 | 3.40 |
| 3.3 | $0.2 + \frac{\pi}{100}$ | 10.29 | 1.4017 | 11.15 | **10.05** | 0.6019 | 4.79 | **10.03** | 0.4217 | 3.38 |
| | $0.5 + \frac{\pi}{100}$ | 10.30 | 1.4138 | 11.17 | **10.05** | 0.6057 | 4.82 | **10.03** | 0.4247 | 3.40 |
| | $0.8 + \frac{\pi}{100}$ | 10.29 | 1.3936 | 11.15 | **10.05** | 0.6026 | 4.79 | **10.03** | 0.4218 | 3.38 |
| | | | | | $\nu = 40$ | | | | | |
| 1 | $0.2 + \frac{\pi}{100}$ | 41.43 | **7.1667** | **13.92** | 40.06 | 3.5357 | 6.91 | **38.96** | 2.4458 | 5.37 |
| | $0.5 + \frac{\pi}{100}$ | **41.46** | 7.1133 | 13.84 | 40.07 | 3.5259 | 6.89 | **38.96** | 2.4277 | 5.31 |
| | $0.8 + \frac{\pi}{100}$ | 41.44 | 7.1590 | 13.91 | **40.06** | **3.5578** | **6.99** | **38.96** | **2.4537** | **5.39** |
| 2 | $0.2 + \frac{\pi}{100}$ | 41.31 | 5.8954 | 11.78 | 40.21 | 2.4368 | 4.91 | 40.15 | 1.7444 | 3.47 |
| | $0.5 + \frac{\pi}{100}$ | 41.30 | 5.8607 | 11.67 | 40.21 | 2.4355 | 4.91 | 40.15 | 1.7436 | 3.47 |
| | $0.8 + \frac{\pi}{100}$ | 41.34 | 5.9101 | 11.77 | 40.21 | 2.4372 | 4.92 | 40.15 | 1.7440 | 3.47 |
| 3 | $0.2 + \frac{\pi}{100}$ | 41.34 | 5.8360 | 11.61 | **40.25** | 2.4270 | 4.98 | **40.13** | 1.7194 | 3.45 |
| | $0.5 + \frac{\pi}{100}$ | 41.30 | **5.7451** | **11.38** | **40.25** | **2.4216** | 4.96 | **40.13** | **1.7191** | **3.44** |
| | $0.8 + \frac{\pi}{100}$ | 41.32 | 5.8244 | 11.60 | **40.25** | 2.4326 | 5.00 | **40.13** | 1.7222 | 3.46 |
| 3.3 | $0.2 + \frac{\pi}{100}$ | **41.13** | 5.8369 | 11.73 | 40.20 | 2.3891 | 4.76 | 40.14 | 1.7499 | 3.54 |
| | $0.5 + \frac{\pi}{100}$ | 41.15 | 5.8131 | 11.58 | 40.20 | 2.3744 | **4.75** | 40.14 | 1.7416 | 3.51 |
| | $0.8 + \frac{\pi}{100}$ | 41.14 | 5.8511 | 11.79 | 40.20 | 2.3882 | 4.76 | 40.14 | 1.7491 | 3.54 |
| | | | | | $\nu = 120$ | | | | | |
| 1 | $0.2 + \frac{\pi}{100}$ | 124.36 | 18.0832 | 12.07 | **121.08** | 9.1306 | 6.08 | 120.38 | 6.9625 | 4.58 |
| | $0.5 + \frac{\pi}{100}$ | **124.41** | 18.1453 | 12.07 | 121.05 | **9.1494** | 6.07 | 120.37 | 6.9620 | 4.57 |
| | $0.8 + \frac{\pi}{100}$ | 124.28 | **18.2570** | **12.15** | 121.06 | 9.1414 | **6.12** | **120.36** | **6.9713** | **4.59** |
| 2 | $0.2 + \frac{\pi}{100}$ | 123.94 | 17.6421 | 11.75 | **120.64** | 7.2900 | 4.90 | 120.48 | 5.2150 | 3.48 |
| | $0.5 + \frac{\pi}{100}$ | 123.94 | 17.5939 | 11.69 | **120.64** | 7.2872 | 4.90 | 120.49 | 5.2172 | 3.48 |
| | $0.8 + \frac{\pi}{100}$ | 123.99 | 17.6759 | 11.77 | **120.64** | 7.2971 | 4.90 | 120.49 | 5.2215 | 3.48 |
| 3 | $0.2 + \frac{\pi}{100}$ | 124.22 | 17.6417 | 11.67 | 120.81 | 7.3282 | 4.97 | **120.52** | 5.2893 | 3.53 |
| | $0.5 + \frac{\pi}{100}$ | 124.32 | 17.5835 | 11.74 | 120.78 | 7.3543 | 4.99 | 120.50 | 5.3007 | 3.53 |
| | $0.8 + \frac{\pi}{100}$ | 124.25 | 17.6947 | 11.72 | 120.81 | 7.3341 | 4.98 | **120.52** | 5.2889 | 3.53 |
| 3.3 | $0.2 + \frac{\pi}{100}$ | 123.87 | 17.6318 | 11.71 | 120.67 | 7.2654 | 4.87 | 120.42 | 5.2162 | 3.47 |
| | $0.5 + \frac{\pi}{100}$ | 123.95 | **17.4404** | **11.60** | 120.67 | **7.1998** | **4.85** | 120.42 | **5.1910** | **3.46** |
| | $0.8 + \frac{\pi}{100}$ | **123.86** | 17.6250 | 11.72 | 120.66 | 7.2592 | 4.87 | 120.42 | 5.2151 | 3.47 |

Table 4: Simulation Results for parameter $\theta$ considering the logistic map with $\theta \in \{1, 2, 3, 3.3\}$ and sample size $n \in \{100, 500, 1000\}$.

| $\theta$ | $u_0$ | $n = 100$ | | | $n = 500$ | | | $n = 1{,}000$ | | |
|---|---|---|---|---|---|---|---|---|---|---|
| | | $\bar{\theta}$ | $sd_\theta$ | MAPE | $\bar{\theta}$ | $sd_\theta$ | MAPE | $\bar{\theta}$ | $sd_\theta$ | MAPE |
| | | | | | $\nu = 10$ | | | | | |
| 1 | $0.2 + \frac{\pi}{100}$ | **1.00** | 0.0056 | **0.44** | **1.00** | **0.0006** | 0.33 | **1.00** | **0.0003** | 0.34 |
| | $0.5 + \frac{\pi}{100}$ | **1.00** | 0.0057 | 0.46 | **1.00** | **0.0006** | 0.33 | **1.00** | **0.0003** | 0.34 |
| | $0.8 + \frac{\pi}{100}$ | **1.00** | **0.0055** | **0.44** | **1.00** | **0.0006** | 0.33 | **1.00** | **0.0003** | 0.34 |
| 2 | $0.2 + \frac{\pi}{100}$ | **2.00** | **0.0595** | **2.37** | **2.00** | **0.0277** | 1.10 | **2.00** | 0.0187 | **0.75** |
| | $0.5 + \frac{\pi}{100}$ | **2.00** | 0.0590 | 2.33 | **2.00** | **0.0277** | **1.11** | **2.00** | **0.0189** | 0.76 |
| | $0.8 + \frac{\pi}{100}$ | **2.00** | 0.0589 | 2.34 | **2.00** | 0.0276 | 1.10 | **2.00** | 0.0187 | **0.75** |
| 3 | $0.2 + \frac{\pi}{100}$ | **2.99** | 0.0472 | 0.82 | **2.99** | 0.0203 | 0.25 | **3.00** | 0.0132 | 0.15 |
| | $0.5 + \frac{\pi}{100}$ | **2.99** | 0.0450 | 0.76 | **3.00** | 0.0171 | **0.21** | **3.00** | 0.0128 | **0.14** |
| | $0.8 + \frac{\pi}{100}$ | **2.99** | 0.0493 | 0.83 | **3.00** | 0.0193 | 0.24 | **3.00** | 0.0127 | **0.14** |
| 3.3 | $0.2 + \frac{\pi}{100}$ | **3.31** | 0.0543 | 1.27 | **3.30** | 0.0238 | 0.57 | **3.30** | 0.0160 | 0.39 |
| | $0.5 + \frac{\pi}{100}$ | **3.31** | 0.0553 | 1.28 | **3.30** | 0.0237 | 0.57 | **3.30** | 0.0160 | 0.39 |
| | $0.8 + \frac{\pi}{100}$ | **3.31** | 0.0529 | 1.25 | **3.30** | 0.0236 | 0.57 | **3.30** | 0.0160 | 0.39 |
| | | | | | $\nu = 40$ | | | | | |
| 1 | $0.2 + \frac{\pi}{100}$ | **1.00** | **0.0038** | **0.30** | **1.00** | **0.0005** | **0.05** | **1.00** | **0.0002** | 0.06 |
| | $0.5 + \frac{\pi}{100}$ | **1.00** | 0.0039 | 0.31 | **1.00** | **0.0005** | **0.05** | **1.00** | **0.0002** | 0.06 |
| | $0.8 + \frac{\pi}{100}$ | **1.00** | **0.0038** | **0.30** | **1.00** | **0.0005** | **0.05** | **1.00** | **0.0002** | 0.06 |
| 2 | $0.2 + \frac{\pi}{100}$ | **2.00** | **0.0306** | **1.23** | **2.00** | **0.0142** | **0.57** | **2.00** | 0.0097 | **0.39** |
| | $0.5 + \frac{\pi}{100}$ | **2.00** | **0.0306** | 1.21 | **2.00** | **0.0142** | **0.57** | **2.00** | **0.0098** | **0.39** |
| | $0.8 + \frac{\pi}{100}$ | **2.00** | **0.0306** | 1.22 | **2.00** | 0.0141 | 0.56 | **2.00** | 0.0097 | 0.38 |
| 3 | $0.2 + \frac{\pi}{100}$ | **3.00** | 0.0095 | 0.24 | **3.00** | 0.0025 | **0.05** | **3.00** | 0.0017 | 0.03 |
| | $0.5 + \frac{\pi}{100}$ | **3.00** | 0.0108 | 0.25 | **3.00** | 0.0019 | **0.05** | **3.00** | 0.0010 | **0.02** |
| | $0.8 + \frac{\pi}{100}$ | **3.00** | 0.0092 | 0.23 | **3.00** | 0.0027 | **0.05** | **3.00** | 0.0017 | 0.03 |
| 3.3 | $0.2 + \frac{\pi}{100}$ | **3.30** | 0.0266 | 0.64 | **3.30** | 0.0124 | 0.30 | **3.30** | 0.0085 | 0.21 |
| | $0.5 + \frac{\pi}{100}$ | **3.30** | 0.0267 | 0.64 | **3.30** | 0.0124 | 0.30 | **3.30** | 0.0085 | 0.21 |
| | $0.8 + \frac{\pi}{100}$ | **3.30** | 0.0263 | 0.64 | **3.30** | 0.0123 | 0.30 | **3.30** | 0.0085 | 0.21 |
| | | | | | $\nu = 120$ | | | | | |
| 1 | $0.2 + \frac{\pi}{100}$ | **1.00** | **0.0023** | 0.19 | **1.00** | **0.0004** | 0.03 | **1.00** | **0.0002** | **0.01** |
| | $0.5 + \frac{\pi}{100}$ | **1.00** | 0.0024 | 0.19 | **1.00** | **0.0004** | 0.03 | **1.00** | **0.0002** | **0.01** |
| | $0.8 + \frac{\pi}{100}$ | **1.00** | **0.0023** | 0.18 | **1.00** | **0.0004** | 0.03 | **1.00** | **0.0002** | **0.01** |
| 2 | $0.2 + \frac{\pi}{100}$ | **2.00** | **0.0178** | **0.71** | **2.00** | 0.0082 | **0.33** | **2.00** | **0.0056** | **0.23** |
| | $0.5 + \frac{\pi}{100}$ | **2.00** | **0.0178** | 0.70 | **2.00** | **0.0083** | **0.33** | **2.00** | 0.0057 | **0.23** |
| | $0.8 + \frac{\pi}{100}$ | **2.00** | **0.0178** | **0.71** | **2.00** | 0.0082 | **0.33** | **2.00** | **0.0056** | 0.22 |
| 3 | $0.2 + \frac{\pi}{100}$ | **3.00** | 0.0051 | 0.14 | **3.00** | 0.0011 | **0.03** | **3.00** | 0.0005 | **0.01** |
| | $0.5 + \frac{\pi}{100}$ | **3.00** | 0.0051 | **0.13** | **3.00** | 0.0011 | **0.03** | **3.00** | 0.0005 | **0.01** |
| | $0.8 + \frac{\pi}{100}$ | **3.00** | 0.0050 | **0.13** | **3.00** | 0.0011 | **0.03** | **3.00** | 0.0005 | **0.01** |
| 3.3 | $0.2 + \frac{\pi}{100}$ | **3.30** | 0.0156 | 0.37 | **3.30** | 0.0072 | 0.18 | **3.30** | 0.0050 | 0.12 |
| | $0.5 + \frac{\pi}{100}$ | **3.30** | 0.0156 | 0.38 | **3.30** | 0.0072 | 0.17 | **3.30** | 0.0050 | 0.12 |
| | $0.8 + \frac{\pi}{100}$ | **3.30** | 0.0155 | 0.37 | **3.30** | 0.0072 | 0.17 | **3.30** | 0.0050 | 0.12 |

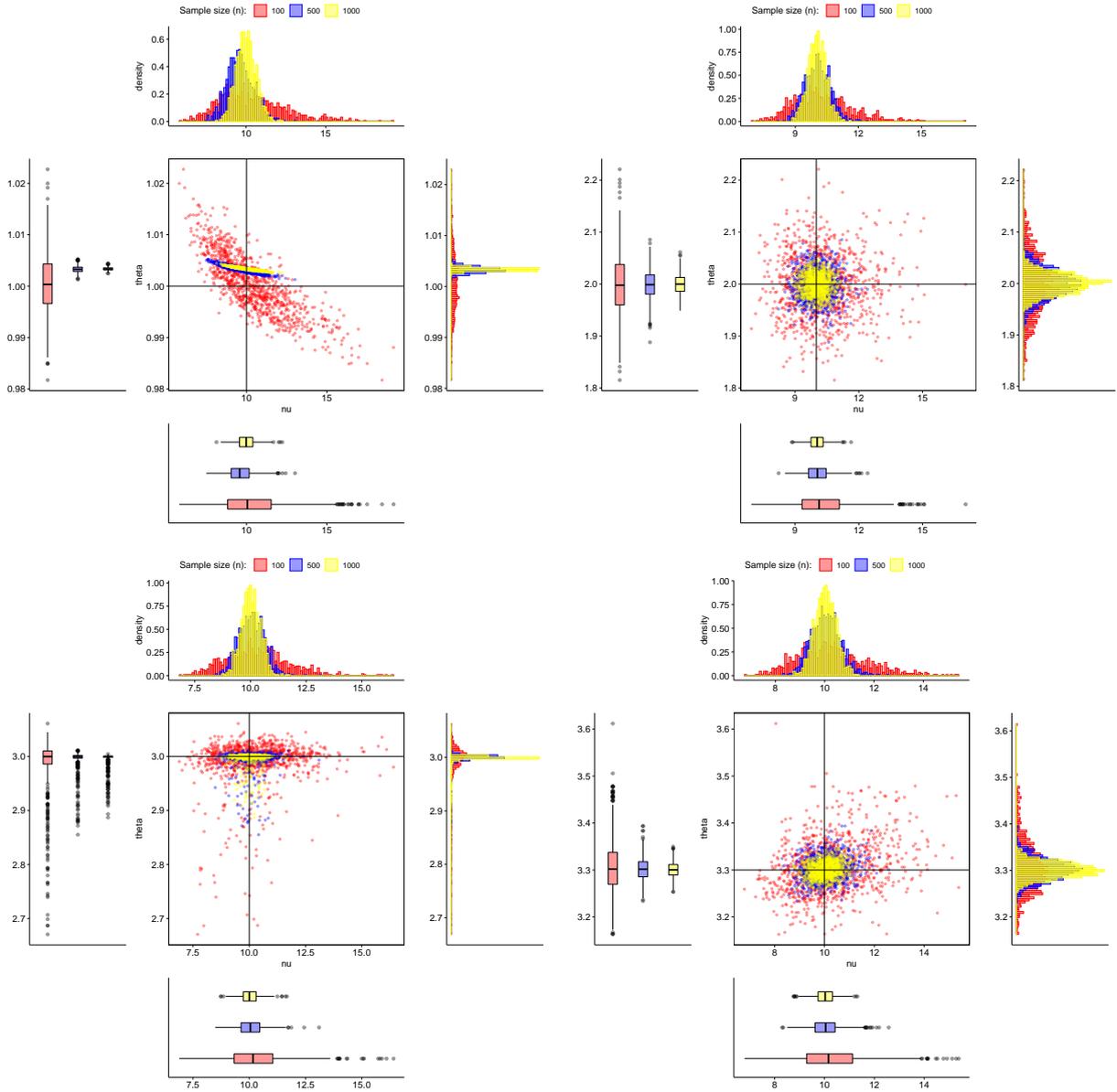

Figure 21: Histogram and Boxplots for $\nu = 10$ and $\theta \in \{1, 2, 3, 3.3\}$.

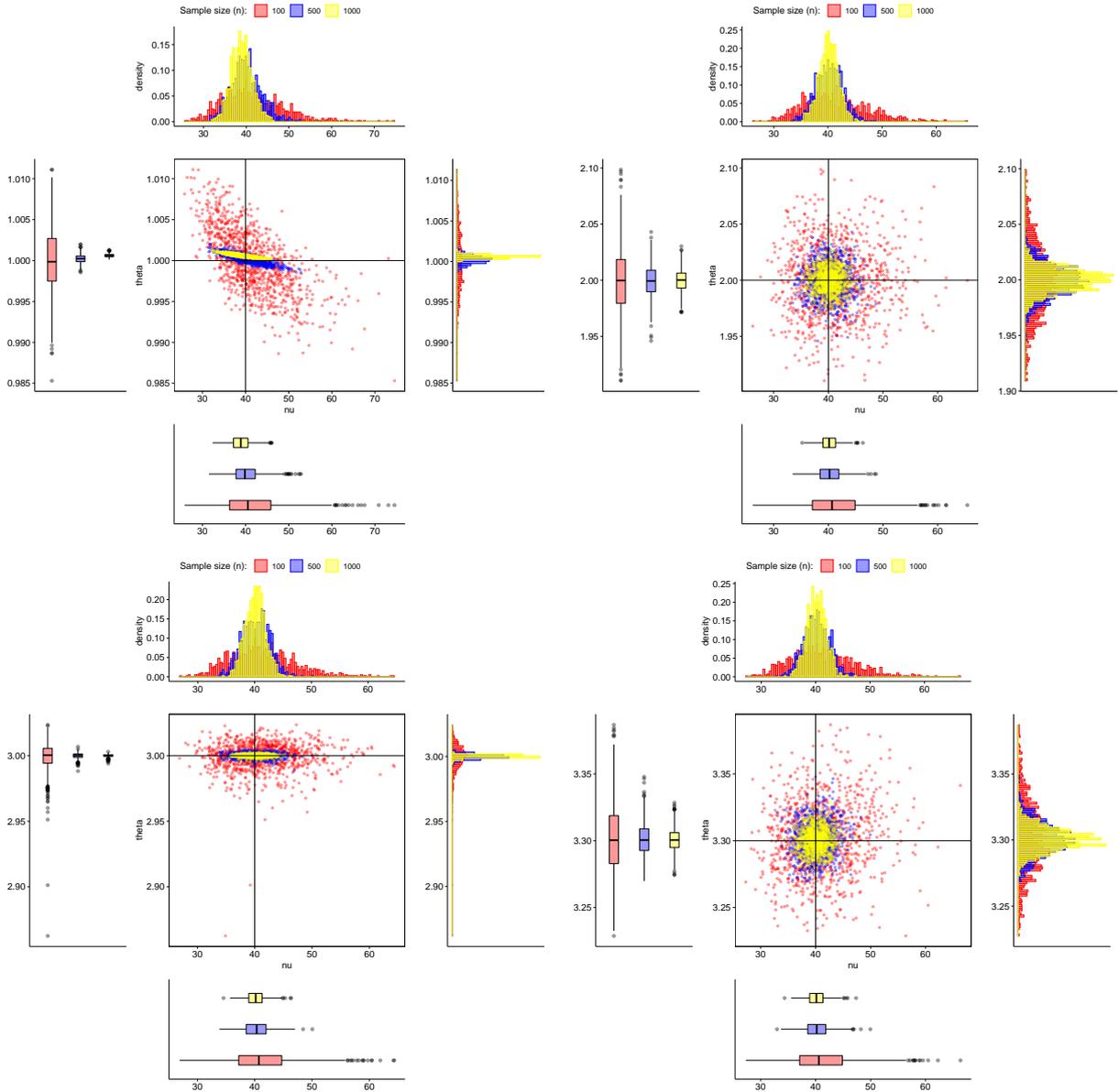

Figure 22: Histogram and Boxplots for $\nu = 40$ and $\theta \in \{1, 2, 3, 3.3\}$.

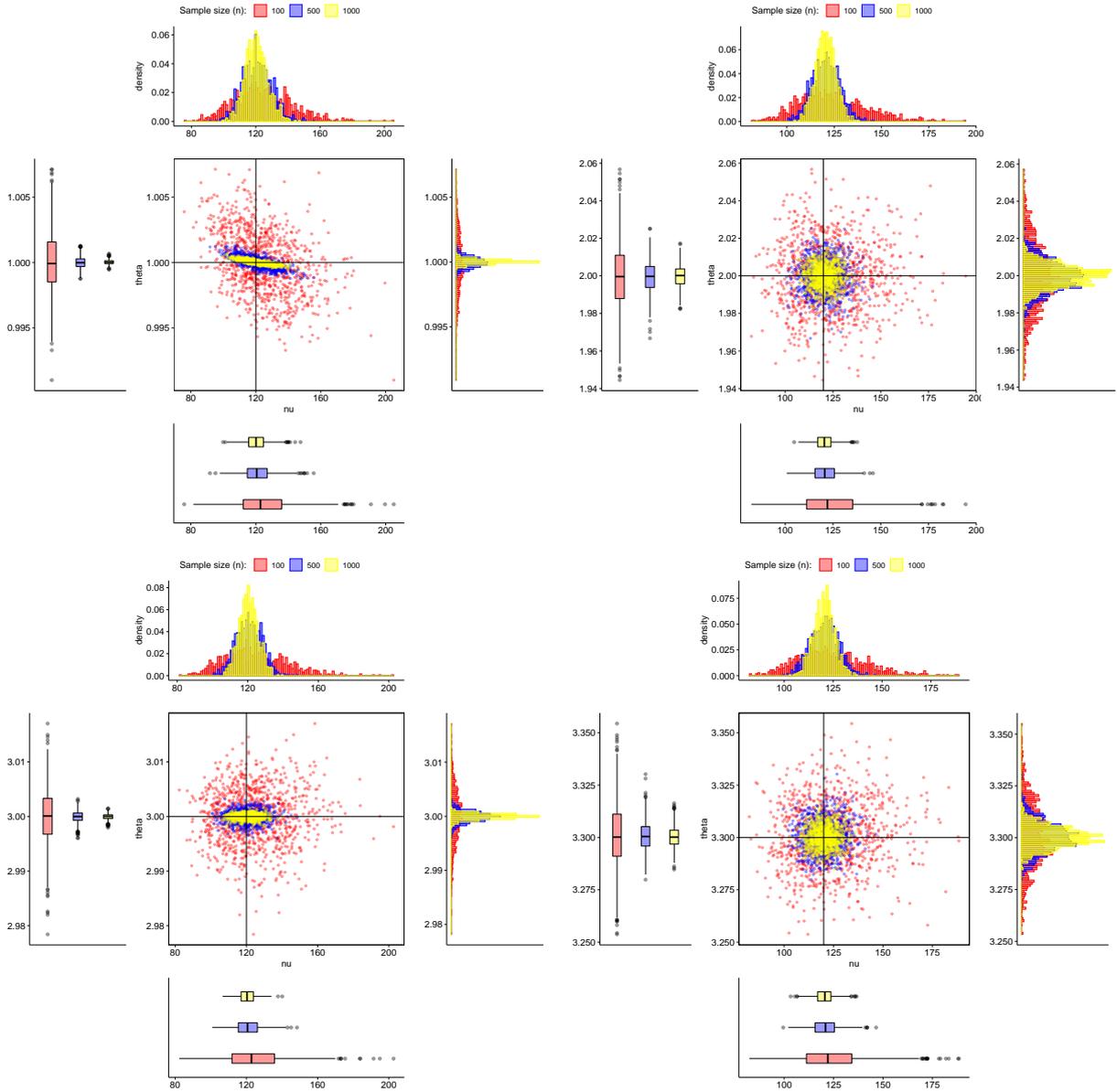

Figure 23: Histogram and Boxplots for $\nu = 120$ and $\theta \in \{1, 2, 3, 3.3\}$.

28

# 3 Simulation results for the $\beta$ARC model with covariates.

In this section we present a Monte Carlo simulation study considering a $\beta$ARC(1) model with a covariate. More specifically, we consider the $\beta$ARC(1) model with systematic component given by

$$g(\mu_t) = \alpha + \beta x_t + \phi\big[g(y_{t-1}) - \beta x_{t-1}\big] + T_\theta^{t-1}(u_0), \quad t = 1, \cdots, 3000, \qquad (3)$$

where $T_\theta$ is the logistic map, $g$ is the logit function, $\alpha = 0.6$, $\beta = 0.5$, $\phi = 0.2$, $\theta = 3.5$, $\nu = 20$, $u_0 = \pi/4$ and

$$x_t = \cos\big(2\pi t/365\big).$$

Figure 24 presents a sample from the $\beta$ARC(1) model considered in (3) along with the covariate, the orbit of $T_\theta$ and $\mu_t$.

For estimation purposes, we initialize the parameter as follows: $\alpha$ and $\beta$ are initialized by performing a least square fit with $g(y_t)$ as response and $x_t$ as the covariate, $\phi = 0.1$, $\nu = 50$ and $\theta = 3.65$. In Table 5 we present descriptive statistics based on 10,000 replications of the $\beta$ARC(1) model with systematic component (3). Presented are the mean, standard deviation, minimum and maximum values, the first, second (median) and third quartiles ($Q_1$, $Q_2$ and $Q_3$, respectively), MSE and MAPE of each parameter. Also presented are the descriptive statistics related to the number of function evaluations (neval) performed during optimization of the partial likelihood. We only consider the cases where the optimization achieved convergence (97.18%). In Figure 3, we show (clockwise from the top left) the scatter plots of $\hat\beta \times \hat\alpha$, $\hat\beta \times \hat\phi$, $\hat\theta \times \hat\phi$ and $\hat\theta \times \hat\alpha$ and along with histograms and boxplots.

Table 5: Descriptive statistics for the simulation results based on 10,000 replications of the $\beta$ARC(1) model with specification (3).

|       | $\alpha$ | $\beta$ | $\phi$ | $\theta$ | $\nu$ | neval |
|-------|--------|--------|--------|--------|---------|-------|
| True  | 0.6    | 0.5    | 0.2    | 3.5    | 20      | -     |
| Mean  | 0.6057 | 0.5006 | 0.1964 | 3.4994 | 20.0194 | 456   |
| sd    | 0.0328 | 0.0195 | 0.0184 | 0.0412 | 0.5370  | 102   |
| min   | 0.4709 | 0.4271 | 0.0904 | 3.2559 | 15.7518 | 224   |
| $Q_1$ | 0.5845 | 0.4879 | 0.187  | 3.4839 | 19.6656 | 393   |
| $Q_2$ | 0.6034 | 0.5003 | 0.1977 | 3.4991 | 20.0159 | 433   |
| $Q_3$ | 0.6232 | 0.5128 | 0.2082 | 3.5192 | 20.3675 | 486   |
| max   | 0.8135 | 0.6963 | 0.3135 | 3.5916 | 22.3958 | 1215  |
| MSE   | 0.0011 | 0.0004 | 0.0004 | 0.0017 | 0.2887  | -     |
| MAPE  | 0.0414 | 0.0299 | 0.0689 | 0.0077 | 0.0211  | -     |

We observe that the estimation performance is very good with small bias and standard deviation. For $\alpha$, $\beta$ and $\phi$ the histogram resembles a normal distribution as well as the associated scatter plots. The estimation of $\theta$ in $(3.4, 3.7)$ range typically yields a very flat likelihood and several local maxima, which makes identification often difficult. To exemplify this issue, Figure 26 present the plots of the likelihood as a function of $\theta$ in the interval $[3.2, 3.6]$ for 3 different samples. The blue line marks the (approximate) global maxima for the function (based on a grid of 200 points), the green line represents the true parameter ($\theta = 3.5$) and the red line marks the estimated value of $\theta$ based on that sample. In the first plot we see that the Nelder-Mead converged to a local maxima ($\hat\theta = 3.3560$). In the middle plot, all lines are close together and this is the typical case ($\hat\theta = 3.4846$). In the plot on the right, the algorithm did converge to the global maxima ($\hat\theta = 3.5703$), which is not located at the true parameter. Finally, Figure 26 reveals that the likelihood function may present several local maxima. Hence, it might be a good practice to run the optimization algorithm with several different starting points.

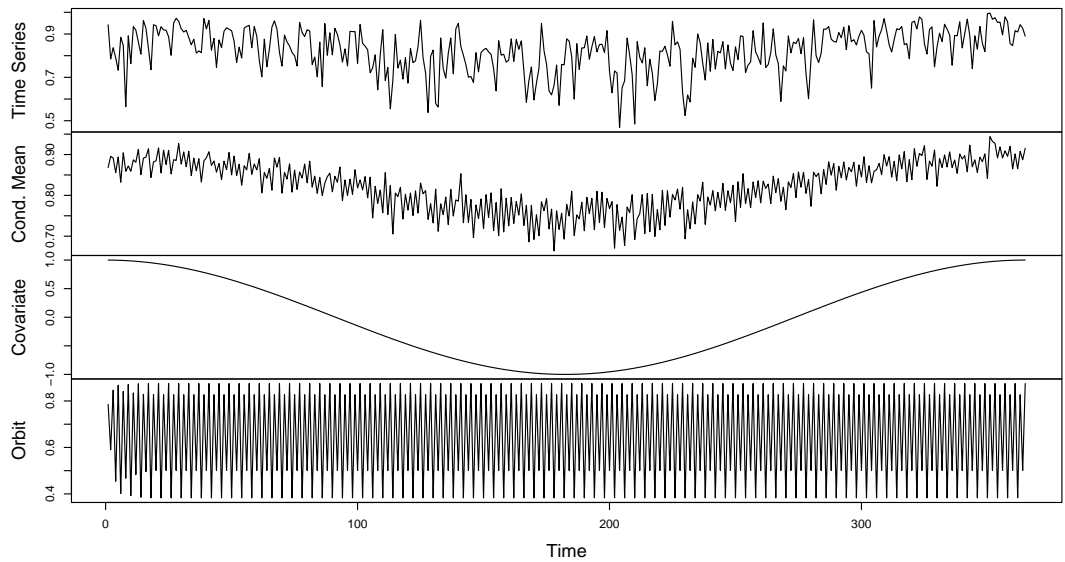

Figure 24: A simulated sample ($n = 365$) from the $\beta$ARC(1) model specified by (3) showing (from top to bottom) the generated time series, the conditional mean $\mu_t$, the effect of the covariate and the orbit of the transformation.

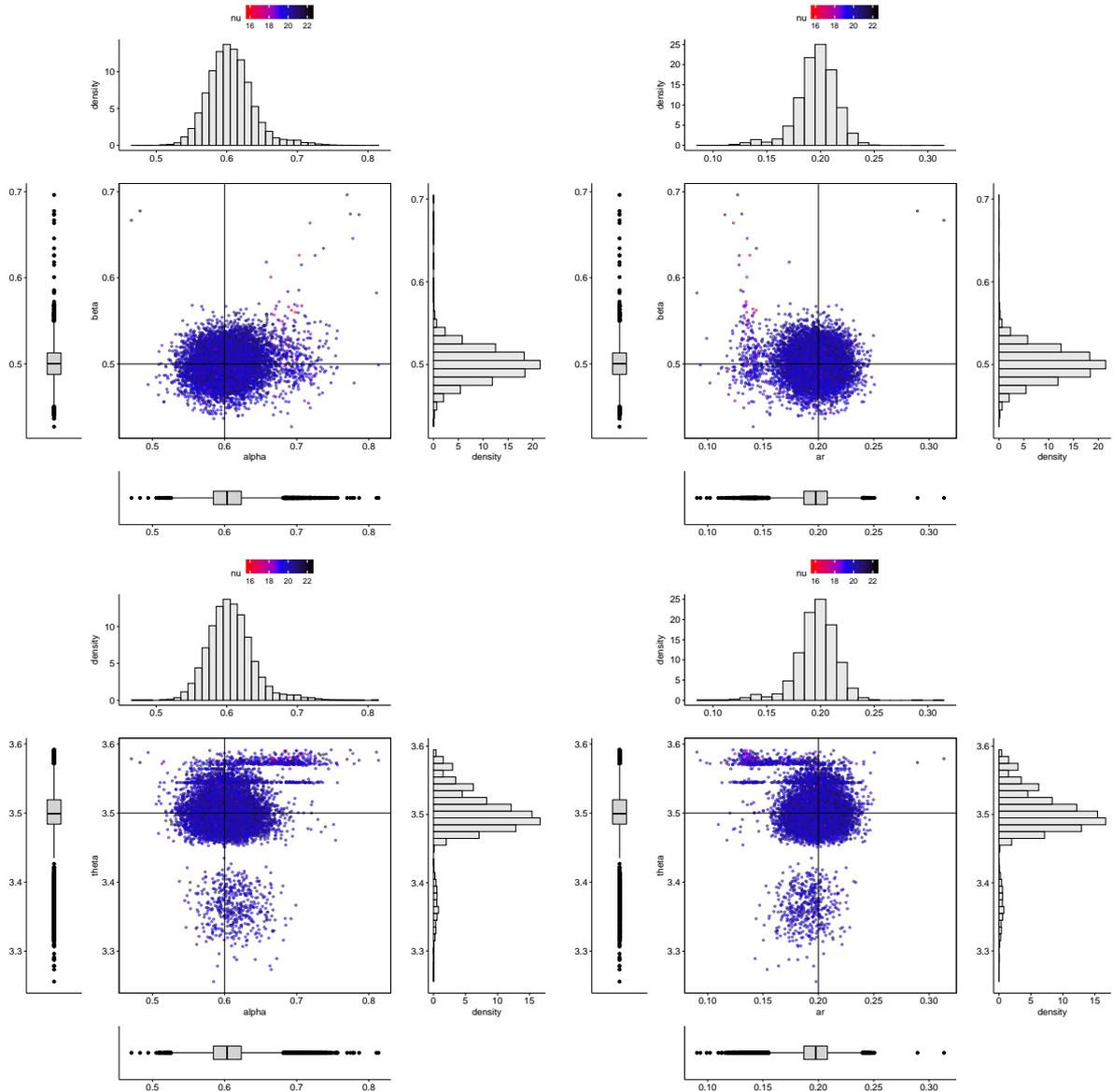

Figure 25: Clockwise from the top left, we present the scatter plots of $\hat{\beta} \times \hat{\alpha}$, $\hat{\beta} \times \hat{\phi}$, $\hat{\theta} \times \hat{\phi}$ and $\hat{\theta} \times \hat{\alpha}$ along with respective histograms and boxplots. The color applied in the scatter plots represent the estimated values of $\nu$, following the scale on the top of each plot.

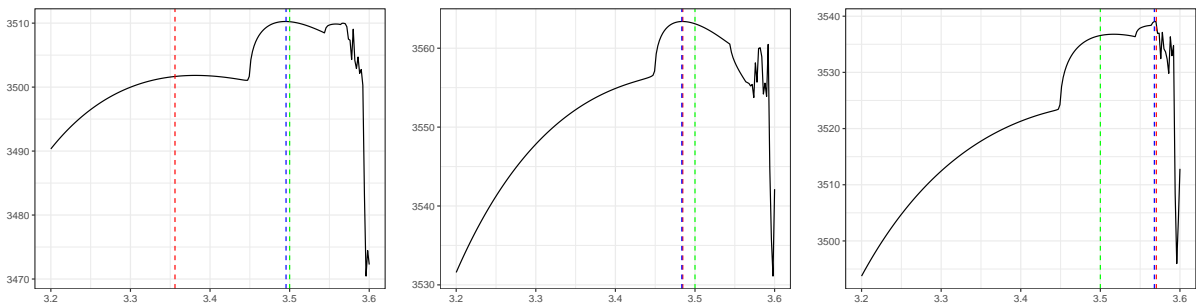

Figure 26: Plots of the likelihood as a function of $\theta$ in the interval $[3.2, 3.6]$ for 3 different samples, showing global maxima for the function (in blue), the true parameter $\theta = 3.5$ (green) and the estimated value of $\theta$ based on that sample (in red).

31



%
%
%

\end{document}